\documentclass[
a4paper,
11pt, reqno]{amsart}
\usepackage{amsmath,amssymb,amscd,amsthm,latexsym}

\usepackage{hyperref}
\usepackage{color, xcolor}
\usepackage{mathtools}
\mathtoolsset{showonlyrefs}

\usepackage{ mathrsfs }

\allowdisplaybreaks

\def\scK{{\mathscr K}}
\def\scG{{\mathscr G}}
\def\scF{{\mathscr F}}

\def\frL{{\mathfrak{L}}}
\def\frF{{\mathfrak{F}}}

\def\1{{\bf 1}}

\def\nn{\nonumber}

\newcommand{\bk}{\color{black}}

 \def\sE {{\cal E}} 
 \def\sF {{\cal F}}
\def\sG {{\cal G}}  
\def\sJ {{\cal J}}

 \def\sN {{\cal N}} 
  \def\sR {{\cal R}} 
\def\sS {{\cal S}}
\def\sT {{\cal T}}

 \def\bE {{\mathbb E}}

 \def\bN {{\mathbb N}} 
\def\bP {{\mathbb P}} 
 \def\bR {{\mathbb R}}

 \def\bZ {{\mathbb Z}}
\def\R {{\mathbb R}}  
 \def\Z {{\mathbb Z}}

\newtheorem{thm}{Theorem}[section]
\newtheorem{lemma}[thm]{Lemma}
\newtheorem{defn}[thm]{Definition}
\newtheorem{definition}[thm]{Definition}
\newtheorem{prop}[thm]{Proposition}
\newtheorem{corollary}[thm]{Corollary}
\newtheorem{remark}[thm]{Remark}
\newtheorem{example}[thm]{Example}
\numberwithin{equation}{section}

\def\qed{{\hfill $\Box$ \bigskip}}
\newcommand{\cal}[1]{\mathcal{#1}}
\def\NN{{\mathcal N}}

\def\FF{{\mathcal F}}
\def\EE{{\mathcal E}}

\def\R{{\mathbb R}}

\def\P{{\mathbb P}}
\def\N{{\mathbb N}}
\def\eps{\varepsilon}

\def\wt{\widetilde}
\def\pf{\noindent{\bf Proof.} }

\def\SCSJ{\mathrm{SCSJ}}
\def\Diff{\mathrm{Diff}}
\def\CSJ{\mathrm{CSJ}}
\def\HK{\mathrm{HK}}
\def\SHK{\mathrm{SHK}}
\def\GHK{\mathrm{GHK}}
\def\UHK{\mathrm{UHK}}
\def\UHKD{\mathrm{UHKD}}

\def\FK{\mathrm{FK}}

\def\E{\mathrm{E}}
\def\J{\mathrm{J}}
\def\Ch{\mathrm{Ch}}
\def\VD{\mathrm{VD}}
\def\RVD{\mathrm{RVD}}
\def\PI{\mathrm{PI}}
\def\Gcap{\mathrm{Gcap}}

\begin{document}
\allowdisplaybreaks

\title
[Heat kernel estimates and their stabilities for symmetric jump processes]{ Heat kernel estimates and their stabilities for symmetric jump processes with general mixed polynomial growths on metric measure spaces
}

\author{  Joohak Bae,  \quad Jaehoon Kang,   \quad Panki Kim  \quad and \quad   Jaehun Lee}
\address[Bae]{Samsung Fire \& Marine Insurance, Seoul 06620, Republic of Korea}\thanks{This research is  supported by the National Research Foundation of Korea(NRF) grant funded by the Korea government(MSIP) (No. 2016R1E1A1A01941893)}
\curraddr{}
\email{juhak88@snu.ac.kr}

\address[Kang]{Department of Mathematical Sciences,
	Seoul National University,
	Building 27, 1 Gwanak-ro, Gwanak-gu
	Seoul 08826, Republic of Korea}\thanks{This research is  supported by the National Research Foundation of Korea(NRF) grant funded by the Korea government(MSIP) (No. 2016R1E1A1A01941893).
}
\curraddr{}
\email{jhnkang@snu.ac.kr}

\address[Kim]{Department of Mathematical Sciences and Research Institute of Mathematics,
	Seoul National University,
	Building 27, 1 Gwanak-ro, Gwanak-gu
	Seoul 08826, Republic of Korea}\thanks{This research is  supported by the National Research Foundation of Korea(NRF) grant funded by the Korea government(MSIP) (No. 2016R1E1A1A01941893).
}
\curraddr{}
\email{pkim@snu.ac.kr}

\address[Lee]{Korea Institute for Advanced Study, 
85 Hoegiro, Dongdaemun-gu
Seoul 02455
Republic of Korea}\thanks{This research is supported by
 the National Research Foundation of Korea (NRF) grant funded by the Korea government (MSIP) : NRF-2016K2A9A2A13003815.
}
\curraddr{}
\email{hun618@kias.re.kr}

\date{}

\begin{abstract}
In this paper, we consider a symmetric pure jump Markov process $X$ on a metric measure space with volume doubling conditions. Our focus is on estimating the transition density $p(t,x,y)$ of $X$ and studying its stability when the jumping kernel exhibits general mixed polynomial growth.

Unlike previous work, in our setting, the rate function governing the jump growth 
may not be comparable to 
 the scale function that determines whether $p(t,x,y)$ has near-diagonal or off-diagonal estimates. Under the assumption that lower scaling index of scale function is greater than $1$, we establish stabilities of heat kernel estimates. Additionally, if the metric measure space admits a conservative diffusion process with a transition density satisfying sub-Gaussian bounds, we generalize heat kernel estimates from \cite[Theorems 1.2 and 1.4]{BKKL} using the rate function and  the function $F$ related to walk dimension of underlying space.
 As an application, we prove the equivalence between a finite moment condition based on $F$ and a generalized Khintchine-type law of iterated logarithm at infinity for symmetric Markov processes.
	\end{abstract}
	
\maketitle	
	\bigskip
	


\allowdisplaybreaks

\begin{doublespace}
	
	\tableofcontents 
	
		\section{Introduction}

The study of heat kernel estimates for Markov processes has a long history, with numerous significant results in this field (see \cite{A, B, BB, BB04, BBK06, BL02, CHXZ, CKS, CK03, CK08, HKE, Gr94, GT, HK99, LY, Mu18, Sa92} and references therein). The intimate interplay between symmetric Markov processes and positive self-adjoint operators highlights the fundamental role of the heat kernel in connecting Probability theory and Partial differential equations. Heat kernel estimates for Markov processes on metric measure spaces not only provide insights into the behavior of these processes but also reveal intrinsic properties such as the walk dimension of the underlying space (\cite{B, BB, GT, LY}).

Recently, there has been significant research interest in studying heat kernel estimates for Markov processes with jumps, owing to their significance in both theory and applications (see \cite{BBCK, BGK, BKKL, BL02, BGR14, Chen, CKK, CK03, CK08, CZ16, HKE, PHI, CZ17, CZ18, GHH17, GHH18, GHL14, GRT18, GS1, GS2, HL, J, KaSz15, KaSz17, KL, KM, KS12, KR, Mal23, M, MS1, MS2, Sz11, Sz17} and references therein). 
The pioneering work in \cite{CK08} introduced heat kernel estimates for Markov processes with jumping kernels that satisfy mixed polynomial growth conditions. 
In \cite{CK08}, the authors considered pure jump symmetric Markov processes on a metric measure space $(M, d, \mu)$, where $M$ is a locally compact separable metric measure space equipped with the metric $d$ and a positive Radon measure $\mu$ satisfying the volume doubling property. The jumping kernel $J(x,y)$ of the Markov process in \cite{CK08} satisfies the following conditions:
	\begin{equation}\label{e:bm}
	\frac{c^{-1}}{V(x,d(x,y)) \phi (d(x,y))} \le J(x, y) \le  \frac{c}{V(x,d(x,y)) \phi (d(x,y))}\ , \quad x,y \in M,
	\end{equation}
	where $V(x,r)=\mu(B(x,r))$ for all $x\in M$ and $r>0$, and $\phi$ is a strictly increasing function on $[0, \infty)$ satisfying  
	\begin{equation}\label{e:bm0}
	c_1 (R/r)^{\beta_1} \le  \phi(R)/\phi(r) \le c_2 (R/r)^{\beta_2}, \quad 0<r<R< \infty
	\end{equation}  with $0<\beta_1\le \beta_2<2.$ Here, we say $\phi$ is {\it the rate function} since $\phi$ {describes the growth} of jumps. Under the assumptions \eqref{e:bm}, \eqref{e:bm0}, reverse volume doubling  property 
	 and $V(x,r)\asymp \wt V(r)$ for a strictly increasing function $\wt V$, 
	 in \cite{CK08}, it is shown that 
	 the transition density $p(t,x,y)$ of Markov process
	has the following estimates: 
	\begin{align}
	\begin{split}\label{stable}
	p(t,x,y)\asymp \left(\frac{1}{V(x,\phi^{-1}(t))}\wedge
	\frac{t}{V(x,d(x,y))\phi (d(x,y))}\right), \quad \text{for any }t>0 \text{ and } x,y \in M.
	\end{split}
	\end{align}
	(See \cite[Theorem 1.2]{CK08}. See also \cite{HKE}  where the extra condition $V(x,r)\asymp \wt V(r)$ is removed). 
	Here and below, we denote $a \wedge b:=\min\{a, b\}$ and $f\asymp g$ if the quotient $f/g$ remains bounded between two positive constants.
	
	We introduce the concept of a function $\Phi$ as "the scale function" for $p(t,x,y)$. Specifically, $\Phi(d(x,y)) = t$ serves as the threshold that determines whether $p(t,x,y)$ exhibits near-diagonal or off-diagonal estimates. 
Notably, in \eqref{stable}, $\phi$ represents the scale function, and it coincides with the rate function when $0<\beta_1\le\beta_2<2$.
	Furthermore, it can be observed that \eqref{e:bm} and \eqref{stable} are equivalent since $p(t,x,y)/t\to J(x,y)$ weakly as $t\to0$. Consequently, for a wide range of pure jump symmetric Markov processes on metric measure spaces that satisfy volume doubling properties, \eqref{e:bm} is equivalent to \eqref{stable} under the condition \eqref{e:bm0} with $0<\beta_1\le\beta_2<2$.
	\bk

	One of major problems in this field is to derive heat kernel estimates for jump processes on metric measure space without the restriction $\beta_2<2$. Recently, several articles
	have discussed this problem (\cite{BKKL, HKE, GHH18, HL, M, MS2}).  In \cite{HKE}, the authors established stability results for heat kernel estimates of the form \eqref{stable} for symmetric jump Markov processes on metric measure spaces that satisfy volume doubling and reverse volume doubling properties.  Results in \cite{HKE} cover metric measure spaces whose walk dimension is bigger than $2$ such as  Sierpinski gasket and Sierpinski carpet. 
	Without the condition $\beta_2<2$, \eqref{e:bm} is generally not equivalent to \eqref{stable}. In this case, it is shown in \cite[Theorem 1.13]{HKE} that 
	\eqref{stable}  is equivalent to
	 the conjunction of \eqref{e:bm} and a cut-off Sobolev type inequality ($\CSJ(\phi)$ in \cite[Definition 1.5]{HKE}).  It is worth noting  that $\CSJ(\phi)$ always holds if $\beta_2<2$ (see \cite[Remark 1.7]{HKE}). As a corollary of main results in   \cite{HKE}, it is also shown in \cite{HKE} that \eqref{e:bm} is equivalent to \eqref{stable} if underlying space 
	 has a walk dimension strictly bigger than $\beta_2$.
		In \cite{GHH18}, a similar equivalence relation was proven for $\phi(r)=r^{\beta}$. Particularly, it is shown in \cite{GHH18} that the conjunction of \eqref{e:bm} and a generalized capacity condition is also equivalent to \eqref{stable}.

	On the other hand, in \cite{BKKL, M},  new forms of  heat kernel estimates for symmetric jump Markov processes in $n$-dimensional  Euclidean spaces were obtained without the condition $\beta_2<2$.  In particular, the outcomes in \cite{BKKL, M} cover Markov processes with a high intensity of small jumps. In such cases, the rate function and the scale function may not be comparable, and the heat kernel estimates are formulated in a more general manner, which includes \eqref{stable} with $V(x,d(x,y))=|x-y|^n$ (see \cite[Theorem 1.2 and Theorem 1.4]{BKKL} and \cite[Theorem 1.2]{M}).

	This paper is a continuation of our journey  on investigating the estimates of the transition densities of jump processes whose jumping kernels have general mixed polynomial growths. Our focus is on a symmetric pure jump Markov process $X$  on a general metric measure space $M$ that satisfies both volume doubling and reverse volume doubling properties.
The objective of this paper is twofold.

	Firstly, we aim to establish various versions of heat kernel estimates for the symmetric pure jump Markov process $X$, where the jumping kernel satisfies mixed polynomial growths, specifically conditions \eqref{e:bm} and \eqref{e:bm0} with $0<\beta_1\le \beta_2<\infty$. It is important to note that in our framework, the rate function and the scale function may not be comparable. This extension will build upon the key findings in \cite{BKKL}.
	
	Secondly, we intend to explore conditions that are equivalent to our heat kernel estimates. This investigation will further expand upon the principal outcomes in \cite{HKE}, where the scale function is comparable to the rate function.

  In this paper, we will address various types of heat kernel estimates, each requiring different assumptions to derive our results.
  To begin, we establish an upper bound on the heat kernel and its stability, which extends \cite[Theorem 1.15]{HKE} and \cite[Theorem 4.5]{BKKL}. This upper bound is obtained under the assumption that the lower scaling index $\delta$ of the scale function is strictly greater than $1$. It should be noted that in the general case, where $M$ may not satisfy the chain condition, the upper and lower bounds in the generalized version of \cite[Theorem 1.4]{BKKL} may take different forms. In order to achieve precise two-sided estimates, we additionally assume that the metric measure space satisfies the chain condition.
  Under the chain condition and with $\delta > 1$, our results in Theorems \ref{t:main10}, \ref{t:main21}, and \ref{t:main12} establish sharp heat kernel estimates and their stability, which generalize  \cite[Theorem 1.13]{HKE} and \cite[Theorem 1.4]{BKKL}.
 
 To extend the heat kernel estimates from \cite[Theorem 1.2]{BKKL} and establish the corresponding stability result, we introduce the assumption that the underlying space admits a conservative diffusion process with a transition density satisfying general sub-Gaussian bounds in terms of an increasing function $F$ (as defined in Definition \ref{d:diffusion}). This function $F$ serves as a generalization of the walk dimension for the underlying space. Notably, in \cite[Theorem 1.2]{BKKL}, the term $(d(x,y)/\Phi^{-1}(t))^2$ appears in the exponential term of the off-diagonal part, where the exponent 2 represents the walk dimension of Euclidean space.
 
 It has been established in \cite{GT} that the general sub-Gaussian bounds for diffusion processes are equivalent to the conjunction of elliptic Harnack inequality and mean exit time estimates for the diffusion process, assuming the volume doubling property a priori. Examples of spaces satisfying our assumption include diffusion processes on the Sierpinski gasket and generalized Sierpinski carpets (\cite{BB, BP88}). Studies on stability of (sub-)Gaussian type heat kernel estimates for diffusion processes on metric measure spaces can be found in \cite{BB04, BBK06, Gr94, Mu18, Sa92}.
 
 Under the assumption of general sub-Gaussian bounds for the diffusion process with the function $F$, we can explicitly define the scale function using the rate function and $F$ (see \eqref{e:dPhi}). In Theorem \ref{t:main32}, we demonstrate that general weak scaling assumptions on the jumping kernel (rate function) are equivalent to the two-sided estimates of $p(t,x,y)$ expressed in terms of the given rate function and $F$. These estimates are sharp when the lower scaling index of the scale function is strictly greater than $1$ and the underlying metric space satisfies the chain condition.

We further extend the result from \cite[Theorem 5.2]{BKKL} to encompass Markov processes on general metric measure spaces. In particular, we establish that if the walk dimension $\gamma$ of the underlying space is strictly greater than $1$, then the $\gamma$-th moment condition for the Markov process is equivalent to the finiteness of $\limsup_{t \to \infty}d(x, X_t)(\log\log t)^{1-1/\gamma}t^{1/\gamma}$. The full version of this result can be found in Theorem \ref{t:lil}. This extension allows us to connect the moment condition for the Markov process with the asymptotic behavior of the distance between the process and a fixed point over time, incorporating the logarithmic correction factor.
\bk
	
		\vspace{3mm}
	
	\noindent
	{\it Notations} : Throughout this paper, the positive constants $a$, $A$, $c_F$, $c_\mu$, $C_1$, $C_2$, $C_\mu$, $C_L$, $C_U$, $\wt C_L$, $\beta_1$, $\beta_2$, $\delta$, $\gamma_1$, $\gamma_2$, $\eta$, $d_1$, $d_2$ will remain the same, whereas  $C$, $c$, and  $c_0$, $a_0$, $c_1$, $a_1$, $c_2$, $a_2$, $\ldots$  represent positive
	constants having insignificant values that may be changed  from one
	appearance to another. The constants $\alpha_1$, $\alpha_2$, $c_L$, $c_U$ remain the same until Section \ref{s:ODLE}, and redefined in Section \ref{s:main3}. All these constants are positive and finite.
	The labeling of the constants  $c_0$, $c_1$, $c_2$, $\ldots$ begins anew in the proof of
	each result.  $c_i=c_i(a,b,c,\ldots)$, $i=0,1,2,  \dots$, denote generic constants depending on $a, b, c, \ldots$. The constant $\bar C$ in \eqref{e:J_psi}
	may not
	be explicitly mentioned.
	We will use ``$:=$" to denote a
definition, which is read as ``is defined to be".
We use notations $a \wedge b := \min \{ a, b\}$, 
	$a \vee b := \max \{ a, b\}$, $\bR_+:=\{r\in \bR:r>0\}$, and $B(x,r):=\{y\in M: d(x,y)<r\}$. We say $f\asymp g$ if the quotient $f/g$ remains bounded between two positive constants. Also, let  $\lceil a \rceil := \sup \{n \in \Z : n \le a  \}$.

\section{Settings and Main results}
\subsection{Settings}
Let $(M, d)$ be a locally compact separable metric space, and $\mu$ be a positive Radon measure on $M$ with full support and $\mu(M)=\infty$. We also assume that every ball in $(M,d)$ is relatively compact.

For $x\in M$ and $r>0$, define $V(x,r):= \mu(B(x,r))$ be the measure of an open ball $B(x,r)$.
\begin{definition}\label{d:VD}
{\rm	(i) We say that the metric measure space $(M,d,\mu)$ satisfies the {\it volume doubling property} $\VD(d_2)$ with index $d_2>0$ if there exists a constant $C_\mu \ge 1$ such that
	\begin{equation*}\label{VD1}
	\frac{V(x,R)}{V(x,r)} \le  C_\mu \left( \frac{R}{r} \right)^{d_2} \quad \mbox{for all }x \in M \mbox{ and } 0<r \le R.
	\end{equation*}
	(ii) We say that $(M,d,\mu)$ satisfies the {\it reverse volume doubling property} $\RVD(d_1)$ with index $d_1>0$ if there exists a constant $c_\mu > 0$ such that 
	\begin{equation*}\label{RVD}
	\frac{V(x,R)}{V(x,r)} \ge c_\mu \left( \frac{R}{r} \right)^{d_1} \quad \mbox{for all } x \in M\mbox{ and }0<r \le R.
	\end{equation*}
}\end{definition}
Note that $V(x,r)>0$ for every $x \in M$ and $r>0$ since $\mu$ has full support on $M$. Also, under $\VD(d_2)$, we have  
\begin{equation}\label{VD2}
\frac{V(x,R)}{V(y,r)} \le \frac{V(y,d(x,y)+R)}{V(y,r)} \le C_\mu \left( \frac{d(x,y)+R}{r} \right)^{d_2} \quad \mbox{for all } x \in M\mbox{ and }0<r \le R.
\end{equation}

\begin{definition}\label{d:chain}
{\rm	We say that a metric space $(M,d)$ satisfies the \textit{chain condition} $\Ch(A)$ if there exists a constant $A \ge 1$ such that, for any $n \in \N$ and $x,y \in M$, there is a sequence $\{z_k\}_{k=0}^n$ of points in $M$ such that $z_0 =x, z_n=y$ and
	$$ d(z_{k-1},z_k) \le A \frac{d(x,y)}{n} \quad \mbox{for all} \quad k=1,\dots, n. $$}
\end{definition}

\begin{definition}\label{d:diffusion}
{\rm	For a strictly increasing function $F:(0,\infty) \to (0,\infty)$, we say that a metric measure space $(M,d,\mu)$ satisfies the condition $\Diff(F)$ if there exists a conservative symmetric diffusion process $Z = (Z_t)_{t \ge 0 }$ on $M$ such that the transition density $q(t,x,y)$ of $Z$ with respect to $\mu$ exists and it satisfies the following estimates: there exist constants  $c\ge1$ and $a_0>0$ such that for all $t>0$ and $x,y\in M$,
	\begin{align}\label{hkediff}
	\frac{c^{-1}}{V(x,F^{-1}(t))} {\bf 1}_{\{F(d(x,y)) \le t\}}\le q(t,x,y)\le \frac{c}{V(x,F^{-1}(t))}\exp\big(-a_0 F_1(d(x,y),t) \big),
	\end{align}
	where the function $F_1$ is defined as	
	\begin{align}\label{d:F1}
	F_1(r,t) : = \sup_{s>0} \left[ \frac{r}{s} - \frac{t}{F(s)} \right].
	\end{align}
}\end{definition}
The function $F_1$ in \eqref{d:F1} has been  already used in \cite{AB, GT}. 
Note that, when $F(r)=r^{\gamma}$ for some $\gamma >1$, then $F_1(r,t) \asymp (r^{\gamma}/t)^{1/(\gamma-1)}$.

We recall the following definitions from \cite{BKKL}.

\begin{defn}
	{\rm Let $g:(0,\infty) \to (0,\infty)$, and $a\in(0, \infty]$,  $\beta_1, \beta_2>0$, and $0< c\le1\le C$.
		\begin{enumerate}
			\item[(i)] For $a \in (0,\infty)$, we say that $g$ satisfies $L_a(\beta_1, c)$ (resp. $L^a(\beta_1, c)$) if
			$$ \frac{g(R)}{g(r)} \geq c \left(\frac{R}{r}\right)^{\beta_1} \quad \text{for all} \quad r\leq R< a\;(\text{resp.}\;a\le r\leq R).$$
			We also say that $L_a(\beta_1,c,g)$ (resp. $L^a(\beta_1,c,g)$) holds.
			\item[(ii)] We say that $g$ satisfies $U_a(\beta_2, C)$ (resp. $U^a(\beta_2, C)$) if
			$$ \frac{g(R)}{g(r)} \leq C \left(\frac{R}{r}\right)^{\beta_2} \quad \text{for all} \quad r\leq R< a\;(\text{resp.}\;a\le r\leq R).$$
			We also say that $U_a(\beta_2,C,g)$ (resp. $U^a(\beta_2,C,g)$) holds.
			\item[(iii)] 
			When $g$ satisfies $L_a(\beta_1, c)$ (resp. $U_a(\beta_2, C)$) with $a=\infty$, then we 
			 use short notations $L(\beta_1, c)$ (resp. $U(\beta_2, C)$) and we say that $g$ satisfies the global lower scaling condition with index  $\beta_1$. (resp. the global upper scaling condition with index  $\beta_2$.) 
 We also say that $L(\beta_1, c, g)$ (resp. $U(\beta_2, C,g)$) holds.
		\end{enumerate}
}\end{defn}

{\it Throughout this paper, we will assume that  $\psi:[0,\infty)\to[0,\infty)$ is a 
	non-decreasing function satisfying $L(\beta_1,C_L)$ and $U(\beta_2,C_U)$ for some $0<\beta_1 \le \beta_2$.}	Note that $L(\beta_1,C_L,\psi)$ implies $\lim_{t \to 0} \psi(t) = 0$.

 We assume that there exists a regular Dirichlet form $(\sE, \sF)$ on $L^2(M,\mu)$, which is given by
\begin{align}\label{e:DF}
\sE(u,v)=\int_{M\times M\setminus \{(x, x): x\in M\}}(u(x)-u(y))(v(x)-v(y)) J(dx,dy),\quad\;\;u, v\in \sF,
\end{align}
where $J(dx,dy)$ is a symmetric and Radon measure on $M \times M \setminus \{(x, x): x\in M\}$.
See \cite[page 6]{FOT} for the definition of the regularity of a Dirichlet form. 

\begin{definition}{\rm We say $\J_\psi$ holds if 
\begin{align}
\label{e:JAC}
J(dx,dy)=J(x,y)\mu(dx)\mu(dy), \quad \text{for
$\mu \times \mu$-almost all $x, y\in M$}
\end{align}
where $J$ is a symmetric and positive Borel measurable function on $M \times M \setminus \{(x, x): x\in M\}$, 
and 
there exists a constant $\bar C>1$ so that for every $x, y \in M$,
		\begin{equation}\label{e:J_psi}
		\frac{\bar C^{-1}}{V(x,d(x, y)) \psi (d(x, y))}\le J(x, y) \le \frac{\bar C}{V(x,d(x, y)) \psi (d(x, y))}.
		\end{equation}
		We say that $\J_{\psi,\le}$ (resp. $\J_{\psi,\ge}$) if \eqref{e:JAC} holds and the upper bound (resp. lower bound) in \eqref{e:J_psi} holds.}
\end{definition}

Associated with the regular Dirichlet form $(\EE, \FF)$ on $L^2(M;\mu)$ is a $\mu$-symmetric Hunt process $X= \{ X_t , t \ge 0; \P^x , x \in M \setminus \NN   \}$. Here $\NN$ is a properly exceptional set for $(\EE,\FF)$ in the sense that $\mu(\NN)=0$ and $\P^x (X_t \in \NN \mbox{ for some } t>0)=0$ for all $x \in M \setminus \NN$. This Hunt process is unique up to a properly exceptional set. (See \cite[Theorem 4.2.8]{FOT}.) 
$(\EE, \FF)$ has jump part only in terms of Beurling-deny formula in \cite[Theorem 3.2]{FOT} and so $X$ is a pure jump process.
We fix $X$ and $\NN$, and write $M_0 := M \setminus \NN.$

For a set
$U\subset M$ and process $X$, define the first exit time $\tau_U = \inf\{ t >0 : X_t \in U^c \}$.

\begin{definition}{\rm For a non-negative function $\phi$, we say that $\E_\phi$ holds if
		there is a constant $c \ge 1$ such that
		$$c^{-1}\phi(r)\le \bE^x [ \tau_{B(x,r)} ] \le c\phi(r) \quad \mbox{for all } x \in M_0, \,\, r>0.$$ We say that $\E_{\phi,\le}$ (resp. $\E_{\phi,\ge}$) holds
		if the upper bound (resp. lower bound) in the inequality above holds. }
\end{definition}

\begin{remark}\label{r:exit}
{\rm	Suppose  $\RVD(d_1)$, $\VD(d_2)$ and $\J_{\psi, \ge}$ hold.
	 Let $x \in M_0$ and $r>0$. By the L\'evy system in \cite[Lemma 7.1]{HKE} and $\J_{\psi, \ge}$, we have
	 that 
\begin{align*}
1 &\ge \P^x(X_{\tau_{B(x,r)}} \in B(x,2r)^c ) = \bE^x \left[\int_0^{\tau_{B(x,r)}} \int_{B(x,2r)^c} J(X_s,y) \mu(dy)ds \right] \\  &\ge
c_0 \bE^x[\tau_{B(x,r)}] \int_{B(x,2r)^c} \frac{1}{V(d(x,y))\psi(d(x,y))} \mu(dy).
\end{align*}
By $\RVD(d_1)$, there exists a constant $c_1>1$ such that $V(x,c_1r) \ge 2V(x,r)$ for any $x \in M$ and $r>0$.	Using this and $U(\beta_2,C_U,\psi)$ we obtain 
\begin{align*} \int_{B(x,2r)^c} \frac{1}{V(d(x,y))\psi(d(x,y))}\mu(dy) &\ge \int_{B(x,2c_1r) \setminus B(x,2r)}\frac{1}{V(d(x,y))\psi(d(x,y))}\mu(dy) \\ &\ge \frac{V(x,2c_1r) - V(x,2r)}{V(x,2c_1r)} \frac{1}{\psi(2c_1r)} \ge \frac{c_2}{\psi(r)}.
\end{align*}
Combining two estimates, we obtain 	$$ \bE^x[\tau_{B(x,r)}] \le c\psi(r), \quad x \in M_0, \,\, r>0,$$
	which implies 
	$\E_{\psi,\le}$.	
}\end{remark}

By Remark \ref{r:exit}, we expect that our scale function with respect to the process $X$, which should be comparable to the exit time $\bE^x[\tau_{B(x,r)}]$, is smaller than $\psi$.

Let $\Phi:(0,\infty)\to(0,\infty)$ be a non-decreasing function satisfying $L(\alpha_1,c_L)$ and $U(\alpha_2, c_U)$ with some $0<\alpha_1 \le \alpha_2$ and $c_L, c_U>0$ and 
\begin{align}\label{comp1}
\Phi(r)< \psi(r), \quad \mbox{for all} \quad r>0.
\end{align}
Since \eqref{comp1} can be relaxed to the condition $\Phi(r) \le c\psi(r)$, 
by the virtue of Remark \ref{r:exit}, the assumption \eqref{comp1} is quite natural for the scale function.

Recall that $\alpha_2$ is the global upper scaling index of $\Phi$. If $\Phi$ satisfies $L_a(\delta,\wt C_L)$, then we have $\alpha_2\ge \delta$. Indeed, if $\delta>\alpha_2$, then for any $0<r \le R < a$, we have
$$ \wt C_L \Big( \frac Rr \Big)^{\delta} \le \frac{\Phi(R)}{\Phi(r)} \le c_U \Big( \frac Rr \Big)^{\alpha_2}, $$
which is contradiction by letting $r \to 0$. 

	Also, we define a function $\Phi_1:(0,\infty)\times(0,\infty)\to\R$ by
\begin{align}\label{d:Phi1}
\Phi_1(r,t) : = \sup_{s>0} \left[ \frac{r}{s} - \frac{t}{\Phi(s)} \right].
\end{align}
(c.f.,\cite{GT}.) See Section \ref{s:tran} for various properties of $\Phi_1$.

For $a_0, t, r>0$ and $x\in M$, we define
\begin{align}\label{d:G}
\sG(a_0, t, x, r)= \sG_{\Phi,\psi}(a_0,t,x,r)&:= \frac{t}{V(x, r)\psi(r)} + \frac{1}{V(x,\Phi^{-1}(t))} \exp{\left(-a_0 \;\Phi_1(r, t) \right)},
\end{align} 
where $\Phi^{-1}$ is the generalized inverse function of $\Phi$, i.e., $\Phi^{-1}(t):= \inf\{s\geq0: \Phi(s)> t   \}$ (with the convention $\inf\emptyset=\infty$).

We say $p(t, x,y)$ is the heat kernel of the semigroup $\{P_t\}$ associated with $(\sE,\sF)$  if (it exists and) 
it is a non-negative symmetric measurable function on $M_0 \times M_0$ for every $t > 0$, 
$\bE^xf(X_t) = P_tf(x) = \int p(t, x, y)f(y) \mu(dy)$ and it satisfies Chapman-Kolmogorov equation.
See \cite[(1.2)--(1.4)]{HKE} for the precise definition. 

\begin{definition}\label{D:1.11}    \begin{itemize}
		\item[(i)] We say that $\HK(\Phi, \psi)$ holds if there exists a heat kernel $p(t, x,y)$
		of the semigroup $\{P_t\}$ associated with $(\sE,\sF)$,
		which has the
		following estimates: there exist $\eta, a_0>0$ and $c \ge1$ such that for all $t>0$ and $x,y\in M_0$,
		\begin{equation}\label{HKjum}
		\begin{split}
		&c^{-1}\left( \frac{1}{V(x, \Phi^{-1}(t))} {\bf 1}_{\{ d(x,y) \le \eta\Phi^{-1}(t)\}}+ \frac{t}{V(x, d(x,y)) \psi(d(x,y))}{\bf 1}_{\{ d(x,y) > \eta\Phi^{-1}(t)\}} \right) \\
		&\le p(t,x,y) \le  c \left( \frac{1}{V(x, \Phi^{-1}(t))} \wedge  \sG\big(a_0, t, x, d(x,y)\big) \right),
		\end{split}
		\end{equation}
		where $\sG = \sG_{\Phi,\psi}$ is the function in \eqref{d:G}.

		\item[(ii)] We say $\UHK(\Phi,\psi)$  holds if the upper bound in \eqref{HKjum} holds.
		
		\item[(iii)]  We say $\UHKD(\Phi)$ holds if there is a constant $c>0$ such that for all $t>0$ and $x\in M_0$,
		$$p(t, x,x)\le \frac c{V(x,\Phi^{-1}(t))}.$$
		
		\item[(iv)] We say that $\SHK(\Phi, \psi)$ holds if  there exist $a_L\ge a_U>0$ and $c \ge1$ such that for all $t>0$ and $x,y\in M_0$,
		\begin{equation*}
		\begin{split}
	\frac{1}{c} \left(\frac{1}{V(x, \Phi^{-1}(t))}\wedge \sG\big(a_L, t, x, d(x,y)\big) \hskip -0.2em \right) 
	\hskip -0.2em	\le p(t,x,y) \le  c \left( \frac{1}{V(x, \Phi^{-1}(t))}\wedge \sG\big(a_U, t, x, d(x,y)\big) \hskip -0.2em \right). 
		\end{split}
		\end{equation*}
		\item[(v)] For two functions $\Phi$ and $\psi$, we say that $\GHK(\Phi, \psi)$ holds if there exist $0<a_U$, $0<\eta$ and $c \ge 1$ such that for all $t>0$ and $x,y\in M_0$,
		\begin{align}
		&c^{-1}V(x, \Phi^{-1}(t))^{-1}{\bf 1}_{\{ d(x,y) \le \eta \Phi^{-1}(t)\}}+ \frac{c^{-1} t}{V(x, d(x,y)) \psi(d(x,y))}{\bf 1}_{\{ d(x,y) \ge \eta\Phi^{-1}(t)\}} \le   p(t, x,y) \nn \\
		& \le \frac{c}{V(x, \Phi^{-1}(t))}\wedge \left(\frac{c\,t}{V(x, d(x,y))\psi(d(x,y))} + \frac{c}{V(x, \Phi^{-1}(t))} e^{-a_U F_1( d(x,y), F(\Phi^{-1}(t))  ) }  \right). \label{GHKjum}
		\end{align}					
	\end{itemize}
\end{definition}

Note that, when $F(r)=r^{\gamma}$ for some $\gamma >1$, then $F_1(d(x,y) ,F(\Phi^{-1}(t))  )  \asymp \big(\frac{d(x,y)}{\Phi^{-1}(t)} \big)^{\gamma/(\gamma-1)}$.

\begin{remark}\label{r:hk}
	{\rm	{ For non-decreasing function} $\Phi : [0,\infty) \to [0,\infty)$ with $\Phi(0)=0$ satisfying $L(\alpha_1,c_L)$ and $U(\alpha_2,c_U)$ and for any $C>1$, the condition $\HK(\Phi, C\Phi)$ is equivalent to the existence of heat kernel $p(t,x,y)$ and the constant $c\ge1$ such that for all $t>0$ and $x,y \in M_0$,
		\begin{equation}\label{e:rhk1}
		\frac{c^{-1}}{V(x,\Phi^{-1}(t))} \land \frac{c^{-1}t}{V(x,d(x,y))\Phi(d(x,y))} \le p(t,x,y) \le \frac{c}{V(x,\Phi^{-1}(t))} \land \frac{ct}{V(x,d(x,y))\Phi(d(x,y))}.
		\end{equation}
		This shows that if $\Phi\asymp \psi$, then the condition $\HK(\Phi,\psi)$ is equivalent to \eqref{e:rhk1}, which is the condition $\HK(\Phi)$ of \cite{HKE}. The proof of the equivalence of \eqref{e:rhk1} and $\HK(\Phi,C\Phi)$ is in Appendix \ref{s:App}. }
\end{remark}

\begin{definition}\label{d:PI}
{\rm	We say that the (weak) Poincar\'e inequality $\PI(\Phi)$ holds if there exist constants $C>0$ and $\kappa \ge 1$ such that for any ball $B_r:=B(x,r)$ with $x \in M$, $r>0$ and for any bounded $f \in \FF$,
	{\begin{equation}\label{e:PI}
	\int_{B_r} (f- \bar{f}_{B_r})^2 d\mu \le C \Phi(r) \int_{B_{\kappa r} \times B_{\kappa r} } (f(y) - f(z))^2 \,J(dy,dz),
	\end{equation}}
	where $\bar{f}_{B_r} = \frac{1}{\mu(B_r)} \int_{B_r} f d\mu$ is the average value of $f$ on $B_r$.
	}\end{definition}

\begin{definition}{\rm
Let $U\subset M$ be an open set, $W$ be any Borel subset of $U$ and $\kappa\ge1$ be a real number. A $\kappa$-cutoff function of pair $(W, U)$ is any function $\varphi\in\sF$ such that $0\le\varphi\le\kappa$ $\mu$-a.e. in $M$,  $\varphi\ge1$ $\mu$-a.e. in $W$ and  $\varphi=0$ $\mu$-a.e. in $U^c$.
We denote by $\kappa$-cutoff$(W, U)$ the collection of all $\kappa$-cutoff function of pair $(W, U)$. Any $1$-cutoff function will be simply referred to as a cutoff function.
}\end{definition}
  Let $\sF':=\{u+b: u\in \sF, b\in\bR\}$.
\begin{definition}[{c.f. \cite[Definition 1.11]{GHH18}}]{\rm We say that $\Gcap(\Phi)$ holds if there exist constants $\kappa\ge1$ and $C>0$ such that for any bounded $u\in \sF'$ and for all  $x_0\in M$ and $R, r>0$, there exists a function $\varphi\in\kappa$-cutoff$(B(x_0, R), B(x_0, R+r))$ such that 
		$$\sE(u^2\varphi,\varphi) \le \frac{C}{\Phi(r)}\int_{B(x_0, R+r)}u^2d\mu.$$  }
\end{definition}
\subsection{Main results}

Recall that we always assume that  $\psi:[0,\infty)\to[0,\infty)$ is a 
	non-decreasing function satisfying $L(\beta_1,C_L)$ and $U(\beta_2,C_U)$ for some $0<\beta_1 \le \beta_2$.

For the function $\Phi$ satisfying \eqref{comp1} and $L^a(\delta,\wt C_L)$ with $\delta>1$, we define 
\begin{align}\label{d:wePhi}
\wt \Phi(s):=c_U^{-1} \frac{\Phi(a)}{a^{\alpha_2}}s^{\alpha_2}\1_{\{s<a\}}+\Phi(s)\1_{\{s\geq a\}}.
\end{align}
Note that for $s \le a$ we have, $\frac{\wt \Phi(s)}{\Phi(s)} = c_U^{-1} \frac{s^{\alpha_2}}{a^{\alpha_2}} \frac{\Phi(a)}{\Phi(s)} \le 1$. Thus,
\begin{equation}\label{e:wtcomp}
\wt \Phi(r) \le \Phi(r) < \psi(r), \quad r>0.
\end{equation}
Also, $L(\delta, \wt C_L,\wt \Phi)$ holds. Indeed, for any $0<r \le a \le R$,
$$ \frac{\wt \Phi(R)}{\wt \Phi(r)}= \frac{\wt \Phi(R)}{\wt \Phi(a)} \frac{\wt\Phi(a)}{\wt\Phi(r)} \ge \wt C_L \Big(\frac Ra \Big)^{\delta} \Big( \frac ar \Big)^{\delta} = \wt C_L \Big( \frac Rr \Big)^{\delta}.$$
The other cases are straightforward. 

By the same way as \eqref{d:Phi1}, let us define
\begin{equation}\label{e:wtPhi1} 
\wt \Phi_1(r,t):=  \sup_{s>0} \left[  \frac{r}{s} - \frac{t}{\wt\Phi(s)}\right].
\end{equation} 
The following Theorem \ref{t:main10}-Corollary \ref{c:main12} are the first set of main results of this paper.
{ Recall that  $(\EE,\FF)$ is the  regular Dirichlet form   on $L^2(M,\mu)$ defined \eqref{e:DF} and $X$ is the corresponding Hunt process.}
\begin{thm}\label{t:main10}
	Assume that the metric measure space $(M,d,\mu)$ satisfies $\VD(d_2)$, and the process $X$ satisfies  $\J_{\psi,\le}$, $\UHKD(\Phi)$ and $\E_{\Phi}$, where $\psi$ is a non-decreasing function satisfying $L(\beta_1,C_L)$ and $U(\beta_2,C_U)$, and $\Phi$ is a non-decreasing function satisfying \eqref{comp1}, $L(\alpha_1, c_L)$ and $U(\alpha_2, c_U)$, where $0<\beta_1 \le \beta_2$ and $0<\alpha_1 \le \alpha_2$.
	
	(i) Suppose that $\Phi$ satisfies $L_a(\delta, \wt C_L)$ with some $a > 0$ and $\delta>1$. Then, for any $T \in (0,\infty)$, there exist constants $a_U>0$ and $c>0$ such that for any $t<T$ and $x,y\in M_0$,
	\begin{equation}\label{e:main10}
	p(t,x,y) 
	\leq 
	 \frac{c}{V(x,\Phi^{-1}(t))} \land \left(
	\frac{c\, t}{V(x,d(x,y))\psi(d(x,y))} + \frac{c}{V(x,\Phi^{-1}(t))} \exp{\big(-a_U \Phi_1(d(x,y),t) \big)} \right). 
	\end{equation}
	Moreover, if $\Phi$ satisfies $L(\delta, \wt C_L)$, then \eqref{e:main10} holds for all $t<\infty$.
	
	(ii) Suppose that $\Phi$ satisfies $L^a(\delta, \wt C_L)$ with some $a > 0$ and $\delta>1$. Then, for any $T \in (0,\infty)$ there exist constants $a_U>0$ and $c>0$ such that for any $t\ge T$ and $x,y\in M_0$,
	\begin{equation}
	\label{e:main101}
	p(t,x,y) 
	\leq  \frac{c}{V(x,\Phi^{-1}(t))} \land \left( \frac{c\, t}{V(x,d(x,y))\psi(d(x,y))} + \frac{c}{V(x,\Phi^{-1}(t))}  \exp{\big(-a_U \wt{\Phi}_1(d(x,y),t) \big)} \right). 
	\end{equation}
\end{thm}

When the scale function $\Phi$ satisfies the global lower scaling condition whose index is great than $1$, 
by combining our Theorem \ref{t:main10} with \cite[Theorem 1.15]{HKE} and \cite{GHH18},
we obtain the following stability results.

\begin{thm} \label{t:main11}
	Assume that the metric measure space $(M, d, \mu)$ satisfies $\RVD(d_1)$ and $\VD(d_2)$. Let
	$\psi$ be a non-decreasing function satisfying $L(\beta_1,C_L)$ and $U(\beta_2,C_U)$, and $\Phi$ be a non-decreasing function satisfying \eqref{comp1}, $L(\alpha_1, c_L)$ and $U(\alpha_2, c_U)$ with $1<\alpha_1 \le \alpha_2$. Then the following are equivalent: \\
	$(1)$ $\UHK(\Phi, \psi)$ and $(\sE, \sF)$ is conservative. \\
	$(2)$ $\J_{\psi,\le}$, $\UHK(\Phi)$ and $(\sE, \sF)$ is conservative. \\
	$(3)$ $\J_{\psi,\le}$, $\UHKD(\Phi)$ and $\E_\Phi$.
\end{thm}

	See \cite[Definitions 1.5 and 1.8]{HKE} for the definitions of $\FK(\Phi)$, $\CSJ(\Phi)$ and $\SCSJ(\Phi)$.

\begin{corollary} \label{c:main11}
	Under the same settings as Theorem \ref{t:main11}, each equivalent condition in above theorem is also equivalent to the following: \\
	$(4)$ $\FK(\Phi)$, $\J_{\psi, \le}$ and $\SCSJ(\Phi)$. \\
	$(5)$ $\FK(\Phi)$, $\J_{\psi, \le}$ and $\CSJ(\Phi)$.\\
	$(6)$ $\FK(\Phi)$, $\J_{\psi, \le}$ and $\Gcap(\Phi)$.
\end{corollary}

When the Poincar\'e inequality $\PI(\Phi)$ holds and our metric measure space satisfies the reverse volume doubling property and the chain condition, we obtain the sharp two sided estimates of  $p(t,x,y)$.

\begin{thm}\label{t:main21}
	Assume that the metric measure space $(M,d,\mu)$ satisfies $\Ch(A)$, $\RVD(d_1)$ and $\VD(d_2)$. Suppose that the process $X$ satisfies $\J_\psi$, $\E_\Phi$ and $\PI(\Phi)$, where $\psi$ is a non-decreasing function satisfying $L(\beta_1,C_L)$ and $U(\beta_2,C_U)$, and $\Phi$ is a non-decreasing function satisfying \eqref{comp1}, $L(\alpha_1, c_L)$ and $U(\alpha_2, c_U)$.
	
	(i) Suppose that $L_a(\delta,\wt C_L,\Phi)$ holds with $\delta>1$. Then, for any $T \in (0,\infty)$, \eqref{e:main10} holds and there exist constants $c>0$ and $a_L>0$ such that for any $x,y \in M_0$ and $t \in (0,T]$,
	\begin{align}\label{e:main211}
	p(t,x,y) \geq \frac{c}{V(x,\Phi^{-1}(t))} \land \left(\frac{ct}{V(x,d(x,y))\psi(d(x,y))} + \frac{c}{V(x,\Phi^{-1}(t))}\exp \left(-a_L \Phi_1(d(x,y),t) \right)\hskip-1mm \right).
	\end{align}
	Moreover, if $L(\delta, \wt C_L,\Phi)$ holds, then \eqref{e:main211} holds for all $t \in (0,\infty)$.
	
	(ii) Suppose that $L^a(\delta,\wt C_L,\Phi)$ holds with $\delta>1$. Then, for any $T \in (0,\infty)$, \eqref{e:main101} holds and there exist constants $c>0$ and $a_L>0$ such that for any $x,y \in M_0$ and $t \ge T$,
	\begin{align}\label{e:main212}
	p(t,x,y) \geq \frac{c}{V(x,\Phi^{-1}(t))} \land \left(\frac{ct}{V(x,d(x,y))\psi(d(x,y))} + \frac{c}{V(x,\Phi^{-1}(t))}\exp \left(-a_L \wt \Phi_1(d(x,y),t) \right)\hskip-1mm \right).
	\end{align}
\end{thm}

\medskip

Recall that $\SHK(\Phi,\psi)$ is defined in Definition \ref{D:1.11} {\rm (vi)}.

\medskip

\begin{thm} \label{t:main12}
	Under the same settings as Theorem \ref{t:main11}, the following are equivalent: \\
	$(1)$ $\HK(\Phi, \psi)$. \\
	$(2)$  $\J_{\psi}$, $\PI(\Phi)$, $\UHK(\Phi)$ and $(\EE,\FF)$ is conservative. \\
	$(3)$ $\J_{\psi}$, $\PI(\Phi)$ and $\E_\Phi$. \\
	If we further assume that $(M,d)$ satisfies $\Ch(A)$ for some $A \ge 1$, then the following is also equivalent to others: \\
	$(4)$ $\SHK(\Phi,\psi)$.
\end{thm}
By Theorem \ref{t:main12} and  Corollary \ref{c:main11}, we also obtain that
\begin{corollary} \label{c:main12}
	Under the same settings as Theorem \ref{t:main11}, each equivalent condition $(1)$-$(3)$ in Theorem \ref{t:main12} is also equivalent to the following: \\
	$(5)$ $\J_{\psi}$, $\PI(\Phi)$ and $\SCSJ(\Phi)$.\\
	$(6)$ $\J_{\psi}$, $\PI(\Phi)$ and $\CSJ(\Phi)$.\\
	$(7)$ $\J_{\psi}$, $\PI(\Phi)$ and $\Gcap(\Phi)$.
\end{corollary}

From now on, we consider a metric measure space $(M,d,\mu)$ satisfying $\Diff(F)$. 
For the rest of this section, let $F$ be a strictly increasing function satisfying $L(\gamma_1,c_F^{-1})$ and $U(\gamma_2,c_F)$.

By the recent result \cite[Corollary 1.10]{Mu20}, we may and do assume that $2 \le \gamma_1 \le \gamma_2$. Indeed, from \cite[Theorem 1.2]{GHL15}, $\Diff(F)$ implies the conditions in \cite[Corollary 1.10]{Mu20}. Note that since $F_1(r,t) = tF_1(r/t,1)$ for any $r,t>0$, our $\Diff(F)$ and
${\rm (UE)}_{F} + {\rm (NLE)}_{F}$
in  \cite{GHL15} are the same.

{ Suppose that our  regular Dirichlet form $(\EE,\FF)$  on $L^2(M,\mu)$, defined \eqref{e:DF}, satisfies $\J_\psi$.}
Then, under $\Diff(F)$, we see by \cite[Theorems 2.2 and 2.3 and Lemma 3.5]{LM21} that the following condition holds:
\begin{align}\label{e:intcon}
\int_0^{1} \frac{dF(s)}{\psi(s)} \,<\infty.
\end{align}
Moreover, according to Lemma \ref{l:Y_psi}, \eqref{e:intcon} also implies the existence of a regular Dirichlet form satisfying $\J_\psi$.

Recall that the function $F_1(r,t) = \sup_{s>0} \big[ \frac{r}{s} - \frac{t}{F(s)} \big]$ has defined in \eqref{d:F1}. Since $L(\gamma_1,c_F^{-1})$ holds with $\gamma_1 \ge 2$, the function $F_1$ satisfies the properties in Subsection \ref{s:tran}. For instance, $F_1(r,t) \in (0,\infty)$ for all $r,t>0$. 
{By the above observation, under $\J_\psi$}, the function
\begin{equation}\label{e:dPhi} \Phi(r):= \frac{F(r)}{\int_0^r \frac{dF(s)}{\psi(s)}}, \quad r>0. \end{equation}
is well-defined. Also, $\Phi$ is a strictly increasing function satisfying \eqref{comp1}, 
 $U(\alpha_2, c_U)$
and $L(\alpha_1, c_L)$ for some $\alpha_2\ge\alpha_1>0$ and $c_U, c_L>0$ (see Section \ref{s:main3}).

We now ready to state the second  set of main results of this paper. Recall that $\GHK(\Phi,\psi)$ and $\SHK(\Phi,\psi)$ are defined in Definition \ref{D:1.11}{\rm (v)} and {\rm (iv)}.
\begin{thm} \label{t:main32}
	{Assume} that the metric measure space $(M,d,\mu)$ satisfies $\RVD(d_1)$, $\VD(d_2)$ and $\Diff(F)$, where $F:(0,\infty) \to (0,\infty)$ is a strictly increasing function satisfying $L(\gamma_1,c_F^{-1})$ and $U(\gamma_2,c_F)$. Let $\psi$ be a non-decreasing function satisfying $L(\beta_1,C_L)$ and $U(\beta_2,C_U)$. Then, the following are equivalent :

{
\noindent
$(1)$  $(\EE,\FF)$ satisfies $\J_\psi$.

 \noindent
$(2)$ $(\EE,\FF)$ satisfies $\GHK(\varphi,\psi)$ with a non-decreasing function $\varphi$ satisfying $U(\delta_2,c)$ for some $\delta_2>0$ and $c \ge 1$.

\noindent
 $(3)$  $(\EE,\FF)$ satisfies $\GHK(\Phi,\psi)$, where $\Phi$ is the function defined in \eqref{e:dPhi}.

When $\Phi$ in \eqref{e:dPhi} satisfies $L(\alpha_1,c_L)$ with $\alpha_1>1$ and $(M,d)$ satisfies $\Ch(A)$ for some $A \ge 1$, then the following is also equivalent to others:

\noindent
$(4)$  $(\EE,\FF)$ satisfies $\SHK(\Phi,\psi)$.}
	\end{thm}
	
{Note that $\GHK(\Phi,\psi)$ with $\Phi$ defined in \eqref{e:dPhi} implies that \eqref{e:intcon} holds. Indeed, if \eqref{e:intcon} does not hold, then $\Phi(r)=0$ for all $r>0$. This and $\GHK(\Phi,\psi)$ give $p(t,x,y)=0$ for any $t>0$ and $x,y\in M_0$, which is contradiction to the existence of our regular Dirichlet form $(\EE,\FF)$.} 

We say that $G(x,y) := \int_0^\infty p(t,x,y)dt$ be the Green function of $(\EE,\FF)$ (if it exists.) 
	
\begin{corollary}\label{c:main33}
	Under the same setting as Theorem \ref{t:main32}, we assume that a regular Dirichlet form $(\EE,\FF)$ satisfies $\J_\psi$. Then

	(i) $(\EE,\FF)$ satisfies $\PI(\Phi)$ and $\E_{\Phi}$, where $\Phi$ is the function in \eqref{e:dPhi}.

	(ii) Suppose that $\Phi$ satisfies $U(\alpha_2,c_U)$ with $\alpha_2<d_1$. Then, $G(x,y)\asymp \frac{\Phi(d(x,y))}{V(x, d(x,y))}$ for all $x,y \in M_0$.
\end{corollary}

It is worth mentioning that thanks to \cite{
Mu20}, in Theorem \ref{t:main32} we do not assume
 the lower scaling index of the scale function being  strictly bigger than $1$.
Note that, $\GHK(\Phi,\psi)$  in Theorem \ref{t:main32}(3) is not sharp in general.	Without the chain condition, even the transition density of the diffusion process may not have the sharp two-sided bounds. However, if the upper scaling index $\beta_2$ of the rate function is strictly less than the walk dimension, our heat kernel estimates $\GHK(\Phi,\psi)$ is equivalent to \eqref{stable} ($\phi$ replaced by $\Phi$). 
See  Remark \ref{r:hk}.

Finally, we now state local estimates of heat kernels.
\begin{corollary}\label{c:main32}
Under the same setting as Theorem \ref{t:main32}, we assume that 
$(M,d)$ satisfies $\Ch(A)$ for some $A \ge 1$ and that the regular Dirichlet form $(\EE,\FF)$ satisfies $\J_\psi$. Let $\Phi$ be the function defined in \eqref{e:dPhi}.

	(i) Assume that $L_a(\delta,\wt C_L,\Phi)$ holds with some $\delta>1$ and $a>0$. Then, for any $T \in (0,\infty)$, there exist constants $0<a_U \le a_L$ and $c>0$ such that \eqref{e:main10} and \eqref{e:main211} holds for all $t \in (0,T]$ and $x,y \in M_0$.

	(ii) Assume that $L^a(\delta,\wt C_L,\Phi)$ holds with some $\delta>1$ and $a>0$. Then, for any $T \in (0,\infty)$, there exist constants $0<a_U \le a_L$ and $c>0$ such that \eqref{e:main101} and \eqref{e:main212} holds for all $t \in [T,\infty)$ and $x,y \in M_0$.	
\end{corollary}

\begin{remark}\rm
After finishing the first version of this paper and submitting  to Arxiv, we  later found that \cite{CKW19} is announced. 
The first set of main results of this paper, Theorem \ref{t:main10}-Corollary \ref{c:main12}, is similar to the main results in \cite{CKW19}. See \cite[Remark 1.18]{CKW19}. The second set of main results of this paper on 
the heat kernel estimates on metric measure space under the assumption $\Diff(F)$, 
is not covered in \cite{CKW19}.

\end{remark}

\section{
HKE and stability on general metric measure space}\label{s:ODE}

\subsection{Implications of $\UHKD(\Phi)$, $\J_{\psi,\le}$ and $\E_\Phi$}\label{s:ODUE}

Recall that the function $\psi$ satisfies the conditions $L(\beta_1, C_L)$ and $U(\beta_2, C_U)$ with $0<\beta_1 \le \beta_2$. 
{\it Throughout this section, we will assume that $\Phi$ satisfies \eqref{comp1}, $L(\alpha_1,c_L)$ and $U(\alpha_2,c_U)$, with some $0 < \alpha_1 \le \alpha_2$. }

\noindent Fix $\rho>0$ and define  a bilinear form $(\sE^{\rho}, \sF)$ by
\begin{align*}
\sE^{\rho}(u,v)=\int_{M\times M \setminus \{(x, x): x\in M\}} (u(x)-u(y))(v(x)-v(y)){\bf 1}_{\{d(x,y)\le \rho\}}\, J(dx,dy).
\end{align*}
Clearly, the form $\sE^{\rho}(u,v)$ is well defined for $u,v\in \sF$, and $\sE^{\rho}(u,u)\le \sE(u,u)$ for all $u\in \sF$. Moreover, if $\J_{\psi,\le}$ holds, then for all $u\in \sF$,
\begin{equation*}\label{EErhocomp}\begin{split}
\sE(u,u)-\sE^{\rho}(u,u)&= \int(u(x)-u(y))^2{\bf 1}_{\{d(x,y)>\rho\}}\,J(x, y)\mu(dx)\mu(dy)\\
&\le 4\int_{M}u^2(x)\,\mu(dx)\int_{B(x,\rho)^c}J(x,y)\,\mu(dy)\le \frac{c_0\|u\|_{2}^2 }{\psi(\rho)},
\end{split}\end{equation*}
where the last inequality follows from \cite[Lemma 2.1]{HKE}.
Thus, $\sE_1 (u, u):= \sE(u,u) + \| u \|_2^2$ is equivalent to $\sE^{\rho}_1(u,u):=
\sE^{\rho}(u,u)+ \|u\|_2^2$ for every $u\in \sF$, which implies that  $(\sE^{\rho},
\sF)$ is also a regular Dirichlet form on $L^2(M ,d\mu)$.
We call $(\sE^{\rho}, \sF)$ the $\rho$-truncated
Dirichlet form. The Hunt process associated with $(\sE^{\rho}, \sF)$ which will be denoted by $X^{\rho}$  can be identified in distribution with the Hunt process of the original Dirichlet form $(\sE, \sF)$ by removing those jumps of size larger than  $\rho$. Let $p^{\rho}(t,x,y)$ be the transition density  of $X^{\rho}$ corresponds to $(\EE^{\rho}, \FF)$. For a set $U\subset M$, we define $\tau_U^{\rho}= \inf\{ t >0 : X^{\rho}_t \in U^c \}$.

For any  open set ${D} \subset M$, $\sF_D$ is defined to be the $\sE_1$-closure in $\sF$ of $\sF\cap C_c(D)$. Let 
$\{P_t^{D}\}$ and $\{P_t^{\rho,{D}}\}$ be the semigroups of $(\sE, \sF_{D})$ and $(\sE^{\rho}, \sF_{D})$, respectively. 

In the next two lemmas, we obtain a priori estimate for the upper bound of heat kernel.
The first one gives a sharper bound for \cite[Proposition 4.24]{HKE}, which we need for the estimate.

\begin{lemma} \label{l:1}
    Assume $\VD(d_2)$, $\J_{\psi, \le}$ and $\E_\Phi$.  Then, there is a constant $c>0$ such that for any $\rho>0$, $t>0$ and $x \in M_0$,
	\begin{equation*}\label{e:4.24.1}
\bE^x \int_0^t \frac{1}{V(X_s^{\rho},\rho)} ds   \le \frac{ct}{V(x,\rho)} \left( 1+ \frac{t}{\Phi(\rho)}   \right)^{d_2+1}. 
	\end{equation*}
\end{lemma}

\pf 
Following the proof of \cite[Proposition 4.24]{HKE}, using $\J_{\psi, \le}$ we have
\begin{align}\begin{split}\label{e:4.24.2}
\bE^x \left[ \int_0^t \frac{1}{V(X_s^{\rho},\rho)} ds \right] =  \sum_{k=1}^\infty \bE^x \left[ \int_0^t \frac{1}{V(X_s^{\rho},\rho)}ds ; \tau_{B(x,k\rho)}^{\rho} > t \ge \tau_{B(x,(k-1)\rho)}^{\rho} \right] := \sum_{k=1}^\infty I_k. \quad \quad 
\end{split}\end{align}
When $t < \tau^{\rho}_{B(x,k\rho)}$, we have $d(X_s^{\rho},x) \le k\rho$ for all $s \le t$. This along with $\VD(d_2)$ yields that for all $ k \ge 1$ and $s \le t < \tau_{B(x,k\rho)}^{\rho}$,
\begin{equation}\label{e:4.241}
\frac{1}{V(X_s^{\rho},\rho) } \le \frac{c_1k^{d_2}}{V(X_s^{\rho},2k\rho)} \le \frac{c_1 k^{d_2}}{\displaystyle \inf_{d(z,x)\le k\rho} V(z,2k\rho)} \le \frac{c_1k^{d_2}}{V(x,\rho)}.
\end{equation}
On the other hand, by \cite[Corollary 4.22]{HKE}, there exist constants $c_i>0$ with $i=2,3,4$ such that for all $t,\rho>0$, $k \ge 1$ and $x \in M_0$,
\begin{equation}\label{l:4.22}
\P^x(\tau_{B(x,k\rho)}^{\rho} \le t ) \le c_2 \exp \Big( - c_3 k + c_4 \frac{t}{\Phi(\rho)} \Big).
\end{equation}
Let $k_0 = \lceil \frac{2c_4}{c_3} \frac{t}{\Phi(\rho)} \rceil + 1$. Using \eqref{e:4.241} and definition of $k_0$, we have
$$ \sum_{k=1}^{k_0} I_k \le \sum_{k=1}^{k_0} \frac{c_1 k^{d_2} t}{V(x,\rho)} \le \frac{c_5 k_0^{d_2+1} t}{V(x,\rho)} \le \frac{c_6 t}{V(x,\rho)} \left( 1+\frac{t}{\Phi(\rho)} \right)^{d_2+1}. $$
On the other hand, for any $k_0 < k$, using \eqref{e:4.241} and \eqref{l:4.22} with  $k_0 = \lceil \frac{2c_4}{c_3} \frac{t}{\Phi(\rho)} \rceil + 1$ we have 
$$ I_k \le \frac{c_1 k^{d_2} t}{V(x,\rho)} \P^x(\tau_{B(x,k\rho)}^{\rho} \le t ) \le \frac{c_1 c_2 k^{d_2} t}{V(x,\rho)} \exp \Big( - c_3 k + c_4 \frac{t}{\Phi(\rho)} \Big) \le \frac{c_1 c_2 k^{d_2} t}{V(x,\rho)} \exp \left(- \frac{c_3}{2} k  \right). $$
Thus, we conclude
$$ \sum_{k=k_0+1}^\infty I_k \le \frac{c_1 c_2 t}{V(x,\rho)}  \sum_{k=k_0+1}^\infty k^{d_2} e^{-\frac{c_3}{2} k}  =  \frac{c_7 t}{V(x,\rho)}. $$
From above two estimates, we obtain $\sum_{i=1}^\infty I_k \le \frac{c_8 t}{V(x,\rho)} (1+ \frac{t}{\Phi(\rho)})^{d_2+1}$. Combining this with  \eqref{e:4.24.2}, we obtain the desired estimate. \qed

\begin{lemma}   \label{l:4.241}
	Assume $\VD(d_2)$, $\J_{\psi, \le}$, $\UHKD(\Phi)$ and $\E_\Phi$. Then, there are constants $c>0$ and $C_1,C_2>0$ such that for any $\rho>0$, $t>0$ and $x,y \in M_0$, 
	\begin{align*}
	p(t,x,y) \le \frac{c}{V(x,\Phi^{-1}(t))}\Big(1+\frac{d(x,y)}{\Phi^{-1}(t)}\Big)^{d_2} \exp \Big( C_1 \frac{t}{\Phi(\rho)} - C_2 \frac{d(x,y)}{\rho} \Big) 
 +\frac{ct}{V(x,\rho)\psi(\rho)} \Big( 1+\frac{t}{\Phi(\rho)} \Big)_.^{d_2+1}
\end{align*}
\end{lemma}
\pf Recall that $X_t^{\rho}$ and $p^{\rho}(t,x,y)$ are the Hunt process and heat kernel correspond to $(\EE^\rho,\FF)$, respectively. Using \cite[Lemma 3.1 and (3.5)]{BGK}, 
\begin{align}\label{compose1}
p(t,x,y) \le p^{\rho}(t,x,y) + \bE^x \left[ \int_0^t \int_M J(X_s^{\rho},z)\1_{\{d(z,X_s^{\rho}) \ge \rho\}}(z) \,  p(t-s,z,y) \mu(dz)ds \right].
\end{align}
Also, using the symmetry of heat kernel, $\J_{\psi, \le}$ and Lemma \ref{l:1} we obtain
\begin{align}\label{compose2}
\begin{split} 
&\bE^x \left[ \int_0^t \int_M J(X_s^{\rho},z)\1_{\{d(z,X_s^{\rho}) \ge \rho\}}(z) \,  p(t-s,z,y) \mu(dz)ds \right] \\ &\le c_1 \bE^x \left[ \int_0^t \frac{1}{V(X_s^{\rho}, \rho)\psi(\rho)} ds \right] \le \frac{c_1 t}{V(x,\rho)\psi(\rho)} \left( 1+ \frac{t}{\Phi(\rho)} \right)^{d_2 +1}. 
\end{split}
\end{align}
Combining the estimates in \cite[Lemma 5.2]{HKE} and \eqref{compose2}, we conclude the proof of the lemma. Note that since $\J_{\psi, \le}$ and \eqref{comp1} imply $\J_{\Phi,\le}$, the conditions in  \cite[Lemma 5.2]{HKE} are satisfied. \qed

The proof of next lemma is same as that of  \cite[Lemma 4.2]{BKKL}. Thus, we skip the proof.
\begin{lemma}\label{l:7.11}
	Let  $r, t, \rho>0$. Assume that 
	\begin{align*}\label{l:7.11ass}\P^x(\tau_{B(x_0,r)}^{\rho} \le t) \le \phi(r,t) \quad \mbox{for all }x_0 \in M_0, \; x \in B(x_0,r/4) \cap M_0,\end{align*}
	where $\phi$ is a non-negative measurable function on $\R_+\times\R_+ $. Then, for any $k \in \N$,
	$$P_t^{\rho} \1_{B(x_0,k(r+\rho))^c}(x) \le \phi(r,t)^k \quad \mbox{for all }x_0\in M_0,  x \in B(x_0,r/4) \cap M_0. $$
\end{lemma}

\begin{lemma}\label{l:tail}
	Assume $\VD(d_2)$, $\J_{\psi, \le}$, $\UHKD(\Phi)$ and $\E_\Phi$. Let $T>0$ and $f:\bR_+\times\bR_+\to\bR_+$ be a measurable function satisfying that $t\mapsto f(r, t)$ is non-increasing for all $r>0$ and that $r\mapsto f(r, t)$ is non-decreasing for all $t>0$. Suppose that the following hold: (i)  For each $b>0$, $\sup_{ t \le T} f(b\Phi^{-1}(t), t) < \infty$ (resp., $\sup_{ t \ge T} f(b\Phi^{-1}(t), t) < \infty$); 
	(ii) there exist $\eta\in(0,\beta_1]$, $a_1>0$ and $c>0$ such that
	\begin{equation}\label{e:tail}
	\P^x\big(d(x,X_t)>r\big)
	\le c({\psi^{-1}(t)}/{r})^\eta+c\exp{\big(-a_1f(r, t)\big)}
	\end{equation} 
	for all $t\in (0,T]$ (resp. $t\in [T,\infty)$), $r>0$ and $x \in M_0$.
	
	Then, there exist constants $k \in \N, c_0>0$ such that 
	\begin{equation*}
	p(t,x,y) \le \frac{c_0\, t}{V(x,d(x,y))\psi(d(x,y))} + \frac{c_0}{V(x,\Phi^{-1}(t))}\left(1+\frac{d(x,y)}{\Phi^{-1}(t)}\right)^{d_2} \exp{\Big(-a_1k f\big(d(x,y)/(16k), t\big)\Big)} 
	\end{equation*}
	for all $t\in (0,T)$ (resp. $t\in [T,\infty)$) and $x, y\in M_0$.\end{lemma}

\pf
Since the proofs for cases $t \in (0, T]$ and $t \in [T,\infty)$ are similar, we only prove for $t\in(0,T]$. For $x_0 \in M_0$, let $B(r)=B(x_0,r) \cap M_0$.  By the strong Markov property, \eqref{e:tail}, and the fact that $t\mapsto f(r, t)$ is non-increasing,
we have that for $x \in B(r/4)$ and $t\in (0,T/2]$,
\begin{align}
\P^x(\tau_{B(r)} \le t)&  \le \P^x( X_{2t} \in B(r/2)^c) +\P^x(\tau_{B(r)} \le t, X_{2t} \in B(r/2))\nn\\ 
&\le \P^x(X_{2t} \in B(x, r/4)^c) + \sup_{ z \in B(r)^c, s \le t} \P^z(X_{2t-s} \in B(z,r/4)^c) \nn\\
&\le c(4{\psi^{-1}(2t)}/r)^\eta+c\exp{\big(-a_1f(r/4, 2t)\big)}. \label{e:tail2}
\end{align}
By \cite[Proposition 4.6]{GHL14} and \cite[Lemma 2.1]{HKE}, we have
\begin{align*}
& \Big|P_t^{B(r)} \1_{B(r)}(x)-P_t^{r,{B(r)}}\1_{B(r)}(x)\Big|
\le2t\,\sup_{z\in M}\int_{B(z,r)^c}J(z,y)\,\mu(dy)\le \frac{c_3 t}{\psi(r)}.
\end{align*}
Combining this with \eqref{e:tail2}, using $L(\beta_1,C_L,\psi)$ and Lemma \ref{lem:inverse}, we see that for all $x\in  B(r/4)$ and $t\in (0,T/2]$,
\begin{align}\label{e:tail3}
&\P^x(\tau^{r}_{B(r)} \le t) = 1-P_{t}^{r,B(r)} \1_{B(r)}(x)
\le 1-P_t^{B(r)}\1_{B(r)}(x) + \frac{c_3t}{\psi(r)}\nn\\
&\le c_1({\psi^{-1}(t)}/{r})^\eta+c_2\exp{\big(-a_1f(r/4, 2t)\big)}+c_3(t/{\psi(r)})
=: \phi(r,t).
\end{align}
Applying Lemma \ref{l:7.11} with $r=\rho$ to \eqref{e:tail3}, we see that for any $t \in (0,T/2]$, $x \in B(r/4)$ and $k \in \N$,
\begin{equation}\label{e:tail4}
\int_{B(x,2kr)^c} p^{r}(t,x,y)\mu(dy) =P_t^{r} \1_{B(x,2kr)^c}(x)  \le \phi(r,t)^k.
\end{equation}
Let $k:=\lceil \frac{\beta_2+2d_2}{\eta} \rceil +1$. 
For $t\in(0, T]$ and $x,y \in M_0$ satisfying $4k\Phi^{-1}(t) \ge d(x,y)$, by using that $r \mapsto f(r,t)$ is non-decreasing and the assumption $(i)$, we have
$f(d(x,y)/(16k), t)\le f(\Phi^{-1}(t)/4,t)\le C<\infty $. {Thus, using \cite[Lemma 5.1]{HKE},
\begin{equation}
\label{e:uhk1}
p(t,x,y) 
\le  \frac{c_5 e^{a_1kC}}{V(x,\Phi^{-1}(t))} \exp{\Big(-a_1kf(d(x,y)/(16k), t)\Big)}.
\end{equation} 
}
For the remainder of the proof, assume $t\in(0, T]$ and $4k\Phi^{-1}(t) < d(x,y)$. Also, denote $r=d(x,y)$ and $\rho={r}/{(4k)}$. Using \cite[(5.2)]{HKE}, \eqref{e:tail4} and \eqref{VD2}, we have
\begin{align}
&p^{\rho}(t,x,y) \le \Big( \int_{B(x,r/2)^c} + \int_{B(y,r/2)^c} \Big) p^{\rho}(t/2,x,z) p^{\rho}(t/2,z,y)\mu(dz) \nn\\
&\le \frac{c_6}{V(y,\Phi^{-1}(t))} \int_{B(x,2k\rho)^c}  p^{\rho}(t/2,x,z)\mu(dz)+\frac{c_6}{V(x,\Phi^{-1}(t))} \int_{B(y,2k\rho)^c}  p^{\rho}(t/2,z,y)\mu(dz)\nn\\
&\le  \frac{c_7}{V(x,\Phi^{-1}(t))}\Big(1+\frac{r}{\Phi^{-1}(t)}\Big)^{d_2} \phi(\rho,t/2)^k \le \frac{c_8}{V(x,\Phi^{-1}(t))}\Big(\frac{r}{\Phi^{-1}(t)}\Big)^{d_2} \phi(\rho,t/2)^k. \label{e:rfee}
\end{align} 
Note that $k\beta_1  \ge k\eta  \ge \beta_2 + 2d_2$ and $\rho \ge \Phi^{-1}(t) > \psi^{-1}(t)$. Thus, by  $L(\beta_1,C_L,\psi)$ we obtain
\begin{align*}
\left( \frac{\psi^{-1}(t)}{\rho}\right)^{\eta k}+\left(\frac{t}{\psi(\rho)}\right)^k  
 \le  \left( \frac{\psi^{-1}(t)}{\rho}\right)^{\beta_2 +2d_2} +c_9\left(\frac{\psi^{-1}(t)}{\rho}\right)^{k \beta_1} \le  c_{10}\left(\frac{\psi^{-1}(t)}{r}\right)^{\beta_2 +2d_2}.
\end{align*}
Applying this to \eqref{e:rfee} and using $\VD(d_2)$ and $U(\beta_2,C_U,\psi)$ we have
\begin{align*}
p^{\rho}(t,x,y)&\le \frac{c_{11}}{ V(x,\Phi^{-1}(t))}\Big(\frac{r}{\Phi^{-1}(t)}\Big)^{d_2} \Big( \Big( \frac{\psi^{-1}(t)}{\rho}\Big)^{\eta k}+\exp{\Big(-a_1kf(\rho/4, t)\Big)}+\Big(\frac{t}{\psi(\rho)}\Big)^k \Big) \\
&\le \frac{c_{12}}{ V(x,r)}
\frac{V(x,r)}{ V(x,\Phi^{-1}(t))}
\Big(\frac{r}{\Phi^{-1}(t)}\Big)^{d_2} \Big( \Big(\frac{\psi^{-1}(t)}{r}\Big)^{\beta_2 + 2d_2} + \exp{\big(-a_1kf(\rho/4, t)\big)} \Big) 
\\
&\le \frac{c_{13} t}{V(x,r)}\left(\frac{\psi^{-1}(t)}{r}\right)^{\beta_2 }
 +  \frac{c_{12}}{V(x,\Phi^{-1}(t))}\left(\frac{r}{\Phi^{-1}(t)}\right)^{d_2}\exp{\Big(-a_1kf(r/(16k), t)\Big)}
\\
&\le \frac{c_{14} t}{V(x,r)\psi(r)} +  \frac{c_{12}}{V(x,\Phi^{-1}(t))}\left(\frac{r}{\Phi^{-1}(t)}\right)^{d_2}\exp{\Big(-a_1kf(r/(16k), t)\Big)}.
\end{align*} 
Thus, by {\eqref{compose1}, \eqref{compose2}} and $U(\beta_2, C_U,\psi)$, 
we conclude that for any $t \in (0,T]$ and $x,y \in M_0$ with $4k\Phi^{-1}(t) < d(x,y)$, 
\begin{align}\label{e:uhk2}
\begin{split}
&p(t,x,y) \le p^{\rho}(t,x,y) + \frac{c_{15}t}{V(x,\rho) \psi(\rho)} \left(1+\frac{t}{\Phi(\rho)}\right)^{d_2 + 1} \\
& \le \frac{ c_{16}}{V(x,\Phi^{-1}(t))}\left(1+\frac{d(x,y)}{\Phi^{-1}(t)}\right)^{d_2}\exp{\Big(-a_1kf\big(\frac{d(x,y)}{16k}, t\big)\Big)} + \frac{c_{16} t}{V(x,d(x,y))\psi(d(x,y))}. 
\end{split}
\end{align}
Here in the second inequality we have used $\Phi(\rho) \ge t$. Now the conclusion follows from \eqref{e:uhk1} and \eqref{e:uhk2}. 

\qed

\begin{lemma}\label{l:exp}
	Suppose $\VD(d_2)$,  $\J_{\psi,\le}$, $\UHKD(\Phi)$ and $\E_\Phi$. Then, there exist constants $a_0,c>0$ and $N \in \N$ such that 
	\begin{equation}\label{e:expN}
	p(t,x,y) 
	\le \frac{c\, t}{V(x,d(x,y))\psi(d(x,y))} + c\,V(x,\Phi^{-1}(t))^{-1} \exp{\left(-\frac{a_0 d(x,y)^{1/N}}{\Phi^{-1}(t)^{1/N}}\right)}, 
	\end{equation}
	for all $t>0$ and $x,y\in M_0$.
\end{lemma}

\pf Let $N:=\lceil\frac{\beta_1 + d_2}{\beta_1} \rceil+1$, and $\eta:=\beta_1-(\beta_1+d_2)/N>0$.  
We first claim that there exist $a_1>0$ and $c_1>0$ such that for any $t,r>0$ and $x \in M_0$,
\begin{equation}\label{e:tailN}
\int_{\{ y:d(x,y) \ge r\}  } p(t, x,y)\mu(dy)\le c_1\left(\frac{\psi^{-1}(t)}{r}\right)^\eta+c_1\exp{\left(-\frac{a_1 r^{1/N}}{\Phi^{-1}(t)^{1/N}}\right)}.
\end{equation} 
 When $r \le \Phi^{-1}(t)$, we immediately obtain \eqref{e:tailN} by letting $c = \exp(a_1)$. Thus, we will only consider the case $r>\Phi^{-1}(t)$.
Fix $\alpha \in (d_2/(d_2 +\beta_1),1)$ and define for $n \in \N$, 
$$
\rho_n=\rho_n(r,t)=2^{n\alpha}r^{1-1/N}\Phi^{-1}(t)^{1/N}.
$$
Since $r>\Phi^{-1}(t)$ , we have
$
\Phi^{-1}(t) < \rho_n \leq 2^nr.
$
In particular, $t \le \Phi(\rho_n)$. Thus, using Lemma \ref{l:4.241} with $\rho=\rho_n$, we have that for every $t>0$ and $x,y \in M_0$ with $2^nr \le d(x,y) < 2^{n+1}r$,
\begin{align*}
&p(t,x,y) \\ &\le \frac{c_2}{V(x,\Phi^{-1}(t))}\left(\frac{2^{n+1} r}{\Phi^{-1}(t)}\right)^{d_2+1} \exp \left( C_1 \frac{t}{\Phi(\rho_n)} - C_2 \frac{d(x,y)}{\rho_n} \right) +\frac{c_2t}{V(x,\rho_n)\psi(\rho_n)}\left( 1+\frac{t}{\Phi(\rho_n)} \right)^{d_2+1} \\
&\le \frac{c_3}{V(x,\Phi^{-1}(t))} \left(\frac{2^n r}{\Phi^{-1}(t)}\right)^{d_2+1}  \exp\left(-C_2 \frac{2^n r}{\rho_n}\right) + \frac{c_3 t}{V(x,\rho_n) \psi(\rho_n)} \,  \\
&= \frac{c_3}{V(x,\Phi^{-1}(t))}\left(\frac{2^n r}{\Phi^{-1}(t)}\right)^{d_2+1} \exp\left(-C_2 \frac{2^{n(1-\alpha)}r^{1/N}} {\Phi^{-1}(t)^{1/N}} \right) + \frac{c_3t}{V(x,\rho_n) \psi(\rho_n)}.
\end{align*}
Using the above estimate and $\VD(d_2)$ we get that
\begin{align*}
&\int_{B(x,r)^c}p(t,x,y)\mu(dy)
=\sum_{n=0}^{\infty}\int_{B(x,2^{n+1}r)\setminus B(x,2^nr)}p(t,x,y)\mu(dy)\\
&\leq c_3\sum_{n=0}^{\infty} \frac{V(x,2^{n+1} r)}{V(x,\Phi^{-1}(t))} \left(\frac{2^n r}{\Phi^{-1}(t)}\right)^{d_2 +1}\exp\left(-C_2 \frac{2^{n(1-\alpha)}r^{1/N}} {\Phi^{-1}(t)^{1/N}} \right) + c_3\sum_{n=0}^\infty \frac{ V(x,2^{n+1} r)}{V(x,\rho_n)} \frac{t}{ \psi(\rho_n)}\\&=:I_1+I_2.
\end{align*}
We first estimate $I_1$. Observe that for any $d_0 \ge 1$, there exists $c_1=c_1(c_0, \alpha)>0$ such that $2n \le \frac{c_0}{2d_0} 2^{n(1-\alpha)} + c_1$ holds for every $n\ge0$. Thus, for any $n\ge0$ and $\kappa \ge1$,
\begin{align}\label{e:3.4.1}
& 2^{nd_0} \exp\left(-C_2 2^{n(1-\alpha)}\kappa \right)\le  2^{-nd_0}\exp\left(2nd_0-C_2 2^{n(1-\alpha)}\kappa \right)\nn \\ 
&\le 2^{-nd_0}\exp\left( \left(\frac{C_2 }{2d_0} 2^{n(1-\alpha)} + c_1 \right)d_0-C_2 2^{n(1-\alpha)}\kappa \right)\nn\\
& \le 2^{-nd_0}\exp\left( \frac{C_2}{2} 2^{n(1-\alpha)}\kappa  + c_1d_0 -C_2 2^{n(1-\alpha)}\kappa \right)  = e^{c_1d_0} 2^{-nd_0} \exp \left(- \frac{C_2 }{2}\kappa  \right).
\end{align}
Using $\Phi^{-1}(t) < r$, $\VD(d_2)$, \eqref{e:3.4.1} with $d_0 = 2d_2+1$ and $\kappa = \big(\frac{r}{\phi^{-1}(t)}\big)^{1/N}$, and the fact that 
$$\displaystyle \sup_{s \ge 1} s^{2d_2+1} \exp\left(-\frac{C_2}{4} s^{1/N}\right):= c_4 < \infty,$$
we obtain
\begin{align}
\begin{split}
I_1&\le c_5\sum_{n=0}^{\infty}\left(\frac{r}{\Phi^{-1}(t)}\right)^{2d_2+1}2^{n(2d_2+1)}\exp\left(-C_2\frac{2^{n(1-\alpha)}r^{1/N}}{\Phi^{-1}(t)^{1/N}}\right) \\
&\leq c_6\left(\frac{r}{\Phi^{-1}(t)}\right)^{2d_2+1}\exp\left(-\frac{C_2}{2}\frac{r^{1/N}}{\Phi^{-1}(t)^{1/N}}\right)\sum_{n=0}^{\infty}2^{-n(2d_2+1)} \leq c_4c_6\exp\left(-\frac{C_2 r^{1/N}}{4\Phi^{-1}(t)^{1/N}}\right). \label{e:wq1}
\end{split}
\end{align}
 We next estimate $I_2$. Note that by \eqref{comp1} and $t<\Phi(\rho_n)$, we have $\psi^{-1}(t) \le \Phi^{-1}(t) \le \rho_n$. 
Thus, using $\VD(d_2)$ and $L(\beta_1,C_L,\psi)$ for the first line and using $\alpha(d_2+\beta_1) > d_2$ for the second line, we obtain 
\begin{align*}
I_2 &
\leq c_{7}\sum_{n=0}^\infty \left(\frac{2^nr}{\rho_n}\right)^{d_2}\left(\frac{\psi^{-1}(t)}{\rho_n}\right)^{\beta_1}= c_7 \left(\frac{\Phi^{-1}(t)}{r}\right)^{-\frac{d_2+\beta_1}{N}}\left(\frac{\psi^{-1}(t)}{r}\right)^{\beta_1}\sum_{n=0}^\infty 2^{n(d_2-\alpha(d_2 + \beta_1))}\\
& \leq c_8 \left(\frac{\psi^{-1}(t)}{r}\right)^{\beta_1-\frac{d_2+\beta_1}{N}}
=c_8 \left(\frac{\psi^{-1}(t)}{r}\right)^{\eta}.
\end{align*}
Thus, by above estimates of $I_1$ and $I_2$, we obtain \eqref{e:tailN}.

 By $\eta<\beta_1$ and \eqref{e:tailN}, assumptions in Lemma \ref{l:tail} hold with $f(r,t):=\big(r/\Phi^{-1}(t)\big)^{1/N}$. Thus, by  Lemma \ref{l:tail}, we have
 $$ p(t,x,y) \le \frac{c_{10}t}{V(x,d(x,y))\psi(d(x,y))} + \frac{c_{10}}{V(x,\Phi^{-1}(t))} \Big( 1+ \frac{d(x,y)}{\Phi^{-1}(t)} \Big)^{d_2} \exp \Big[ -a_1 k \Big( \frac{d(x,y)}{16k\Phi^{-1}(t)} \Big)^{1/N} \Big]. $$
 Using the fact that $\sup_{s>0} (1+s)^{d_2} \exp (- c s^{1/N}  ) < \infty$ for every $c>0$, we conclude \eqref{e:expN}. \qed

The next lemma will be used in the next subsection and Section \ref{s:main3proof}.

\begin{lemma}\label{l:exp1}
Suppose $\VD(d_2)$, $\J_{\psi, \le}$, $\UHKD(\Phi)$ and $\E_\Phi$. Then, for any $\theta>0$ and $c_0, c_1 \ge 1$, there exists $c>0$ such that for any $x \in M_0$, $t>0$ and $r \ge c_0 \frac{\Phi^{-1}( c_1t)^{1+\theta}}{\psi^{-1}(c_1t)^{\theta}}$,
\begin{equation*}\label{e:tail1}
\int_{B(x,r)^c} p(t,x,y)\mu(dy) \le c \left(\frac{\psi^{-1}(t)}{r} \right)^{\beta_1}.
\end{equation*}
\end{lemma}
\pf 
Denote $t_1 = c_1 t$ and let $a_0,N$ be the constants in Lemma \ref{l:exp}. By \eqref{comp1} we have that for any $y \in M_0$ with $d(x,y) > r$, there exists $\theta_0 \in (\theta,\infty)$ satisfying $d(x,y)= c_0\Phi^{-1}(t_1)^{1+\theta_0}/\psi^{-1}(t_1)^{\theta_0}$. Note that there exists a positive constant $c_2=c_2(\theta)$ such that for any $s>0$,
\begin{equation}\label{e:tail1.1}
s^{-d_2 - \beta_2 - \beta_2/\theta} \ge c_2\exp(-a_0 s^{1/N}).
\end{equation}
Also, since $c_0 \ge 1$ we have
$$  \psi^{-1}(t_1) \le c_0\psi^{-1}(t_1) \le c_0 \Phi^{-1}(t_1)< d(x,y) = c_0 \Phi^{-1}(t_1)^{1+\theta_0} / \psi^{-1}(t_1)^{\theta_0}. $$ 
Thus, using $\VD(d_2)$ and $U(\beta_2, C_U,\psi)$ for the first inequality and \eqref{e:tail1.1} for the second, we have
\begin{align*}
&\frac{t}{V(x,d(x,y))\psi(d(x,y))}
=\frac{c_1^{-1}}{V(x,c_0\Phi^{-1}(t_1))}\frac{V(x,c_0\Phi^{-1}(t_1))}{V(x,d(x,y))}\frac{\psi(  \psi^{-1}(t_1))}{\psi(d(x,y))}\nn\\ 
&\ge \frac{ c_1^{-1}}{V(x,c_0\Phi^{-1}(t_1))} C_\mu^{-1} \left( \frac{c_0 \Phi^{-1}(t_1)}{d(x,y)} \right)^{d_2} C_U^{-1} \left( \frac{\psi^{-1}(t_1)}{d(x,y)} \right)^{\beta_2} \\
&= \frac{ c_1^{-1}c_0^{-\beta_2}C_\mu^{-1} C_U^{-1}}{V(x,c_0\Phi^{-1}(t_1))}\left(\frac{\psi^{-1}(t_1)}{\Phi^{-1}(t_1)}\right)^{d_2\theta_0} \left(\frac{\psi^{-1}(t_1)}{\Phi^{-1}(t_1)}\right)^{(1+\theta_0)\beta_2}\nn \\
&= \frac{c_1^{-1} c_0^{-\beta_2}C_U^{-1} C_\mu^{-1}}{V(x,c_0\Phi^{-1}(t_1))}\left(\left(\frac{\Phi^{-1}(t_1)}{\psi^{-1}(t_1)}\right)^{\theta_0}\right)^{-d_2-\beta_2-\beta_2/\theta_0} \nn\\
&\ge \frac{c_2 c_1^{-1}c_0^{-\beta_2}C_U^{-1}C_\mu^{-1}c_U^{-1}}{V(x,c_0\Phi^{-1}(t_1))}\exp{\left(-\frac{a_0\Phi^{-1}(t_1)^{\theta_0/N}}{\psi^{-1}(t_1)^{\theta_0/N}}\right)}= \frac{c_2 c_1^{-1}c_0^{-\beta_2}C_U^{-1}C_\mu^{-1}c_U^{-1}}{V(x,c_0\Phi^{-1}(t_1))}\exp{\left(-\frac{a_0d(x,y)^{1/N}}{\Phi^{-1}(t_1)^{1/N}}\right)}.
\end{align*}
Applying Lemma \ref{lem:inverse} for $L(\alpha_1,c_L,\Phi)$, we have $U( 1/\alpha_1,c_L^{-1/\alpha_1}, \Phi^{-1})$, which yields $\Phi^{-1}(t) \le \Phi^{-1}(t_1) \le c_L^{-1/\alpha_1} c_1^{1/\alpha_1}\Phi^{-1}(t)$. Thus, using this and $\VD(d_2)$ again, we have
\begin{align*}
& \frac{1}{V(x,c_0\Phi^{-1}(t_1))} \exp{\Big(-\frac{a_0 d(x,y)^{1/N}}{\Phi^{-1}(t_1)^{1/N}}\Big)} \ge \frac{1}{V(x,c_0c_L^{-1/\alpha_1}c_1^{1/\alpha_1}\Phi^{-1}(t))}\exp{\Big(-a_0\frac{  d(x,y)^{1/N}}{\Phi^{-1}(t)^{1/N}}\Big)} \\
& \ge C_\mu^{-1} (c_0c_L^{-1/\alpha_1}c_1^{1/\alpha_1})^{d_2} \frac{1}{V(x,\Phi^{-1}(t))} \exp \Big( -a_0 \frac{d(x,y)^{1/N}}{\Phi^{-1}(t)^{1/N}}\Big). 
\end{align*}
Thus, by Lemma \ref{l:exp} and above two estimates, we have that for every $y \in M_0$ with $d(x,y)>r$,
\begin{align*}
p(t,x,y) &\le \frac{c \, t}{V(x,d(x,y))\psi(d(x,y))} + \frac{c}{V(x,\Phi^{-1}(t))} e^{-a_0\frac{ d(x,y)^{1/N}}{\Phi^{-1}(t)^{1/N}} } \le \frac{c_2\, t}{V(x,d(x,y))\psi(d(x,y))}.
\end{align*}
Using this, \cite[Lemma 2.1]{HKE} and $L(\beta_1,C_L, \psi)$ with the fact that $r> c_0\psi^{-1}(c_1t)$ which follows from \eqref{comp1}, we conclude that
$$
\int_{B(x,r)^c}p(t,x,y)\mu(dy) \le c_2 \int_{B(x,r)^c} \frac{t}{V(x,d(x,y)) \psi(d(x,y))} \mu(dy) 
\le c_3 \frac{t}{\psi(r)} \le c_4\left(\frac{\psi^{-1}(t)}{r}\right)^{\beta_1}.
$$
This proves the lemma. \qed

\subsection{Basic properties of $\sT(\phi)$}\label{s:tran}

In this subsection, we assume that $\phi : \R_+ \rightarrow \R_+$ is a non-decreasing function.
Recall that $\phi^{-1}(t):= \inf \{ s \ge 0 : \phi(s) >t\}$ is the generalized inverse function of $\phi$.
\begin{lemma}\label{l:ginv}
Suppose that  $\phi$  satisfies $U(\alpha_2, c_U)$ with some $\alpha_2>0$. Then 
	\begin{equation}\label{e:ginv}
	c_U^{-1} t \le \phi(\phi^{-1}(t)) \le c_U t, \qquad \text{ for all } t>0.
	\end{equation}
\end{lemma}
\pf Let $\phi^{-1}(t) = s$. Since $\phi$ is non-decreasing, we have $\phi(u) \le t$ for all $u < s$. Thus, using $U(\alpha_2,c_U,\phi)$ we have 
$ \phi(s) \le c_U ( s/u )^{\alpha_2} \phi(u) \le c_U ( {s}/{u})^{\alpha_2} t. $
Letting $u \uparrow s$ we obtain $\phi(s) \le c_U t$. 

By the similar way, using the fact that $\phi(u) \ge t$ for all $u > s$, we have
$ \phi(s) \ge c_U^{-1} (s/u)^{\alpha_2} \phi(u) \ge c_U^{-1}  (s/u)^{\alpha_2} t,$
which yields $\phi(s) \ge c_U^{-1} t$ by letting $u\downarrow s$. 
\qed

 For the remainder of this subsection, we assume that  $\phi$ satisfies $L(\alpha_1,c_L)$ and $U(\alpha_2, c_U)$ with some $0<\alpha_1 \le \alpha_2$, and satisfies $L_a(\delta,\wt C_L)$ for some $a>0$ and $\delta>1$. For such function $\phi$, we define \begin{equation}\label{e:trans}
\sT(\phi)(r,t):= \sup_{s>0} \left[ \frac{r}{s} - \frac{t}{\phi(s)} \right], \quad r,t>0. 
\end{equation}
Note that from $L(\alpha_1,c_L,\phi)$ and $L_a(\delta, \wt C_L)$, we obtain $\displaystyle \lim_{s \to  \infty}\phi(s)=\infty$ and $\displaystyle \lim_{s \to 0} \frac{\phi(s)}{s} = 0$, respectively. This concludes that $\sT(\phi)(r,t) \in [0,\infty)$ for all $r,t>0$. Also, comparing the definitions in \eqref{d:Phi1} and \eqref{e:trans}, we see that $\sT(\Phi)= \Phi_1$ for instance. It immediately follows from the definition of $\sT(\phi)$ that for any $c,r,t>0$,
$
\sT(\phi)(cr,ct) = c\sT(\phi)(r,t).
$

We first observe when the supremum in \eqref{e:trans} occurs.
\begin{lemma} \label{l:max}
 Let $\delta_1 := \frac{1}{\delta-1}$. For any $T \in (0,\infty)$, there exists constant $b \in (0,1)$ such that for any $r>0$, $t \in (0,T]$ with $r \ge 2c_U \phi^{-1}(t)$, 
 \begin{equation}\label{e:max}
	 \sT(\phi)(r,t)= \sup_{s \in [ br^{-\delta_1}\phi^{-1}(t)^{\delta_1+1}, 2\phi^{-1}(t)]} \left[ \frac{r}{s} - \frac{t}{\phi(s)} \right] 
 \ge \frac{r}{2\phi^{-1}(t)}.	\end{equation}	
Moreover, if $L(\delta, \wt C_L,\phi)$ holds, \eqref{e:max} holds for all $t \in (0,\infty)$.
\end{lemma}
\pf From Remark \ref{mwsc}, we may and do assume $a  =\phi^{-1}(T)$ without loss of generality.  \\
Denote $b = (c_U^{-1} \wt C_L)^{\delta_1} \in (0,1)$.
Fix $r>0$ and $t \in (0,T]$ with $r \ge 2c_U\phi^{-1}(t)$ and let us define $$\displaystyle \sT_1(\phi)(r,t):=\sup_{s \in [br^{-\delta_1}\phi^{-1}(t)^{\delta_1+1}, 2\phi^{-1}(t)]} \left[ \frac{r}{s} - \frac{t}{\phi(s)} \right].$$
  Since $r \ge 2\phi^{-1}(t)$, we have
$
  r^{-\delta_1} \phi^{-1}(t)^{\delta_1 +1}  \le \phi^{-1}(t), 
$ 
which yields $\phi^{-1}(t) \in [ br^{-\delta_1}\phi^{-1}(t)^{\delta_1 +1} , 2\phi^{-1}(t)]$. 

\noindent Now, taking $s = \phi^{-1}(t)$ for \eqref{e:max}. Using \eqref{e:ginv} and $r \ge 2c_U\phi^{-1}(t)$ we have
\begin{equation}\label{e:max2}
\sT_1(\phi)(r,t) \ge \frac{r}{\phi^{-1}(t)} - \frac{t}{\phi(\phi^{-1}(t))} \ge \frac{r}{\phi^{-1}(t)} - c_U \ge   \frac{r}{2\phi^{-1}(t)}.
\end{equation}
Assume $s > 2\phi^{-1}(t)$. Then, we have
$ \frac{r}{s} - \frac{t}{\phi(s)} \le \frac{r}{2\phi^{-1}(t)}. $
Thus, by \eqref{e:max2} we obtain
$ \frac{r}{s} - \frac{t}{\phi(s)} \le \sT_1(\phi)(r,t) $
for $s > 2\phi^{-1}(t)$.

 Now assume $s < br^{-\delta_1} \phi^{-1}(t)^{\delta_1+1}$. Since $s \le \phi^{-1}(t) \le \phi^{-1}(T)$, using \eqref{e:ginv}, $L_{\phi^{-1}(T)}(\delta,\wt C_L,\phi)$ and $r \ge 2c_U\phi^{-1}(t) $ we have
\begin{align*}
\frac{r}{s} - \frac{t}{\phi(s)} &\le \frac{r}{s} - c_U^{-1}\frac{\phi(\phi^{-1}(t))}{\phi(s)} \le \frac{r}{s} - c_U^{-1}\wt C_L \Big( \frac{\phi^{-1}(t)}{s} \Big)^\delta \le \frac{\phi^{-1}(t)^\delta}{s} \Big( \frac{r}{\phi^{-1}(t)^\delta} - c_U^{-1}\wt C_L s^{1- \delta} \Big) \\ &\le \frac{\phi^{-1}(t)^\delta}{s} \Big( \frac{r}{\phi^{-1}(t)^\delta} - c_U^{-1}\wt C_L (br^{-\delta_1}\phi^{-1}(t)^{\delta_1+1} )^{1- \delta} \Big)  = 0,
\end{align*}
where we have used $c_U^{-1}\wt C_L b^{1-\delta} = 1$ for the last line.
Therefore, by \eqref{e:max2} we obtain $\frac{r}{s} - \frac{t}{\phi(s)} \le 0 \le \sT_1(\phi)(r,t)$. Combining above two cases and the definition of $\sT_1(\phi)$ we conclude that $\frac{r}{s} - \frac{t}{\phi(s)} \le \sT_1(\phi)(r,t)$ for every $s>0$. This concludes the lemma. \qed
\begin{lemma}\label{l:sT}
	(i) For any $T>0$ and $c_1,c_2>0$, there exists a constant $c>0$ such that for any $r>0$ and $t \in (0,T]$ with $r \ge 2c_U\phi^{-1}(t)$,
	\begin{equation}\label{e:sT}
	\sT(\phi)(c_1r,c_2 t) \le c\sT(\phi)( r,  t). 
	\end{equation}
	(ii) For any $T>0$ and $c_3>0$, there exists a constant $\tilde{c}>0$ such that for any $t \in (0,T]$ and $r \le c_3\phi^{-1}(t)$,
	\begin{equation}\label{e:sT2}
		\sT(\phi)(r,t) \le \tilde{c}.
	\end{equation}
Moreover, if $L(\delta, \wt C_L,\phi)$ holds, both \eqref{e:sT} and \eqref{e:sT2} hold for all $t \in (0,\infty)$.
\end{lemma}
\pf Without loss of generality we may and do assume $a = 2\phi^{-1}(T)$.

\noindent (i) Since
$$ \sT(\phi)(c_1 r, c_2 t) = c_2 \sT(\phi)(\frac{c_1}{c_2} r, t),$$
it suffices to show that for any $c_4>0$, there exists a constant $c(c_4)>0$ such that for any $r>0$ and $t \in(0,T]$ with $r \ge 2c_U\phi^{-1}(t)$,
\begin{equation}\label{e:sT1}
 \sT(\phi)(c_4 r, t) \le c \sT(\phi)(r,t). 
\end{equation}
Also, since $\sT(\phi)(r,t)$ is increasing on $r$, we may and do assume that $c_4 \ge 1$. Since $c_4 r \ge r \ge 2c_U\phi^{-1}(t)$, by Lemma \ref{l:max} we have
\begin{equation*}
 \sT(\phi)(c_4 r, t) =\sup_{s \in (0, 2\phi^{-1}(t)]} \left[ \frac{c_4 r}{s} - \frac{t}{\phi(s)} \right]. 
\end{equation*} 
Let $\theta = (c_4^{-1} \wt C_L)^{1/(\delta-1)} \le 1$, which satisfies $ \theta = c_4\wt C_L^{-1} \theta^\delta$. Firstly, for any $s \in (0, 2\theta \phi^{-1}(t)]$, we have
\begin{align*}
\frac{r}{s/\theta} - \frac{t}{\phi(s/\theta)} = \frac{\theta r}{ s} - \frac{\phi(s)}{\phi(s/\theta)} \frac{t}{\phi( s)} \ge \frac{\theta r}{ s} - \wt C_L^{-1} \theta^\delta \frac{t}{\phi( s)} = \wt C_L^{-1} \theta^\delta \Big( \frac{c_4 r}{s} - \frac{t}{\phi(s)} \Big),
\end{align*}
where we have used $L_{2\phi^{-1}(T)}(\delta,\wt C_L,\phi)$ and $s \le s/\theta \le 2\phi^{-1}(t)$. Thus, we have
$$ \sup_{s \in (0,2\theta \phi^{-1}(t)]}\left[ \frac{c_4 r}{ s} - \frac{t}{\phi( s)} \right] \le \wt C_L \theta^{-\delta} \sT(\phi)(r,t). $$
Also, using \eqref{e:max} and $r \ge 2c_U\phi^{-1}(t)$, for any $s \in (2\theta \phi^{-1}(t), 2\phi^{-1}(t)]$ we have
$$ \frac{c_4 r}{s} - \frac{t}{s} \le \frac{c_4 r}{2\theta \phi^{-1}(t)} \le \frac{c_4}{\theta} \sT(\phi)(r,t). $$
Combining above two inequalities, we obtain the desired estimate.
 
\noindent (ii) Since $\sT(\phi)(r,t)$ is increasing in $r$, we may and do assume that $r = c_3\phi^{-1}(t)$ and $c_3 \ge 2c_U$. Observe that by $c_3 \ge 2c_U$, \eqref{e:ginv}, Lemma \ref{l:max} and $L_{2\phi^{-1}(T)}(\delta, \wt C_L,\phi)$ we have that for any $t \le T$,
\begin{align*}
& \sT(\phi)(c_3\phi^{-1}(t),t) = \sup_{s \le 2\phi^{-1}(t)} \left[ \frac{c_3\phi^{-1}(t)}{s} - \frac{t}{\phi(s)} \right] 
\le \sup_{s \le 2\phi^{-1}(t)} \left[ \frac{c_3\phi^{-1}(t)}{s} - c_U^{-1}\frac{\phi(\phi^{-1}(t))}{\phi(2\phi^{-1}(t))} \frac{\phi(2\phi^{-1}(t))}{\phi(s)} \right] 
 \\ &\le \sup_{s \le 2\phi^{-1}(t)} \left[ \frac{c_3\phi^{-1}(t)}{s} - c_5  \Big( \frac{ 2\phi^{-1}(t)}{s} \Big)^\delta \right] = \sup_{u>0} \big( c_3u - c_52^\delta u^\delta \big) := \tilde c < \infty.
 \end{align*}
 Here in the last line we have used $\delta>1$ to obtain $\tilde c<\infty$. This proves \eqref{e:sT2}.  \qed

\subsection{Proofs of Theorems \ref{t:main10} and \ref{t:main11}}

In this subsection, we prove our first main results. We start with the proof of Theorem \ref{t:main10}. \\

\noindent \textit{Proof of Theorem \ref{t:main10}.} Note that under the condition $L^a(\delta, \wt C_L, \Phi)$, we have $\alpha_2\ge \delta\vee \alpha_1$. Let  $C_1$ and $C_2$ be the constants in Lemma \ref{l:4.241} and,  without loss of generality, we may and  do assume that $C_1 \ge 2$ and $C_2 \le 1$. Take $$ 
\theta:=\frac{(\delta-1)\beta_1}{\delta(2d_2+\beta_1) + (\beta_1+2\alpha_2+2d_2\alpha_2)} \in (0, \delta-1) \quad \text{and}\quad C_0= \frac{4c_U}{C_2}.$$  Let $\alpha$ be a number in $(\frac{d_2}{d_2+\beta_1},1)$.

\noindent (i) We will show that there exist $a_1>0$ and $c_1>0$   such that for any $t\leq T$, $x \in M_0$ and $r>0$,
\begin{equation}\label{e:tail-F}
\int_{B(x,r)^c} p(t, x,y)\,\mu(dy)\le c_1\bigg(\frac{\psi^{-1}(t)}{r}\bigg)^{\beta_1/2}+c_1\exp\left(-a_1\Phi_1(r,t)\right).
\end{equation}
Firstly, since $\Phi_1(r,t)$ is increasing on $r$, by  \eqref{e:sT} and \eqref{e:sT2} we have that for $r \le C_0\Phi^{-1}(C_1t)$,
$$ \Phi_1(r,t) \le \Phi_1(C_0 \Phi^{-1}(C_1t), t) \le c_2 \Phi_1(C_0 \Phi^{-1}(C_1 t), C_1t) \le c_3. $$
Here for the second inequality, $C_0  \ge 2c_U$ yields the condition in \eqref{e:sT}. Thus, for any $x \in M_0$ and $r \le C_0\Phi^{-1}(C_1 t)$ we have
\begin{equation*}
\label{tail1}
\int_{B(x,r)^c} p(t,x,y) \mu(dy) \le 1 \le  e^{a_1c_3} \exp{ \left(-a_1 \Phi_1(r,t)\right)}. 
\end{equation*}
Also, when $r >  C_0 \Phi^{-1}(C_1 t)^{1+\theta} / \psi^{-1}(C_1t)^\theta$, \eqref{e:tail-F} immediately follows from Lemma \ref{l:exp1} and the fact that $r>\psi^{-1}(t)$.

Now consider the case $C_0\Phi^{-1}(C_1t) < r \le  C_0\Phi^{-1}(C_1t)^{1+\theta} / \psi^{-1}(C_1t)^{\theta}$. In this case, there exists $\theta_0 \in (0,\theta]$ such that $r= C_0 \Phi^{-1}(C_1t)^{1+\theta_0}/\psi^{-1}(C_1t)^{\theta_0}$ by \eqref{comp1}. Since $C_0 = \frac{4c_U}{C_2}$, applying Lemma \ref{l:max} with the constant $C_1T$ we have $\rho \in [b (C_2r/2)^{-\delta_1} \Phi^{-1}(C_1t)^{\delta_1+1}, 2\Phi^{-1}(C_1t)]$ such that
$$\Phi_1 ( C_2r/2,C_1 t) - \frac{C_0C_2}{8} \le \frac{C_2 r}{2\rho}- \frac{C_1 t}{\Phi(\rho)} \le \Phi_1 ( C_2r/2,C_1 t), $$
where $\delta_1 = \frac{1}{\delta -1}$. Also, let $\rho_n = 2^{n\alpha} \rho$ for $n \in \N_0$.
Then, we have
\begin{align*}
\frac{C_2 2^n r}{\rho_n} = \frac{C_2 r}{2\rho}+\frac{C_2 r}{\rho} (2^{n(1-\alpha)}-\frac{1}{2})\ge \frac{C_2r}{2\rho}+\frac{C_2r}{2\rho} 2^{n(1-\alpha)}\ge \frac{C_2r}{2\rho} + \frac{C_2 r}{4\Phi^{-1}(C_1 t)} 2^{n(1-\alpha)}.
\end{align*}
Using this, $r>C_0\Phi^{-1}(C_1 t)$ and \eqref{e:sT} with $U(1/\alpha_1, c_L^{-1/\alpha_1},\Phi^{-1})$, which follows from $L(\alpha_1,c_L,\Phi)$ and Lemma \ref{lem:inverse}, yield that for any $n \in \N_0$, 
\begin{align}
&\frac{C_1 t}{\Phi(\rho_n)} - \frac{C_2 2^n r}{\rho_n} \le \frac{C_1t}{\Phi(\rho)} - \frac{C_2 r}{2\rho} - \frac{C_2 r}{4\Phi^{-1}(C_1t)}2^{n(1-\alpha)}  
\\&\le - \Phi_1(C_2r/2,C_1t) + \frac{C_0C_2}{8}- \frac{C_2 r}{4\Phi^{-1}(C_1t)}2^{n(1-\alpha)}  \le - \Phi_1(C_2r/2,C_1t) - \frac{C_2 r}{8\Phi^{-1}(C_1t)}2^{n(1-\alpha)} \\
&\le - \Phi_1(C_2r/2,C_1t) - \frac{C_2}{8}(c_L/C_1)^{1/{\alpha_1}} 2^{n(1-\alpha)}\frac{r}{\Phi^{-1}(t)}  \le -c_4 \Phi_1 (r,t)  -c_5 2^{n(1-\alpha)}\frac{r}{\Phi^{-1}(t)}.
\label{e:A8.1}
\end{align}
By \eqref{e:A8.1} and Lemma \ref{l:4.241} with $\rho=\rho_n$, we have that for $t \in (0,T]$ and $y \in B(x,2^{n+1}r)\setminus B(x,2^n r)$,
\begin{align}\label{e:11}
&p(t,x,y) \nn\\
&\le \frac{c_6t}{V(x,\rho_n) \psi(\rho_n)}\Big(1+\frac{t}{\Phi(\rho_n)}\Big)^{d_2+1} + \frac{c_7}{V(x,\Phi^{-1}(t))}\Big(1+\frac{2^{n+1}r}{\Phi^{-1}(t)}\Big)^{d_2} \exp\Big(C_1 \frac{t}{\Phi(\rho_n)} -C_2 \frac{2^n r}{\rho_n}\Big) \nn \\ 
&\le \frac{c_6t}{V(x,\rho_n) \psi(\rho_n)}\Big(1+\frac{t}{\Phi(\rho_n)}\Big)^{d_2+1}
+ \frac{c_7}{V(x,\Phi^{-1}(t))} \Big(1+\frac{2^{n+1}r}{\Phi^{-1}(t)}\Big)^{d_2} \exp \Big( -c_4\Phi_1(r,t) - c_5  \frac{2^{n(1-\alpha)}r}{\Phi^{-1}(t)} \Big) \nn \\
&\le \frac{c_6t}{V(x,\rho_n) \psi(\rho_n)}\Big(1+\frac{t}{\Phi(\rho_n)}\Big)^{d_2+1}
+ \frac{c_7}{V(x,\Phi^{-1}(t))}  \exp \Big( -c_4\Phi_1(r,t) - \frac{c_5}{2}  \frac{2^{n(1-\alpha)}r}{\Phi^{-1}(t)} \Big),
\end{align}
where for the last inequality we have used the fact that $\frac{r}{\Phi^{-1}(t)}>C_0 \frac{\Phi^{-1}(C_1t)}{\Phi^{-1}(t)}  \ge 1$ and
\begin{align*}
 \sup_{n \in \N} \sup_{s>1} (1+2^{n+1} s)^{d_2} \exp\Big[-\frac{c_5}{2} 2^{n(1-\alpha)}s \Big]
  \le \sup_{s_1 >1} (1+2s_1)^{d_2} \exp\Big( -\frac{c_5}{2} s_1^{1-\alpha} \Big)  < \infty.
\end{align*}
With estimates in \eqref{e:11}, we get that
\begin{align*}
\int_{B(x,r)^c} p(t,x,y)\mu(dy) &\le \sum_{n=0}^\infty \int_{B(x,2^{n+1}r) \setminus B(x,2^n r)} p(t,x,y) \mu(dy) \\
&\le  c_8\sum_{n=0}^{\infty} \frac{V(x,2^n r)}{V(x,\Phi^{-1}(t))} \exp \left( -c_4 \Phi_1(r,t) - \frac{c_5}{2}\frac{2^{n(1-\alpha)}r}{\Phi^{-1}(t)}  \right)\nn\\
&\qquad+ c_8\sum_{n=0}^\infty  \frac{tV(x,2^nr)}{V(x,\rho_n) \psi(\rho_n)}\left(1+\frac{t}{\Phi(\rho_n)}\right)^{d_2+1}:= c_8(I_1+I_2).
\end{align*}
Using $r> C_0\Phi^{-1}(C_1t) \ge \Phi^{-1}(t)$, 
we obtain upper bound of $I_1$ by following the calculations in \eqref{e:wq1}. Thus,  
we have $
I_1 \le c_9 \exp\left( -a_2 \Phi_1(r,t) \right).
$

Next, we estimate $I_2$. 
Since $r= C_0 \Phi^{-1}(C_1t)^{1+\theta_0}/\psi^{-1}(C_1t)^{\theta_0}$, $\psi^{-1}(t)<\Phi^{-1}(t)$ and $\theta_0\le \theta<1/\delta_1$, we obtain
\begin{align*}
\frac{\Phi^{-1}(C_1t)}{\rho}\le \frac{\Phi^{-1}(C_1t) b^{-1}(C_2r/2)^{\delta_1}}{\Phi^{-1}(C_1t)^{\delta_1+1}}= b^{-1}(C_0C_2/2)^{\delta_1}\left(\frac{\Phi^{-1}(C_1t)}{\psi^{-1}(C_1t)}\right)^{\delta_1\theta_0}\le  b^{-1}(C_0C_2/2)^{\delta_1}\frac{\Phi^{-1}(C_1t)}{\psi^{-1}(C_1t)}. 
\end{align*}
Thus, $\frac{b \psi^{-1}(C_1 t)}{(C_0C_2/2)^{\delta_1}} \le \rho \le 2\Phi^{-1}(C_1t) \le r$. Using this, $\VD(d_2)$, \eqref{e:ginv}, $L(\beta_1,C_L, \psi)$ and $U(\alpha_2,c_U,\Phi)$ we have 
\begin{align*}
I_2 & =\sum_{n=0}^\infty \frac{V(x,2^nr)}{V(x,\rho_n)}\frac{t}{\psi(\rho_n)} \left(1+\frac{t}{\Phi(\rho_n)}\right)^{d_2+1} \\
&\le c_{10} \sum_{n=0}^\infty \frac{V(x,2^nr)}{V(x,\rho_n)}\frac{C_1 t}{\psi(b^{-1}(C_0C_2/2)^{\delta_1}\rho_n)} \frac{\psi(b^{-1}(C_0C_2/2)^{\delta_1}\rho_n)}{\psi(\rho_n)} \left(\frac{C_1 t}{\Phi(\rho/2)}\right)^{d_2+1} \\
&\le c_{11} \sum_{n=0}^\infty \Big( \frac{2^n r}{\rho_n} \Big)^{d_2} \Big( \frac{\psi^{-1}(C_1 t)}{\rho_n} \Big)^{\beta_1} \Big( \frac{\Phi^{-1}(C_1 t)}{\rho/2} \Big)^{\alpha_2(d_2 +1)} \\
&\le c_{12} \sum_{n=0}^\infty 2^{n(d_2 - \alpha(d_2+\beta_1))} \Big(\frac{r}{\rho}\Big)^{d_2} \Big( \frac{\psi^{-1}(C_1 t)}{\rho} \Big)^{\beta_1}\Big( \frac{\Phi^{-1}(C_1 t)}{\rho}\Big)^{\alpha_2(d_2+1)} \\
&= c_{13}r^{d_2} \psi^{-1}(C_1t)^{\beta_1}\Phi^{-1}(C_1t)^{\alpha_2(d_2+1)} \rho^{-d_2-\beta_1-\alpha_2(d_2+1)}.
\end{align*}
Since $br^{-\delta_1} \Phi^{-1}(C_1t)^{\delta_1+1} \le \rho$, we conclude that
\begin{equation}\label{e:I2}
I_2 \le c_{13} b^{-d_2-\beta_1-\alpha_2(d_2+1)} \Big( \frac{\psi^{-1}(C_1 t)}{r} \Big)^{\beta_1} \Big(\frac{\Phi^{-1}(C_1 t)}{r} \Big)^{-\delta_1(d_2\alpha_2+\alpha_2+\beta_1+d_2)-(d_2+\beta_1)}.
\end{equation}
Using $r=  C_0\Phi^{-1}(C_1t)^{1+\theta_0}/\psi^{-1}(C_1t)^{\theta_0}$, we have $ C_0\psi^{-1}(C_1t)< C_0\Phi^{-1}(C_1t) <r$. Since $\theta_0 \le \theta$, we have
$$\frac{ C_0\Phi^{-1}(C_1t)}{r} = \left( \frac{\psi^{-1}(C_1t)}{\Phi^{-1}(C_1t)} \right)^{\theta_0} = \left( C_0\frac{\psi^{-1}(C_1t)}{r} \right)^{\theta_0/(1+\theta_0)} \ge \left( C_0 \frac{\psi^{-1}(C_1t)}{r} \right)^{\theta/(1+\theta)}.$$
By using $\frac{\theta}{1+\theta}\big[\delta_1(d_2\alpha_2+\alpha_2+\beta_1+d_2) +(d_2+\beta_1)  \big] =\beta_1/2$, we have
\begin{align*} 
\left( \frac{\Phi^{-1}(C_1t)}{r} \right)^{-\delta_1(d_2\alpha_2+\alpha_2+\beta_1+d_2)-(d_2+\beta_1)  }
 \le c_{14}\left( \frac{\psi^{-1}(C_1t)}{r} \right)^{-\beta_1/2}.
\end{align*}
Therefore, using \eqref{e:I2} we obtain 
$I_2 \le c_{15}\left({\psi^{-1}(t)}/{r}\right)^{\beta_1/2}.$
By the estimates of $I_1$ and $I_2$, we arrive
\begin{equation*}\label{e:2p3}
\int_{B(x,r)^c} p(t,x,y)dy \le c_7(I_1 + I_2) \le c_{9}\bigg(\frac{\psi^{-1}(t)}{r}\bigg)^{\beta_1/2}+c_{15}\exp{\big(-a_2\Phi_1(r,t)\big)}. 
\end{equation*}
Combining all the cases, we obtain \eqref{e:tail-F}. Thus the assertions of Lemma \ref{l:tail} holds with $f(r,t):=\Phi_1(r,t)$. Thus, using Lemma \ref{l:tail} we have constants $k, c_0>0$ such that 
\begin{equation}\label{e:UHK}
p(t,x,y) \le \frac{c_0\, t}{V(x,d(x,y))\psi(d(x,y))} + \frac{c_0}{V(x,\Phi^{-1}(t))}\left(1+\frac{d(x,y)}{\Phi^{-1}(t)}\right)^{d_2} \exp{\big(-a_2k \Phi_1(d(x,y)/(16k), t) \big)} 
\end{equation}
for all $t\in (0,T]$ and $x,y \in M_0$. Recall that $2c_U>0$ is the constant in Lemma \ref{l:sT} with $\phi=\Phi$. When $d(x,y) \le 32c_Uk\Phi^{-1}(t)$, using $\UHKD(\Phi)$ and \eqref{VD2} we have
$$ p(t,x,y) \le p(t,x,x)^{1/2}p(t,y,y)^{1/2} \le \frac{c}{V(x,\Phi^{-1}(t))^{1/2}V(y,\Phi^{-1}(t))^{1/2}} \le \frac{c_{16}}{V(x,\Phi^{-1}(t))}. $$
Thus, by \eqref{e:sT2} and $d(x,y)/16k \le 2c_U\Phi^{-1}(t)$ we have
\begin{equation*}\label{e:UHK1}
p(t,x,y) \le \frac{c_{16}}{V(x,\Phi^{-1}(t))} \le \frac{c_{16}e^{a_2 c_{17} k}}{V(x,\Phi^{-1}(t))}\exp{\big(-a_2k \Phi_1(d(x,y)/(16k), t) \big)},
\end{equation*}
which yields \eqref{e:main10} for the case $d(x,y) \le 32c_Uk\Phi^{-1}(t)$. Also, for $r > 32c_Uk\Phi^{-1}(t)$ with $0<t \le T$, using \eqref{e:A8.1} with $n=0$ and \eqref{e:sT} we have
$$  c_4 \Phi_1(r,t) + c_5 \frac{r}{\Phi^{-1}(t) } \le \Phi_1(C_2 r, C_1 t) \le c_{17}\Phi_1(r/16k , t). $$
Therefore, using \eqref{e:UHK} we obtain
\begin{align*}
&p(t,x,y) \\ &\le \frac{c_0\, t}{V(x,d(x,y))\psi(d(x,y))} + \frac{c_0}{V(x,\Phi^{-1}(t))}\Big(1+\frac{d(x,y)}{\Phi^{-1}(t)}\Big)^{d_2} \exp{\big(-a_2k \Phi_1(d(x,y)/(16k), t) \big)}  \\
&\le \frac{c_0\, t}{V(x,d(x,y))\psi(d(x,y))} + \frac{c_0}{V(x,\Phi^{-1}(t))}\Big(1+\frac{d(x,y)}{\Phi^{-1}(t)}\Big)^{d_2} \exp{\Big(-a_3 \Phi_1(d(x,y), t) - c_{18}\frac{d(x,y)}{\Phi^{-1}(t)} \Big)} \\
&\le  \frac{c_0\, t}{V(x,d(x,y))\psi(d(x,y))} + \frac{c_{19}}{V(x,\Phi^{-1}(t))}\exp{\big(-a_3 \Phi_1(d(x,y), t) \big)},
\end{align*}
where we have used $\sup_{s>0} (1+s)^{d_2} \exp(-c_{18}s)<\infty$ for the last line. 
Combining two cases and $\UHKD(\Phi)$, we obtain \eqref{e:main10}.

\noindent (ii) Again we will show that there exist $a_1>0$ and $c_1>0$   such that for any $t\ge T$, $x \in M_0$ and $r>0$,
\begin{equation}\label{e:tail-F2}
\int_{B(x,r)^c} p(t, x,y)\,\mu(dy)\le c_1\bigg(\frac{\psi^{-1}(t)}{r}\bigg)^{\beta_1/2}+c_1\exp\left(-a_1\wt\Phi_1(r,t)\right).
\end{equation}
Note that using \eqref{e:wtcomp}, the proof of \eqref{e:tail-F2} for the case $r \le  C_0\Phi^{-1}(C_1t)$ and $r >  C_0\frac{\Phi^{-1}(C_1t)^{1+\theta}}{\psi^{-1}(C_1t)^\theta}$ are the same as that for (i). 

Without loss of generality we may assume $a = \Phi(T)$. Then for $t \ge T$, we have $\Phi^{-1}(t) = \wt \Phi^{-1}(t)$. Applying this and \eqref{e:wtcomp} for Lemma \ref{l:4.241}, we have for any $t \ge T$,
\begin{align*}
 p(t,x,y) 
 &\le \frac{c_1}{V(x,\Phi^{-1}(t))} \Big(1+\frac{d(x,y)}{ \Phi^{-1}(t)} \Big)^{d_2} \exp \Big( C_1 \frac{t}{ \Phi(\rho)} - C_2 \frac{d(x,y)}{\rho} \Big) + \frac{c_1 t}{V(x,\rho)\psi(\rho)} \Big( 1+ \frac{t}{\Phi(\rho)} \Big)^{d_2+1}         \\
 &\le \frac{c_1}{V(x,\wt\Phi^{-1}(t))} \Big(1+\frac{d(x,y)}{\wt \Phi^{-1}(t)} \Big)^{d_2} \exp \Big( C_1 \frac{t}{\wt \Phi(\rho)} - C_2 \frac{d(x,y)}{\rho} \Big) + \frac{c_1 t}{V(x,\rho)\psi(\rho)} \Big( 1+ \frac{t}{\wt\Phi(\rho)} \Big)^{d_2+1}. 
\end{align*}
Since $L(\delta, \wt C_L,\wt \Phi)$ holds, following the proof of (i) we have for any $t>0$ and $r>0$,
$$ \int_{B(x,r)^c} p(t,x,y) \mu(dy) \le c_1 \Big( \frac{\psi^{-1}(t)}{r} \Big)^{\beta_1 /2} + c_1 \exp(-a_1 \wt \Phi_1(r,t)). $$
 Since the assumptions in Lemma \ref{l:tail} follows from \eqref{e:sT2} and the fact that $\Phi^{-1}(t) = \wt\Phi^{-1}(t)$ for $t \ge T$, we obtain that for any $t \ge T$ and $x,y \in M$,
 \begin{equation*}\label{e:UHK3}
  p(t,x,y) \le \frac{c_0 t}{V(x,d(x,y))\psi(d(x,y))} + \frac{c_0}{V(x,\Phi^{-1}(t))} \Big( 1+ \frac{d(x,y)}{\Phi^{-1}(t)} \Big)^{d_2} \exp \big( - a_U \wt\Phi_1(d(x,y),t) \big). 
 \end{equation*}
  Here in the last term we have used \eqref{e:sT}. With the aid of $L(\delta,\wt C_L,\wt \Phi)$, the remainder is same as the proof of (i).  \qed

Now we give the proof of Theorem \ref{t:main11} and Corollary \ref{c:main11}. 

\smallskip

\noindent \textit{Proof of Theorem \ref{t:main11}.} 
By \eqref{comp1}, $\J_{\psi,\le}$ implies $\J_{\Phi,\le}$. Thus, \cite[Theorem 1.15]{HKE} yields that (2) implies (3) (where we have used $\RVD(d_1)$), and $(3)$ implies the conservativeness of $(\EE,\FF)$. Thus, by Theorem \ref{t:main10}, (3) implies (1). It remains to prove that (1) implies (2). By \eqref{comp1} and Remark \ref{r:hk}, $\UHK(\Phi,\psi)$ implies $\UHK(\Phi)$. Also, following the proof of \cite[Proposition 3.3]{HKE}, we easily prove that $\UHK(\Phi,\psi)$ also implies $\J_{\psi,\le}$. \qed

\noindent\textit{Proof of Corollary \ref{c:main11}.} 
Since $\J_{\psi,\le}$ implies $\J_{\Phi,\le}$, 
 \cite[Theorem 1.15]{HKE} implies 
the equivalence between the conditions in  Theorem \ref{t:main11}(3) and  Corollary \ref{c:main11}(4) and (5).
We now  prove the equivalence between the conditions in  Theorem \ref{t:main11}(3) and  Corollary \ref{c:main11}(6). To do this, we will use the results in \cite{HKE, GHH18}.

Suppose that  $\J_{\psi,\le}$, $\UHKD(\Phi)$ and $\E_\Phi$ hold. By \cite[Proposition 7.6]{HKE}, we have $\FK(\Phi)$. Since we have $\E_\Phi$, the condition $\textrm{EP}_{\Phi,\le, \eps}$ in \cite[Definition 1.10]{HKE} holds by \cite[Lemma 4.16]{HKE}. Since $\textrm{EP}_{\Phi,\le, \eps}$ implies the condition $(S)$ in \cite[Definition 2.7]{GHH18} with $r<\infty$ and $t<\delta \Phi(r)$, we can follow the proof of \cite[Lemma 2.8]{GHH18} line by line (replacing $r^{\beta}$ to $\Phi(r)$) and obtain $\Gcap(\Phi)$. 

Now, suppose that $\FK(\Phi)$, $\J_{\psi, \le}$ and $\Gcap(\Phi)$ hold. Then, by \cite[Lemma 4.14]{HKE}, we have $\E_{\Phi,\le}$. To obtain $\E_{\Phi,\ge}$, we first show that \cite[Lemma 4.15]{HKE} holds under our conditions. i.e., by using $\Gcap(\Phi)$ instead of $\CSJ(\Phi)$, we derive the same result in \cite[Lemma 4.15]{HKE}. To show \cite[Lemma 4.15]{HKE}, we give the main steps of the proof only.
Recall that for any $\rho>0$, $(\sE^{\rho}, \sF)$ is $\rho$-truncated Dirichlet form.  
For $\rho$-truncated Dirichlet form, we say $\textrm{AB}^{\rho}_{\zeta}(\Phi)$ holds if the inequality \cite[(2.1)]{GHH18} holds with $R'<\infty$, $\Phi(r\wedge\rho)$  and $J(x,y){\bf 1}_{\{d(x,y)<\rho\}}$ instead of $R'<\overline R$, $r^{\beta}$ and $j$ respectively. Similarly, we define $\textrm{AB}_{\zeta}(\Phi)=\textrm{AB}^{\infty}_{\zeta}(\Phi)$.
Then, by $\VD(d_2)$, $\J_{\psi,\le}$, \cite[Lemma 2.1]{HKE} and \eqref{comp1}, we can follow the proof of \cite[Lemma 2.4]{GHH18} line by line (replacing $r^{\beta}$ to $\Phi(r)$) to obtain $\textrm{AB}_{\zeta}(\Phi)$. To get
$\textrm{AB}_{1/8}(\Phi)$, we use the proof of \cite[Lemma 2.9]{GHH18}. Here, we take different $r_n, s_n, b_n, a_n$ from the one in the proof of \cite[Lemma 2.9]{GHH18}. Let $\lambda>0$ be a constant which will be chosen later. Take $s_n=c r e^{-n\lambda/2\alpha_2}$ for $n\ge1$, where $c=c(\lambda)$ is chosen so that $\sum^{\infty}_{n=1}s_n=r$ and $\alpha_2$ is the global upper scaling index of $\Phi$. Let $r_n=\sum^{n}_{k=1}s_k$ for $n\ge1$ and $r_0=0$. We also take $b_n=e^{-n\lambda}$ for $n\ge0$ and $a_n=b_{n-1}-b_n$ for $n\ge1$. (c.f. \cite[Proposition 2.4]{HKE}.) With these $r_n, s_n, b_n, a_n$, we can follow the proof of \cite[Lemma 2.9]{GHH18} line by line and obtain $\textrm{AB}_{1/8}(\Phi)$ by choosing small $\lambda>0$. Moreover,  using the argument in the proof of \cite[Proposition 2.3]{HKE}, we also obtain $\textrm{AB}^{\rho}_{1/8}(\Phi)$ which yields \cite[(4.8)]{HKE} for $\rho$-truncated Dirichlet form. Thus, we get  \cite[Corollary 4.12]{HKE}. For open subsets $U, B$ of $M$ with $U\subset B$, and for any $\rho>0$, define $\text{Cap}^{\rho}(U, B)=\inf\{\sE^{\rho}(\varphi,\varphi):\varphi\in \text{cutoff}(U,B)\}$. By $\Gcap(\Phi)$ with $u=1$(c.f. \cite[Definition 1.13]{GHH18} and below), we have 
$$\text{Cap}^{\rho}(B(x,R), B(x,R+r))\le \text{Cap}(B(x,R), B(x,R+r))\le  c\frac{V(x,R+r)}{\Phi(r)}\le c\frac{V(x,R+r)}{\Phi(r\wedge\rho)},$$
which implies the inequalities in \cite[Proposition 2.3(5)]{HKE}. Having this and \cite[Corollary 4.12]{HKE} at hand, we can follow the proof and get the result of \cite[Lemma 4.15]{HKE}. Now $\E_{\Phi,\ge}$ follows from the proof of \cite[Lemma 4.17]{HKE}. Since we now have $\E_{\Phi}$,  $\UHKD(\Phi)$ holds by \cite[Theorem 4.25]{HKE}. \qed

\subsection{Proofs of Theorems \ref{t:main21} and \ref{t:main12}}\label{s:ODLE}

Throughout this subsection, we will assume that all conditions of Theorems \ref{t:main21} and \ref{t:main12} hold. 
From $\J_\psi$ and $\VD(d_2)$, we immediately see that 
there is a constant $c>0$ such that for all $x, y\in M_0$ with $x\not= y$,
\begin{equation}\label{ujs}
J(x,y)\le  \frac{c}{V(x,r)}\int_{B(x,r)}J(z,y)\,\mu(dz)
\quad\hbox{for every }
0<r\le   d(x,y) /2.\tag{UJS}
\end{equation}
(See \cite[Lemma 2.1]{CKK}).

Recall that for any open set ${D} \subset M$, $\sF_D$ is $\sE_1$-closure in $\sF$ of $\sF\cap C_c(D)$. We use  $p^{D}(t,x,y)$ to denote the transition density function corresponding to the regular Dirichlet form $(\sE, \sF_{D})$.

Note  that $(\sE, \sF)$ is a conservative Dirichlet form by  \cite[Lemma 4.21]{HKE}. Thus, by  \cite[Theorem 1.15]{HKE}, we see that 
$\text{CSJ}(\Phi)$ defined in \cite{HKE} holds. Thus, using $\J_{\psi,\le}$, $\PI(\Phi)$, $\text{CSJ}(\Phi)$, $\eqref{ujs}$ and \eqref{comp1}, we have (7) in  \cite[Theorem 1.20]{PHI}.

Therefore, by \cite[Theorem 1.20]{PHI}, $\UHK(\Phi)$ and the following joint H\"older regularity hold for parabolic functions. We refer \cite[Definition 1.13]{PHI} for the definition of parabolic functions.
Note that, by a standard argument, we now can take the continuous version of  parabolic functions (for example, see \cite[Lemma 5.12]{GHH18}).
 Let $Q(t, x,r,R):=(t, t+r)\times B(x,R)$.

\begin{thm}			\label{t:PHR}
	There exist constants $c>0$, $0<\theta<1$ and $0<\epsilon<1$  such that for all $x_0\in M$, $t_0\geq 0$, $r>0$ and for every bounded measurable function $u=u(t,x)$ that is parabolic in $Q(t_0,x_0,\Phi(r),r)$,  the following parabolic H\"older regularity holds:
	\begin{align*}
	|u(s,x)-u(t,y)| \leq c\left(\frac{\Phi^{-1}(|s-t|)+d(x,y)}{r} \right)^{\theta}\sup_{[t_0, t_0+\Phi(r)]\times M}|u|
	\end{align*}
	for every $s,t\in (t_0,t_0+\Phi(\epsilon r))$ and $x,y \in B(x_0,\epsilon r)$.
\end{thm}

Since $p^{D}(t,x,y)$ is parabolic,
from now on, we assume $\sN=\emptyset$ and  take  the joint continuous versions of $p(t,x,y)$ and $p^{D}(t,x,y)$. (c.f., \cite[Lemma 5.13]{GHH18}).)

Again,  by \cite[Theorem 1.20]{PHI} we have the interior near-diagonal lower bound of $p^{B}(t, x,y )$ and parabolic Harnack inequality. 

\begin{thm}\label{t:NDL} 
	There exist $\eps\in (0,1)$ and $c_1>0$ such that for any $x_0\in M$, $r>0$, $0<t\le \Phi(\eps r)$ and $B=B(x_0,r)$,
$
	p^{B}(t, x,y )\ge {c_1}{V(x_0, \Phi^{-1}(t))^{-1}}$ for all  $x ,y\in B(x_0,\eps\Phi^{-1}(t)).$
\end{thm}

The proof of the next result is similar to that of  \cite[Proposition 5.4]{HKE}
\begin{prop}\label{l:LHK-psi} 
	There exists $\eta>0$ and $C_3>0$ such that for any $t>0$,
	\begin{equation}\label{e:GLHK1}
	p(t,x,y) \ge C_3 V(x,\Phi^{-1}(t))^{-1}, \quad x,y \in M \mbox{ with } d(x,y) \le \eta\Phi^{-1}(t), 
	\end{equation}
	and
	\begin{equation*}\label{e:GLHK2}
p(t,x,y) \ge \frac{C_3t}{V(x,d(x,y))\psi(d(x,y))}, \quad x,y \in M \mbox{ with } d(x,y) \ge \eta\Phi^{-1}(t).
	\end{equation*}
\end{prop}

\pf The proof of the proposition is standard. For example, see \cite[Proposition 5.4]{HKE}. 

Let $\eta=\eps/2<1/2$ where $\eps$ is the constant in Theorem  \ref{t:NDL}. Then 
by Theorem \ref{t:NDL},
\begin{align}
\label{e:nbe1}
p(t,x,y) \ge  p^{B(x,\Phi^{-1}(t)/\eps)}(t,x,y) \ge \frac{c_0}{V(x, \Phi^{-1}(t))} \quad \text{ for all } d(x,y) \le \eta\Phi^{-1}(t). 
\end{align}
Note that in the beginning of this section we have mentioned that $\UHK(\Phi)$ holds under $\J_{\psi,\le}$, $\PI(\Phi)$ and $\E_\Phi$. Thus, by \cite[Lemma 2.7]{HKE} and $\UHK(\Phi)$, we have
$$ \P^x(\tau_{B(x,r)} \le t) \le \frac{c_1 t}{\Phi(r)}, \quad r>0,\,\,t>0, \,\, x \in M. $$
Let $\eta_1:= (C_L/2)^{1/\beta_1}\eta\in (0,\eta)$ so that $\eta \Phi^{-1}((1-b)t) \ge \eta_1 \Phi^{-1}(t) $ holds for all $b \in(0, 1/2]$.
Then choose 
$\lambda \le c_1^{-1}C_U^{-1} (2\eta_1/3)^{\beta_2}/2 <1/2$
small enough so that
$\frac{c_1 \lambda t}{\Phi(2\eta_1 \Phi^{-1}(t)/3)}  \le \lambda c_1C_U (2\eta_1/3)^{- \beta_2}\le 1/2$. 
Thus we have  $\lambda \in (0,1/2)$ and $\eta_1 \in (0,\eta)$ (independent of $t$) such that
\begin{align}
\label{e:nbe3}
\eta \Phi^{-1}((1-\lambda)t) \ge \eta_1 \Phi^{-1}(t),  \quad \text{for all } t>0,
\end{align}
and
\begin{align}
\label{e:nbe4}
\P^x(\tau_{B(x,2\eta_1 \Phi^{-1}(t)/3)} \le \lambda t) \le 1/2, \quad \text{for all } t>0 \text{ and } x \in M.\end{align}
For the remainder of the proof we assume that $d(x,y) \ge \eta\Phi^{-1}(t)$.
Since, using \eqref{e:nbe1} and \eqref{e:nbe3}, 
\begin{align*}
p(t,x,y) &\ge \int_{B(y,\eta \Phi^{-1}((1-\lambda)t))} p(\lambda t,x,z) p((1-\lambda)t,z,y)\mu(dz) \\
&\ge \inf_{z \in B(y,\eta \Phi^{-1}((1-\lambda)t))} p((1-\lambda)t,z,y) \int_{B(y,\eta \Phi^{-1}((1-\lambda)t))} p(\lambda t,x,z)\mu(dz) \\
&\ge \frac{c_0}{V(y,\Phi^{-1}(t))} \P^x ( X_{\lambda t} \in B(y,\eta_1 \Phi^{-1}(t))), 
\end{align*}
it suffices to prove 
\begin{align}
\label{e:nbe2}
\P^x ( X_{\lambda t} \in B(y,\eta_1 \Phi^{-1}(t)))  \ge c_2\frac{tV(y,\Phi^{-1}(t))}{V(x,d(x,y))\psi(d(x,y))}.
\end{align}
The proof of \eqref{e:nbe2} is the same as the proof (ii) of \cite[Proposition 5.4]{HKE}. We skip it
\qed

Recall that the chain condition $\Ch(A)$ was defined in Definition \ref{d:chain}.

\begin{lemma}\label{l:dlhk}
Assume that the metric measure space $(M,d)$ satisfies $\Ch(A)$ and that the function $\Phi$ satisfies $L_a(\delta, \wt C_L)$ with $a>0$ and $\delta>1$. Then, for any $T>0$ and $C>0$, there exists a constant $c_1>0$ such that for any $t \in (0,T)$ and $x,y \in M$ with $d(x,y) \le C\Phi^{-1}(t)$,
\begin{equation}\label{e:GLHK3}
p(t,x,y) \ge \frac{c_1}{V(x,\Phi^{-1}(t))}.
\end{equation}
In particular, if $L(\delta, \wt C_L,\Phi)$ holds, then we may take $T=\infty$.
\end{lemma}
\pf Without loss of generality we may and do assume $a = \Phi^{-1}(T)$. Fix $t>0$ and $x,y \in M$ with $d(x,y) \le C\Phi^{-1}(t)$. Let $N:= \lceil \big( \frac{3AC}{\eta} \big)^{\frac{\delta}{\delta-1}} \wt C_L^{-\frac{1}{\delta-1}} \rceil +1 \in \N$. Then, by $\Ch(A)$ there exists a sequence $\{z_k\}_{k=0}^{N}$ of the points in $M$ such that $z_0=x,z_{N}=y$ and $d(z_k,z_{k+1}) \le A \frac{d(x,y)}{N}$ for all $k=0, \dots, N-1$. Note that by Lemma \ref{lem:inverse} and $L_{\Phi^{-1}(T)}(\delta, \wt C_L,\Phi)$ we have $U_T(1/\delta, \wt C_L^{-1/\delta},\Phi^{-1})$. Using this and the definition of $N$, we have
$$ A \frac{d(x,y)}{N} \le \frac{AC\Phi^{-1}(t)}{N}\le \frac{AC}{N} \wt C_L^{-1/\delta} N^{1/\delta} \Phi^{-1}(t/N) \le \frac{\eta}{3}\Phi^{-1}(t/N). $$
For $k=1, \dots, N$, let $B_k:= B(z_{k},\eta \Phi^{-1}(t/N)/3)$. Then, for any $0 \le k \le N-1$, $\xi_k \in B_k$ and $\xi_{k+1} \in B_{k+1}$. So we have
$$ d(\xi_k, \xi_{k+1}) \le d(\xi_k,z_k) + d(z_k,z_{k+1}) + d(z_{k+1}, \xi_{k+1}) \le \frac{2Ad(x,y)}{N} + \frac{\eta \Phi^{-1}(t/N)}{3} \le \eta \Phi^{-1}(t/N). $$
Thus, by \eqref{e:GLHK1} and \eqref{VD2} with $\xi_{k+1} \in B_{k+1}$, we have for any $k=0,\dots, N$, $\xi_k \in B_k$ and $\xi_{k+1} \in B_{k+1}$,
\begin{align*}
p(t/N,\xi_{k},\xi_{k+1}) &\ge \frac{c_1}{V(\xi_{k+1},\Phi^{-1}(t/N))} \ge \frac{c_1C_\mu^{-1}}{V(z_{k+1},\Phi^{-1}(t/N))}\Big(\frac{\Phi^{-1}(t/N)}{d(z_{k+1},\xi_{k+1})+\Phi^{-1}(t/N)}\Big)^{d_2} \\
&\ge \frac{c_2}{V(z_{k+1},\Phi^{-1}(t/N))}.
\end{align*}
Using above estimates and $\VD(d_2)$, we conclude
\begin{align*}
p(t,x,y) &\ge \int_{B_1} \dots \int_{B_{N-1}} p(t/N,x,\xi_1)p(t/N,\xi_1,\xi_2) \dots p(t/N, \xi_{N-1},y) \mu(d\xi_1)\mu(d\xi_2) \dots \mu(d\xi_{N-1}) \\
&\ge \int_{B_1} \dots \int_{B_{N-1}} \prod_{k=0}^{N-1} \frac{c_2}{V(z_k,\Phi^{-1}(t/N))} \mu(d\xi_1)\mu(d\xi_2) \dots \mu(d\xi_{N-1}) \\
&\ge c_2^{N} \prod_{k=1}^{N-1} V(z_k,\eta \Phi^{-1}(t/N)/3) \prod_{k=0}^{N-1} V(z_k,\Phi^{-1}(t/N))^{-1} \\
&\ge c_2^{N} c_3^{N-1} \prod_{k=1}^{N-1} V(z_k,\Phi^{-1}(t/N)) \prod_{k=0}^{N-1} V(z_k,\Phi^{-1}(t/N))^{-1} \ge c_4 V(x,\Phi^{-1}(t))^{-1}.
\end{align*}
This proves the lemma. \qed

\begin{prop}\label{l:LHK-exp1}
Assume that the metric measure space $(M,d)$ satisfies $\Ch(A)$. 

\noindent
(i) Suppose that $L_a(\delta,\wt C_L,\Phi)$ holds with $\delta>1$. Then, for any $T \in (0,\infty)$, there exist constants $c>0$ and $a_L>0$ such that for any $x,y \in M$ and $t \in (0,T]$,
	\begin{align}\label{proplow0}
	p(t,x,y)\geq cV(x,\Phi^{-1}(t))^{-1}\exp \left(-a_L \Phi_1(d(x,y),t) \right).
	\end{align}
Moreover, if $L(\delta, \wt C_L,\Phi)$ holds, then \eqref{proplow0} holds for all $t \in (0,\infty)$.
	
	\noindent
(ii) Suppose that $L^a(\delta,\wt C_L,\Phi)$ holds with $\delta>1$. Then, for any $T \in (0,\infty)$, there exist constants $c>0$ and $a_L>0$ such that for any $x,y \in M$ and $t \ge T$,
	\begin{align}\label{proplowinf}
p(t,x,y)\geq cV(x,\Phi^{-1}(t))^{-1}\exp \big(-a_L \wt\Phi_1(d(x,y),t) \big).
\end{align}
\end{prop}
\pf (i) Without loss of generality we may and do assume that $a =  \Phi^{-1}(T)$. Note that by \eqref{e:GLHK3}, we have a constant $c_1>0$ such that for any $t \in (0,T]$ and $x,y \in M$ with $d(x,y) \le 2c_U\Phi^{-1}(t)$,
\begin{equation}\label{e:LHK-exp1}
 p(t,x,y) \ge \frac{c_1}{ V(x,\Phi^{-1}(t))}.
 \end{equation}
Note that if $t \in (0,T]$ and $d(x,y) \le 2c_U \Phi^{-1}(t)$, \eqref{proplow0} immediately follows from \eqref{e:LHK-exp1} since $\Phi_1(d(x,y),t) \ge 0$. Now we consider $x,y \in M$ and $t \in (0,T]$ with $d(x,y) > 2c_U \Phi^{-1}(t)$.
Let $r:= d(x,y)$ and $\theta :=  \frac{ \wt C_L c_U}{2A} \land 2$. Define
$$ \eps = \eps(t,r) := \inf \Big\{ s>0 : \frac{\Phi(s)}{s} \ge \theta \frac{t}{r}  \Big \}. $$
Note that by \eqref{e:ginv} and $\theta \le 2$, we have $$\frac{\Phi(\Phi^{-1}(t))}{\Phi^{-1}(t)} \ge  \frac{c_U^{-1}t}{(2c_U)^{-1} r} \ge \theta \frac{t}{r},$$ 
which implies $\eps(t,r) \le \Phi^{-1}(t)$. Also, using $\lim_{s \to 0} \frac{\Phi(s)}{s} =0$ we have $\eps(t,r)>0$. Observe that by the definition of $\eps$, we have a decreasing sequence $\{s_n\}$ converging to $\eps$ satisfying
$\frac{\Phi(s_n)}{s_n} \ge \theta \frac{t}{r}$ for all $n \in \N$. Using $U(\alpha_2,c_U,\Phi)$ we have
$$ c_U \Big(\frac{\eps}{s_n}\Big)^{\alpha_2 -1} \frac{\Phi(\eps)}{\eps} \ge \frac{\Phi(s_n)}{s_n} \ge \theta \frac{t}{r} \quad \mbox{for all} \,\, n \in \N $$
Letting $n \to \infty$ we obtain
\begin{equation}\label{e:proplow1}
 \frac{\theta t}{c_Ur} \le \frac{\Phi(\eps)}{\eps}.
\end{equation}
By a similar way, using $L_{\Phi^{-1}(T)}(\delta, \wt C_L,\Phi)$ and $\frac{\Phi(s)}{s} \le \frac{\theta t}{r}$ for any $s < \eps$ we have
\begin{equation}\label{e:proplow1-1}
\frac{\Phi(\eps)}{\eps} \le \frac{\theta t}{\wt C_L r}.
\end{equation}
Also, \eqref{e:proplow1} yields that
\begin{align*}
\Phi_1(2c_U r,\theta t) \ge \frac{2c_U r}{\eps} - \frac{\theta t}{\Phi(\eps)} \ge \frac{r}{\eps} \Big( 2c_U- \frac{\eps}{\Phi(\eps)} \frac{\theta t}{r} \Big) \ge \frac{c_U r}{\eps}.
\end{align*}
Thus, using Lemma \ref{l:sT}(i) with the fact that $r \ge 2c_U\Phi^{-1}(t)$, we have a constant $c_1>0$ satisfying
\begin{equation}\label{e:proplow2}
\frac{r}{\eps} \le c_U^{-1}  \Phi_1(2c_Ur, \theta t) \le c_1 \Phi_1(r,t). 
\end{equation}
Define $N= N(t,r) := \lceil \frac{3A r}{2c_U \eps} \rceil+1$. Since $r \ge 2c_U \Phi^{-1}(t) \ge 2c_U \eps$, we have $N \ge \lceil 3A\rceil +1 \ge 4$. Observe that by $\frac{3Ar}{2c_U \eps} \le N \le \frac{2Ar}{c_U \eps}$ and \eqref{e:proplow1-1} with $\theta \le \frac{ \wt C_L c_U}{2A}$,
\begin{align*}
\Phi\Big(\frac{3Ar}{2c_U N}\Big) \le \Phi(\eps) \le \frac{ \eps \theta t}{r \wt C_L } \le \frac{2A\theta}{c_U \wt C_L } \frac{t}{N} \le \frac{t}{N}.
\end{align*}
This implies $\frac{Ar}{N} \le \frac{2}{3}c_U \Phi^{-1}(\frac{t}{N})$. On the other hand, since $(M,d)$ satisfies $\Ch(A)$, we have a sequence $\{z_l\}_{l=0}^N$ of points in $M$ such that $z_0=x$, $z_N=y$ and $d(z_{l-1},z_l) \le A \frac{r}{N}$ for any $l\in\{1,\dots, N\}$. Thus, for any $\xi_l \in B(z_l, \frac{2}{3}c_U\Phi^{-1}(\frac{t}{N}))$ and $\xi_{l-1} \in B(z_{l-1}, \frac{2}{3}c_U\Phi^{-1}(\frac{t}{N}))$, we have
\begin{align*}d(\xi_l,\xi_{l-1})&\leq d(\xi_l,z_l) + d(z_l,z_{l-1})+d(z_{l-1},\xi_{l-1}) \\ &\leq \frac{2}{3}c_U\Phi^{-1}(t/N) + \frac{Ar}{N} +  \frac{2}{3}c_U\Phi^{-1}(t/N)  \le 2c_U \Phi^{-1}(t/N). \end{align*}
Therefore, using semigroup property, \eqref{e:proplow2} and \eqref{e:GLHK3} with $N \le \frac{2A r}{ c_U \eps}$ we have
\begin{align}
&p(t,x,y) \geq \int_{B(z_{N-1},\frac{\eta}{3}\Phi^{-1}(t/N))} \cdots\int_{B(z_{1},\frac{\eta}{3}\Phi^{-1}(t/N))}p(\tfrac{t}{N},x,\xi_1)\cdots p(\tfrac{t}{N},\xi_{N-1},y)\mu(d\xi_1)\cdots \mu(d\xi_{N-1})\nn \\
&\geq c_2^N \prod_{l=0}^{N-1} V(z_l,\Phi^{-1}(t/N))^{-1} \prod_{l=1}^{N-1} V(z_l,\Phi^{-1}(t/N)) = c_3c_4^N V(x,\Phi^{-1}(t/N))^{-1}\nn \\
&\geq c_3\left(\frac{c_4 C^{d_1}}{3^{d_1}}\right)^N V(x,\Phi^{-1}(t))^{-1}
\geq c_3V(x,\Phi^{-1}(t))^{-1}\exp\left({-c_5 N}\right)\label{3} \\
&\geq c_3V(x,\Phi^{-1}(t))^{-1}\exp(-c_6 \frac{r}{\eps}) \geq c_3V(x,\Phi^{-1}(t))^{-1}\exp\left({-c_7\Phi_1(r,t)}\right) . \nn
\end{align}
This concludes \eqref{proplow0}. Now assume that $\Phi$ satisfies $L(\delta, \wt C_L)$. Note that the case $d(x,y) \le 2c_U\Phi^{-1}(t)$ is same, since we have \eqref{e:LHK-exp1} for every $t>0$. Also, by the similar way we obtain $0<\eps(t,r) \le \Phi^{-1}(t)$, \eqref{e:proplow1} and \eqref{e:proplow1-1} for all $t \in (0,\infty)$ and $r > 2c_U\Phi^{-1}(t)$. Following the calculations in \eqref{3} again, we conclude \eqref{proplow0} for every $t>0$ and $x,y \in M$ with $d(x,y) > 2c_U\Phi^{-1}(t)$.

(ii) Without loss of generality, we may assume $a = \Phi(T)$. Then, it suffices to prove
$$ p(t,x,y) \ge cV(x,\wt \Phi^{-1}(t))^{-1} \exp(-a_L\wt \Phi_1(d(x,y),t)), \quad t \ge T,\,\, x,y \in M.$$
Indeed, $\Phi^{-1}(t) = \wt \Phi^{-1}(t)$ for $t \ge T.$ 
Note that for the proof of \eqref{proplow0} with $T=\infty$, we only used  near-diagonal estimate in \eqref{e:GLHK1} and $L(\delta, \wt C_L, \Phi)$ with semigroup property. Since $L(\delta, \wt C_L,\wt\Phi)$ holds, \eqref{proplowinf} follows from \eqref{e:GLHK1} and \eqref{e:wtcomp}.
\qed

\noindent\textit{Proof of Theorem \ref{t:main21}.} Combining 
Theorem \ref{t:main10} and Propositions \ref{l:LHK-psi} and \ref{l:LHK-exp1}, we obtain our desired result. \qed

\noindent\textit{Proof of Theorem \ref{t:main12}.}  First we assume (2). Using Theorem \ref{t:main11} we obtain $\UHK(\Phi,\psi)$. Also, by $\UHK(\Phi)$, $\J_{\psi,\le}$ and the conservativeness of $(\EE,\FF)$ with Theorem \ref{t:main11} we have $\E_\Phi$. Now, the lower bound of $\HK(\Phi,\psi)$ follows from Proposition \ref{l:LHK-psi}. Therefore, (2) implies both (1) and (3).

Now we assume (1). The implication (1) $\Rightarrow$ $\J_\psi$ is the same as that in the proof of Theorem \ref{t:main11}. Since $\UHK(\Phi)$ holds, using \cite[Theorem 1.20 (3) $\Rightarrow$ (7)]{PHI} we obtain $\PI(\Phi)$. The conservativeness follows from \cite[Proposition 3.1]{HKE}.

Applying \cite[Theorem 1.20]{PHI} and \cite[Lemma 4.21]{HKE}, respectively, we easily see that (3) with \eqref{comp1} implies $\UHK(\Phi)$ and the conservativeness of $(\EE,\FF)$.

 If we further assume $\Ch(A)$, Theorem \ref{t:main21} yields that (3) $\Rightarrow$ (4).  Also, (4) $\Rightarrow$ (1) is straightforward. \qed

\section{HKE and stability on metric measure space with sub-Gaussian estimates for diffusion process} \label{s:main3}
In this section, we consider a metric measure space having sub-Gaussian estimates for 
 diffusion process. We  aim to establish an  equivalence relation similar to Theorems \ref{t:main11} and \ref{t:main12} without assuming that the index $\delta$ in $L_a(\delta,\wt C_L,\Phi)$ is strictly bigger than 1. 
We Recall that if $\Diff(F)$ holds, then there exists a conservative symmetric diffusion process $Z = (Z_t)_{t \ge 0 }$ on $M$ such that the transition density $q(t,x,y)$ of $Z$ with respect to $\mu$ exists and satisfies the estimates in \eqref{hkediff}. \bk

{\it  Throughout this section, we assume $\VD(d_2)$ and $\Diff(F)$ for the metric measure space $(M,d,\mu)$, where  $F:(0,\infty)\to(0,\infty)$ is a strictly increasing function satisfying \eqref{e:intcon}, $L(\gamma_1,c_F^{-1})$ and $U(\gamma_2,c_F)$, that is, 
\begin{align}\label{wsF}
c_F^{-1}\left(\frac{R}{r}\right)^{\gamma_1}\le\frac{F(R)}{F(r)}\le c_F \left(\frac{R}{r}\right)^{\gamma_2}, \quad 0<r \le R.
\end{align}
}

Note that, by Lemma \ref{lem:inverse}, $F^{-1}$ satisfies $L(1/\gamma_2, c_F^{-1/\gamma_2})$ and $U(1/\gamma_1, c_F^{1/\gamma_1})$. { Recall that under $\Diff(F)$, the constant $\gamma_1$ is greater or equal to $2$. }
Define $
\Phi(r)= {F(r)}\big/{\int_0^{r} \frac{dF(s)}{\psi(s)} }
$ as \eqref{e:dPhi}. Recall that we always assume that $\psi:(0,\infty)\to(0,\infty)$ is a non-decreasing function satisfying $L(\beta_1,C_L)$ and $U(\beta_2,C_U)$.

 Since $\psi$ is non-decreasing and $\displaystyle \lim_{s \to 0} \psi(s) = 0$, we easily observe that 
\begin{equation*}\label{e:Phi1}
 \psi(r) = \frac{F(r)}{ \int_0^r \frac{dF(s)}{\psi(r)} } >   \frac{F(r)}{\int_0^r \frac{dF(s)}{\psi(s)}} = \Phi(r), \quad r>0, 
 \end{equation*}
 and
 \begin{equation}\label{e:Phi2}
  \frac{\Phi(R)}{\Phi(r)} = \frac{F(R)}{F(r)} \cdot \frac{\int_0^r \frac{dF(s)}{\psi(s)}}{\int_0^R \frac{dF(s)}{\psi(s)}} \le \frac{F(R)}{F(r)}, \quad 0<r \le R. 
  \end{equation}
{ Thus, \eqref{comp1} holds and $\Phi$ satisfies $U(\gamma_2, c_F)$. If $\beta_2<\gamma_1$, then for $r>0$,
\begin{align*}
\frac{\psi(r)}{F(r)}\int^{r}_{0}\frac{dF(s)}{\psi(s)}
\le c \frac{r^{\beta_2}}{F(r)}\sum^{\infty}_{k=0}\int^{2^{-k}r}_{2^{-k-1}r}\frac{dF(s)}{s^{\beta_2}}
\le c \frac{r^{\beta_2}}{F(r)}\sum^{\infty}_{k=0}\frac{F(2^{-k}r)}{2^{-(k+1)\beta_2}r^{\beta_2}}
\le c \sum^{\infty}_{k=0}{2^{-k(\gamma_1-\beta_2)}}\le c.
\end{align*}
This shows that $\psi\asymp \Phi$ if $\beta_2<\gamma_1$. 
In particular,  $U(\beta_2, c_U, \Phi)$ holds for some $c_U>0$  if $\beta_2< \gamma_1$.}

 Recall that $F_1 = \sT(F)$.  
Note that $F_1(r,t) \in (0,\infty)$ for every $r,t>0$ under \eqref{wsF}.
Here we record \cite[Lemma 3.19]{GT} for the next use. Since $F$ is strictly increasing and satisfying \eqref{wsF}, we have that for any $r,t>0$, 
\begin{equation}\label{Phi0lower}
F_1(r,t)\geq \left(\frac{F(r)}{t}\right)^{\frac{1}{\gamma_1-1}}\wedge \left(\frac{F(r)}{t}\right)^{\frac{1}{\gamma_2-1}} \ge \left(\frac{F(r)}{t}\right)^{\frac{1}{\gamma_2-1}} -1.\end{equation}

\begin{lemma}\label{l:Phiwsc}
	$\Phi$ is strictly increasing. Moreover, $L(\alpha_1,c_L,\Phi)$ holds for some $\alpha_1,c_L>0$.	
\end{lemma}

\pf 
Since $\psi$ is non-decreasing, we observe that for any $0 \le a<b$, 
$$ \frac{F(b) - F(a)}{\psi(b)} \le \int_a^b \frac{dF(s)}{\psi(s)} \le \frac{F(b) - F(a)}{\psi(a)}, $$
regarding $\frac{1}{\psi(0)} = \infty$. Thus, there exists $a_* \in (a,b)$ such that $\int_a^b \frac{dF(s)}{\psi(s)} = \frac{F(b)-F(a)}{\psi(a_*)}$. For any $r<R$, let $r_* \in (0,r)$ and $R_* \in (r,R)$ be the constants satisfying
$$ \int_0^r \frac{dF(s)}{\psi(s)} = \frac{F(r)}{\psi(r_*)} \quad \mbox{and} \int_r^R \frac{dF(s)}{\psi(s)} = \frac{F(R) - F(r)}{\psi(R_*)}. $$
Then, since $\psi$ is non-decreasing, $$\Phi(R) = \frac{F(R)}{\int_0^r \frac{dF(s)}{\psi(s)} + \int_r^R \frac{dF(s)}{\psi(s)} } = \frac{F(R)}{\frac{F(r)}{\psi(r_*)} + \frac{F(R)-F(r)}{\psi(R_*)}} \ge \frac{F(R)}{\frac{F(r)}{\psi(r_*)} + \frac{F(R)-F(r)}{\psi(r_*)}} = \psi(r_*) = \Phi(r) .$$
Thus, $\Phi$ is also non-decreasing. Now suppose that the equality of above inequality holds. Then, since $F(R) -F(r)>0$, we have $\psi(r_*) = \psi(R_*)$, which implies that $\psi(r_*) = \psi(r)$ since $\psi$ in non-decreasing. Thus, we conclude $\int_0^r \frac{dF(s)}{\psi(s)} = \frac{F(r)}{\psi(r)}$, which is contradiction since $\lim_{s \to 0} \psi(s) = 0$. Therefore, $\Phi$ is strictly increasing.

Using $L(\gamma_1,c_F^{-1},F)$ and $L(\beta_1,C_L,\psi)$, there is a constant $C>1$ such that 
\begin{equation}\label{e:C}
F(Cr) \ge 4F(r) \quad \mbox{and} \quad \psi(Cr) \ge 4\psi(r), \quad \mbox{all} \quad r>0.
\end{equation}
For $r>0$, let $r_1 \in (0,r)$, $r_2 \in (r,Cr)$, $r_3 \in (Cr,C^2 r)$ be the constants satisfying
$$ \int_0^r \frac{dF(s)}{\psi(s)} = \frac{F(r)}{\psi(r_1)}, \quad \int_r^{Cr} \frac{dF(s)}{\psi(s)} = \frac{F(Cr)- F(r)}{\psi(r_2)} \quad \mbox{and} \quad \int_{Cr}^{C^2r} \frac{dF(s)}{\psi(s)} = \frac{F(C^2 r) - F(Cr)}{\psi(r_3)}. $$
Then,
$\Phi(r) = \psi(r_1)$ and
$$ \Phi(C^2 r) = \frac{F(C^2r)}{\int_0^{C^2 r} \frac{dF(s)}{\psi(s)}} = \frac{F(C^2 r)}{\frac{F(r)}{\psi(r_1)} + \frac{F(Cr) - F(r)}{\psi(r_2)} + \frac{F(C^2 r) - F(Cr)}{\psi(r_3)}} \ge \frac{F(C^2 r)}{\frac{F(r)}{\psi(r_1)} + \frac{F(Cr) - F(r)}{\psi(r_1)} + \frac{F(C^2 r) - F(Cr)}{\psi(r_3)}}. $$
By \eqref{e:C} and the fact that $r_1 \le r \le Cr \le r_3$, we have
$ \frac{\psi(r_1)}{\psi(r_3)} \le \frac{1}{4}$ and $\frac{F(Cr)}{F(C^2 r)} \le \frac{1}{4}. $
Therefore, for any $r>0$ we have
\begin{align}\label{e:C1}
\frac{\Phi(C^2 r)}{\Phi(r)} \ge \frac{F(C^2 r)}{F(Cr) +  \frac{\psi(r_1)}{\psi(r_3)}(F(C^2 r)- F(Cr))  } \ge \frac{F(C^2 r)}{F(Cr) + \frac{1}{4}(F(C^2 r)- F(Cr))   } \ge 2.
\end{align}
Using \eqref{e:C1} we easily prove that $L(\alpha_1,c_L,\Phi)$ holds with $\alpha_1 = \frac{\log 2}{2\log C}$ and $c_L = \frac{1}{2}$. \qed

\subsection{$\PI(\Phi)$, $\E_{\Phi}$ and upper heat kernel estimate via subordinate diffusion processes}\label{s:Sub}
Let $\phi$ be the function defined by
\begin{align}\label{e:phi}
\phi(\lambda)=\int^{\infty}_{0}(1-e^{-\lambda t})\frac{d t}{t\psi(F^{-1}(t))}.
\end{align}
Note that by \eqref{e:intcon}, $L(\beta_1,C_L,\psi)$ and $U(\gamma_2,c_F,F)$,
\begin{align*}
\int^{\infty}_{0}(1\wedge t)\frac{d t}{t\psi(F^{-1}(t))} = \int_0^1  \frac{dF(s)}{\psi(s)} + \int_{F(1)}^\infty \frac{dt}{t\psi(F^{-1}(t))} <\infty.
\end{align*}
Thus, there exists a subordinator $S=(S_t, t>0)$ which is independent of $Z$ and whose Laplace exponent is $\phi$. Then, the process $Y$ defined by $Y_t:=Z_{S_t}$ is pure jump process whose jump kernel is given by 
\begin{align*}
\sJ_\psi(x,y)=\int^{\infty}_{0}q(t,x,y)\frac{1}{t\psi(F^{-1}(t))} dt.
\end{align*}
Also, the transition density $p^Y(t,x,y)$ of $Y$ can be written by
\begin{align}\label{hkform}
p^Y(t,x,y)=\int^{\infty}_{0}q(s,x,y)\bP(S_t\in ds).
\end{align}
With sub-Gaussian estimates \eqref{hkediff}, we obtain the following lemma.
\begin{lemma}\label{l:Y_psi}
$\sJ_{\psi}$ satisfies \eqref{e:J_psi}. In other words, $\J_\psi$ holds for $Y$. 
\end{lemma}

\pf Fix $x,y \in M$ and denote $r:= d(x,y)$. We first observe that $F_1(r,t) \ge 0$ for any $r,t>0$. Also, for $r>0$ and $t \le F(r)$, we have
$ F_1(r,t) = \sup_{s>0} \big[ \frac{r}{s} - \frac{t}{F(s)} \big] \ge \frac{r}{F^{-1}(t)} - \frac{t}{F(F^{-1}(t))} = \frac{r}{F^{-1}(t)} - 1. $
By \eqref{hkediff} and the inequality $F_1(r,t) \ge 0 \lor ( \frac{r}{F^{-1}(t)}-1 )$, we have
\begin{align*}
&\sJ_{\psi}(x,y)=\int^{\infty}_{0}\frac{1}{t\psi(F^{-1}(t))}q(t,x,y)dt \le \int^{\infty}_{0}\frac{1}{t\psi(F^{-1}(t))} \frac{c}{V(x,F^{-1}(t))} \exp\big(-a_0 F_1(r,t)\big)dt\\
&\le \int^{\infty}_{F(r)}\frac{c}{t\psi(F^{-1}(t))V(x,F^{-1}(t))}dt+\int^{F(r)}_{0}\frac{c}{t\psi(F^{-1}(t))V(x,F^{-1}(t))}\exp\Big(-\frac{a_0r}{F^{-1}(t)} +a_0 \Big)dt \\
& =:I+II.
\end{align*}
We first consider $I$. Using  $L(\beta_1,C_L,\psi)$ and $L(1/\gamma_2, c_F^{-1/\gamma_2}, F^{-1})$ we have
\begin{align*}
I&\le \int^{\infty}_{F(r)}\frac{c}{t\psi(F^{-1}(t))V(x,F^{-1}(t))}dt
\le \frac{c}{V(x,r)}\int^{\infty}_{F(r)}\frac{1}{t\psi(F^{-1}(t))}dt\nn\\
&\le  \frac{c}{V(x,r)\psi( r)}\int^{\infty}_{F(r)}\frac{\psi(F^{-1}(F(r)))}{t\psi(F^{-1}(t))}dt \le \frac{c_1}{V(x,r)\psi(r)}\int^{\infty}_{F(r)}\frac{1}{t}\left(\frac{F(r)}{t}\right)^{\beta_1/\gamma_2}dt = \frac{c_1}{V(x,r)\psi(r)} \frac{\gamma_2}{\beta_1}.
\end{align*}
Now we obtain the upper bound of $II$. Assume $t\le F(r)$. Then, by $U(\beta_2,C_U,\psi)$, $\VD(d_2)$ and $U(\gamma_2,c_F, F)$,
\begin{align}\begin{split}\label{e:psi1}
\frac{1}{t\psi(F^{-1}(t))V(x,F^{-1}(t))}&=\frac{1}{F(r)\psi(r)V(x,r)}\frac{F(r)}{t}\frac{\psi(r)}{\psi(F^{-1}(t))}\frac{V(x,r)}{V(x,F^{-1}(t))}\\
&\le \frac{c_2}{F(r)\psi(r)V(x,r)} \Big( \frac{r}{F^{-1}(t)} \Big)^{\gamma_2+ \beta_2+d_2}.
\end{split}\end{align}
Since the function $s\mapsto s^{\gamma_2+\beta_2+d_2}e^{-a_0s}$ is uniformly bounded on $[1,\infty)$, using \eqref{e:psi1} we have
\begin{align*}
II&\le  e^{a_0}\int^{F(r)}_{0}\frac{c_3}{F(r)\psi(r)V(x,r)} \Big( \frac{r}{F^{-1}(t)} \Big)^{\gamma_2 + \beta_2 + d_2} \exp\Big(-\frac{a_0r}{F^{-1}(t)}\Big)dt \le  \frac{c_4}{\psi(r)V(x,r)}.
\end{align*}
Thus, $\sJ_\psi(x,y) \le I + II \le \frac{c_5}{\psi(r)V(x,r)}$. For the lower bound, we use \eqref{hkediff}, $\VD(d_2)$ and the fact that $\psi$ is non-decreasing to obtain that
\begin{align*}
\sJ_{\psi}(x,y) \ge \int^{F(r)}_{F(r)/2}\frac{c^{-1}}{t\psi(F^{-1}(t))V(x,F^{-1}(t))}dt \ge  \frac{c^{-1}}{2\psi(r)V(x,r)}.
\end{align*}
Now the conclusion follows. \qed

\begin{lemma}\label{l:subup}
	There exists $c>0$ such that for any $\lambda >0$,
\begin{equation}\label{e:phicomp}
\frac{1}{2\Phi(F^{-1}(\lambda^{-1}))} \le \phi(\lambda) \le \frac{c}{\Phi(F^{-1}(\lambda^{-1}))}.
\end{equation} 
\end{lemma}
\pf
Using \eqref{e:phi} and \eqref{e:dPhi},
\begin{align*}
\phi(\lambda)=\int^{\infty}_{0}(1-e^{-\lambda t})\frac{d t}{t\psi(F^{-1}(t))} \geq \int^{1/\lambda}_{0}\frac{\lambda t}{2t\psi(F^{-1}(t))}d t = \int_0^{F^{-1}(1/\lambda)} \frac{\lambda dF(s)}{2\psi(s)}  = \frac{1 }{2\Phi(F^{-1}(\lambda^{-1}))},
\end{align*}
and by \eqref{e:dPhi} and \eqref{comp1},
\begin{align*}
\phi(\lambda)=\int^{\infty}_{0}(1-e^{-\lambda t})\frac{d t}{t\psi(F^{-1}(t))} \leq \int^{1/\lambda}_{0}\frac{\lambda t}{t\psi(F^{-1}(t))}d t+\int^{\infty}_{1/\lambda}\frac{d t}{t\psi(F^{-1}(t))}\leq \frac{c}{\Phi(F^{-1}(\lambda^{-1}))}.
\end{align*}
From the above two inequalities we conclude the lemma.
\qed

{  For any open set $D\subset M$,  we define $\wt\tau_D:=\inf\{t>0: Z_t\in D^c\}$ and
$$q^{D}(t,x,y) :=q(t,x,y) - \bE^x[q(t-\wt\tau_{D},Z_{\wt\tau_{D}},y) : \wt\tau_{D} < t].$$
Then, by the strong Markov property, $q^{D}(t,x,y)$ is the transition density function of $Z^D$, the subprocess of $Z$ killed upon leaving $D$.} 
Following the proof of \cite[Lemma 2.3]{CKK09}, we obtain a consequence of $\Diff(F)$.

\begin{lemma}\label{l:Dlower}
	Assume that the metric measure space $(M,d,\mu)$ satisfies $\RVD(d_1)$, $\VD(d_2)$ and $\Diff(F)$. Then, there exist $ c_1, \theta>0$ and $\eps \in (0,1)$ such that
	\begin{equation}\label{e:Dlower}
	q^{B(x_0,r)}(t,x,y) \ge \frac{c_1}{V(x_0,F^{-1}(t))} \quad \mbox{for all} \quad x_0 \in M_0, \,r>0, \,\, x,y \in B(x_0,\eps F^{-1}(t)), \,\, t \in (0,\theta F(r)].
	\end{equation}
\end{lemma}
\pf Assume $\theta \in(0, 1]$. Let $x_0 \in M_0$, $r>0$ and denote $B_r:= B(x_0,r)$. Using \eqref{hkediff}, for any $x,y \in B(x_0,\eps F^{-1}(t)) $ and $t \in (0, \theta F(r)]$,
\begin{align*}
q^{B_r}(t,x,y) &=q(t,x,y) - \bE^x[q(t-\wt\tau_{B_r},Z_{\wt\tau_{B_r}},y) : \wt\tau_{B_r} < t] \\
&\ge \frac{c^{-1}}{V(x_0,F^{-1}(t))} - \bE^x\Big[ \frac{c}{V(x_0,F^{-1}(t-\wt\tau_{B_r}))} \exp\big( -a_0 F_1(d(Z_{\wt\tau_{B_r}},y),t-\wt\tau_{B_r})\big) : \wt\tau_{B_r} < t \Big] \\
&\ge \frac{c^{-1}}{V(x_0,F^{-1}(t))} - \bE^x\Big[ \frac{c}{V(x_0,F^{-1}(t-\wt\tau_{B_r}))} \exp\big( -a_0 F_1((1-\eps)r,t-\wt\tau_{B_r})\big) : \wt\tau_{B_r} < t \Big] \\
&\ge \frac{c^{-1}}{V(x_0,F^{-1}(t))}- c\,\P^x(\wt\tau_{B_r} < t) \sup_{0<s \le t} \frac{1}{V(x_0,F^{-1}(s))} \exp \big( -a_0 F_1((1-\eps)r,s)\big) \\
&\ge \frac{c^{-1}}{V(x_0,F^{-1}(t))}- c \sup_{0<s \le t} \frac{1}{V(x_0,F^{-1}(s))} \exp \big( -a_0 F_1((1-\eps)r,s)\big).
\end{align*}
By \eqref{Phi0lower}, we also have 
\begin{align*}
&\sup_{0<s \le t} \frac{1}{V(x_0,F^{-1}(s))} \exp \big( - a_0 F_1((1-\eps)r,s) \big) \\
&\le \sup_{0<s \le t} \frac{e^{a_0}}{V(x_0,F^{-1}(s))} \exp \Big(-a_0 \Big(\frac{F((1-\eps)r)}{s}\Big)^{\frac{1}{\gamma_2 -1}} \Big) \\
&=\sup_{0<s \le t}  \frac{e^{a_0}}{V(x_0,F^{-1}(t))} \frac{V(x_0,F^{-1}(t))}{V(x_0,F^{-1}(s))} \exp \left( -a_0 \Big(\frac{t}{s}\Big)^{\frac{1}{\gamma_2-1}} \Big( \frac{F((1-\eps)r)}{t} \Big)^{\frac{1}{\gamma_2 -1}} \right) \\
&\le \frac{e^{a_0}C_\mu c_F^{d_2/\gamma_1}}{V(x_0,F^{-1}(t))} \sup_{0<s \le t} (\frac{t}{s} )^{d_2/\gamma_1} \exp \left( -a_0c_F^{-1/(\gamma_2-1)} \Big(\frac{(1-\eps)^{\gamma_2}}{\theta}\Big)^{\frac{1}{\gamma_2 - 1}}  \Big(\frac{t}{s}\Big)^{\frac{1}{\gamma_2 -1}} \right) \\
&= \frac{e^{a_0}C_\mu c_F^{d_2/\gamma_1}}{V(x_0,F^{-1}(t))} \sup_{1 \le u} u^{d_2/\gamma_1} \exp \left( -a_0c_F^{-1/(\gamma_2-1)}\Big(\frac{(1-\eps)^{\gamma_2}}{\theta}\Big)^{\frac{1}{\gamma_2 - 1}}  u^{\frac{1}{\gamma_2 -1}} \right) := \frac{C(\theta)e^{a_0}C_\mu c_F^{d_2/\gamma_1}}{V(x_0,F^{-1}(t))}.
\end{align*}
Since $C(\theta) \to 0$ as $\theta \downarrow 0$, there exists a constant $\theta>0$ such that $C(\theta) \le \frac{1}{2c^2 e^{a_0} C_\mu c_F^{d_2/\gamma_1}}$. With this, we obtain
$$ q^{B_r}(t,x,y) \ge \frac{c^{-1}}{V(x_0,F^{-1}(t))}- c \sup_{0<s \le t} \frac{1}{V(x_0,F^{-1}(s))} \exp \big( -a_0 F_1((1-\eps)r,s)\big) \ge \frac{c^{-1}}{2V(x_0,F^{-1}(t))}. $$
This concludes the lemma. \qed

Let $ p^{Y,D}(t,x,y)$ be the transition density function of $Y^D$, the subprocess of $Y$ killed upon leaving $D$.

\begin{lemma}\label{l:PI}
{ Suppose that the metric measure space $(M,d,\mu)$ satisfies $\RVD(d_1)$, $\VD(d_2)$  and $\Diff(F)$ where $F:(0,\infty) \to (0,\infty)$ a strictly increasing function  satisfying $\eqref{e:intcon}$, $L(\gamma_1,c_F^{-1})$ and $U(\gamma_2,c_F)$.}
There exist constants $c>0$ and $\hat \eps \in (0,1)$ such that for any $x_0 \in M_0$, $r>0$, $0<t \le \Phi(\hat \eps r)$ and $x,y \in B(x_0,\hat \eps \Phi^{-1}(t))$,
$$ p^{Y,B(x_0,r)}(t,x,y) \ge \frac{c}{V(x_0,\Phi^{-1}(t))}. $$
\end{lemma}
\pf Recall that we have defined $Y_t=Z_{S_t}$, where $S_t$ is a subordinator independent of $Z$ and whose Laplace exponent is the function $\phi$ in \eqref{e:phi}. Also, by \eqref{e:phicomp} we have $\frac{c_1^{-1}}{\Phi(F^{-1}(\lambda^{-1}))} 
\le \phi(\lambda) \le \frac{c_1}{\Phi(F^{-1}(\lambda^{-1}))}$. Take $\lambda$ by $F(\Phi^{-1}(c_1^{-1}t^{-1}))^{-1}$ and $F(\Phi^{-1}(c_1t^{-1}))^{-1}$, and using the fact that $\Phi$ and $F$ are strictly increasing we obtain that for any $t>0$
\begin{equation}\label{e:PI1}
 F(\Phi^{-1}(c_1^{-1} t)) \le \phi^{-1}( t^{-1})^{-1} \le  F(\Phi^{-1}(c_1 t)). 
\end{equation}
By \cite[Proposition 2.4]{M}, there exist $\rho,c_2>0$ such that
\begin{equation}\label{e:ante} 
\P\Big( \frac{1}{2\phi^{-1}(t^{-1})} \le S_t \le \frac{1}{\phi^{-1}(\rho t^{-1})} \Big) \ge c_2.
\end{equation}
Choose $\hat \eps >0$ such that 
$$ \hat \eps \Phi^{-1}(t) \le \eps F^{-1}\Big(\frac12 F(\Phi^{-1}(c_1^{-1} t))\Big) \quad \mbox{and} \quad  F\big(\Phi^{-1}(c_1 \rho^{-1}\Phi(\hat \eps r))\big) \le \theta F( r), $$
where $\eps \in (0,1)$ and $\theta$ are the constants in Lemma \ref{l:Dlower}.
Then, by \eqref{e:PI1}, we see that for $0<t \le \Phi(\hat \eps r)$ and $s \in [ \frac{1}{2\phi^{-1}(t^{-1})}, \frac{1}{\phi^{-1}(\rho t^{-1})}]$, we have
$$ s \le \frac{1}{\phi^{-1}(\rho t^{-1})} \le F(\Phi^{-1}(c_1 \rho^{-1} t)) \le F(\Phi^{-1}(c_1 \rho^{-1}\Phi(\hat \eps r)))  \le \theta F(r)  $$
and
$$ \hat \eps \Phi^{-1}(t)  \le \eps F^{-1}\Big(\frac{1}{2} F(\Phi^{-1}(c_1^{-1} t))\Big)  \le \eps F^{-1} \Big( \frac{1}{2\phi^{-1}(t^{-1})}\Big)\le \eps F^{-1}(s). $$
Thus, by \cite[Proposition 3.1]{SV08}, Lemma \ref{l:Dlower}, \eqref{e:ante}, \eqref{e:PI1}, $\VD(d_2)$, \eqref{wsF} and $U(1/\alpha_1,c_L^{-1/\alpha_1}, \Phi^{-1})$, we see that for $0<t \le \Phi(\hat \eps r)$ and $x,y \in B(x_0, \hat \eps \Phi^{-1}(t))$
\begin{align*}
p^{Y,B(x_0,r)}(t,x,y) &\ge \int_0^\infty q^{B(x_0,r)}(s,x,y) \P(S_t \in ds) \\
& \ge \int_{\frac{1}{2\phi^{-1}(t^{-1}) }}^{ \frac{1}{\phi^{-1}(\rho t^{-1}) }} q^{B(x_0,r)}(s,x,y) \P(S_t \in ds) \\
& \ge  \frac{c_3}{V(x_0,F^{-1}(\phi^{-1}(\rho t^{-1})^{-1}))} \P\Big( \frac{1}{2\phi^{-1}(t^{-1})} \le S_t \le \frac{1}{\phi^{-1}(\rho t^{-1})} \Big) \\
& \ge \frac{c_2 c_3}{V(x_0,F^{-1}( F(\Phi^{-1}(c_1\rho^{-1}t))))} \ge \frac{c_4}{V(x_0,\Phi^{-1}(t))}.
\end{align*}
This finishes the proof of lemma. 
\qed

We now combine our Lemmas \ref{l:Y_psi} and \ref{l:PI} with \cite[Theorem 1.20]{PHI} and obtain a preliminary 
upper bound of $p(t,x,y)$.

\begin{thm}\label{t:UHK-Phi} 
{ Suppose that the metric measure space $(M,d,\mu)$ satisfies $\RVD(d_1)$, $\VD(d_2)$ and $\Diff(F)$ where $F:(0,\infty) \to (0,\infty)$ is a strictly increasing function  satisfying $\eqref{e:intcon}$, $L(\gamma_1,c_F^{-1})$ and $U(\gamma_2,c_F)$.}
Assume that $(\EE,\FF)$ is a regular Dirichlet form on $(M,d,\mu)$ satisfying $\J_\psi$. Then,	there exists a constant $c>0$ such that
	\begin{align}\label{UHK_metric} 
	p(t,x,y)\le c\left(\frac{1}{V(x,\Phi^{-1}(t))}\wedge \frac{t}{V(x,d(x,y))\Phi(d(x,y))}\right)
	\end{align}
	for all $t>0$ and $x,y\in M$. Moreover, $\E_{\Phi}$ and $\PI(\Phi)$ holds for $(\EE,\FF)$. 
\end{thm}
\pf Note that from Lemmas \ref{l:Y_psi} and \ref{l:PI}, the condition (4) in \cite[Theorem 1.20]{PHI} holds for the process $Y$. In particular, using \cite[Theorem 1.20]{PHI}, the conditions $\CSJ(\Phi)$ and $\PI(\Phi)$ holds for the process $Y$. Since the jump kernel of $X$ and $Y$ are comparable by Lemma \ref{l:Y_psi}, the conditions $\PI(\Phi)$ and $\CSJ(\Phi)$ also hold for $X$. In particular, the process $X$ satisfies condition (7) in \cite[Theorem 1.20]{PHI}. Now, using \cite[Theorem 1.20]{PHI} again we obtain $\E_\Phi$ and $\UHK(\Phi)$. This completes the proof. 
\qed

\subsection{Proofs of Theorem \ref{t:main32}} \label{s:main3proof}

In this section, we prove Theorem \ref{t:main32}. 
\begin{thm} \label{t:main31}
{	Assume that the metric measure space $(M,d,\mu)$ satisfies $\RVD(d_1)$ and $\VD(d_2)$. Also, assume further that $\Diff(F)$ holds for a strictly increasing function $F:(0,\infty) \to (0,\infty)$ satisfying $L(\gamma_1,c_F^{-1})$ and $U(\gamma_2,c_F)$. Suppose that the regular Dirichlet form $(\EE,\FF)$ satisfies $\J_\psi$, where $\psi$ is a non-decreasing function satisfying $L(\beta_1,C_L)$ and $U(\beta_2,C_U)$.} Then, there exist constants $c>0$ and $a_1>0$ such that for all $t>0$ and $x,y \in M$,
	\begin{equation*}\label{e:main31}
	p(t, x,y) \le \frac{c}{V(x, \Phi^{-1}(t))}\wedge \left(\frac{c\,t}{V(x, d(x,y))\psi(d(x,y))} + \frac{c}{V(x, \Phi^{-1}(t))} e^{-a_1 F_1( d(x,y), F(\Phi^{-1}(t))  ) }  \right).
	\end{equation*}
\end{thm}
\pf  Note that  $\gamma_1\ge2$ and \eqref{e:intcon} hold by the assumption. 
Since the proof is similar to that of Theorem \ref{t:main10}, we just give the difference. As in the proof of Theorem \ref{t:main10}, we will show that there exist $a_1>0$ and $c_1>0$  such that for any $t>0$ and $r>0$,
\begin{equation}\label{3tail}
\int_{B(x,r)^c} p(t, x,y)\,\mu(dy) \le c_1\bigg(\frac{\psi^{-1}(t)}{r}\bigg)^{\beta_1/2}+c_1\exp{\left(-a_1 F_1 \big(r,F(\Phi^{-1}(t))   \big)   \right)}.
\end{equation}
Let $\gamma := \frac{1}{\gamma_1 -1}$,  
$$\theta:=\frac{(\gamma_1-1)\beta_1}{\gamma_1(2d_2+\beta_1) + (\beta_1+2\beta_2+2d_2\beta_2)} \in (0, \gamma_1-1)$$ and $C_0 = \frac{4c_F}{C_2}$, where $C_1$ and $C_2$ are the constants in Lemma \ref{l:4.22}. Note that  we may and do assume that $C_1 \ge 1$ and $C_2 \le 1$ without loss of generality. Assume $r \le C_0\Phi^{-1}(C_1t)$. Then, using $L(\alpha_1,c_L,\Phi)$ and Lemma \ref{lem:inverse} we have $$ r \le C_0 \Phi^{-1}(C_1t) \le C_0 C_1^{1/\alpha_1} c_L^{-1/\alpha_1} \Phi^{-1}(t) \le c_1 F^{-1}(F(\Phi^{-1}(t))).$$
Thus, by \eqref{e:sT2}, there is a constant $c_2 \ge 0$ such that
$$ F_1\big(r, F(\Phi^{-1}(t)) \big)\le c_2 $$
for any $r,t>0$ with $r \le C_0\Phi^{-1}(t)$. Thus, we have
\begin{equation*}
\int_{B(x,r)^c} p(t,x,y) \mu(dy) \le 1 \le  e^{a_1c_2} \exp \big(-a_1 F_1(r,t)\big).
\end{equation*}
This implies \eqref{3tail} for $r \le C_0\Phi^{-1}(C_1t)$. Also, \eqref{3tail} for $r \ge C_0\Phi^{-1}(C_1t)^{1+\theta} / \psi^{-1}(C_1t)^{\theta}$ follows from Lemma \ref{l:exp1}.

Now consider the case $C_0\Phi^{-1}(C_1t) < r \le C_0\Phi^{-1}(C_1t)^{1+\theta} / \psi^{-1}(C_1t)^{\theta}$. Since \eqref{comp1} holds, there exists $\theta_0 \in (0,\theta]$ satisfying $r= C_0\Phi^{-1}(C_1t)^{1+\theta_0}/\psi^{-1}(C_1t)^{\theta_0}$. By Lemma \ref{l:max} and $C_0 \ge 2c_F>0$ with $L(\gamma_1,c_F^{-1},F)$, there is $\rho \in [b(\frac{C_2 r}{2c_F})^{-\gamma} \Phi^{-1}(C_1 t)^{\gamma+1}, 2\Phi^{-1}(C_1 t)]$ such that
$$ \frac{C_2r}{2 c_F\rho} - \frac{F(\Phi^{-1}(C_1 t))}{F(\rho)} + 1 \ge F_1 \Big(  \frac{C_2}{2c_F} r, F(\Phi^{-1}(C_1t))\Big). $$
Let us define $r_n= 2^nr$ and 
$\rho_n = 2^{n\alpha} \rho$ for $n \in \N_0$, with some $\alpha \in (\frac{d_2}{d_2+\beta_1},1)$. Since $b(\frac{C_2r}{2c_F})^{-\gamma} \Phi^{-1}(C_1t)^{\gamma+1}\le \rho  \le 2\Phi^{-1}(C_1t)$, using \eqref{e:ginv}, \eqref{e:Phi2} and \eqref{e:sT} we have
\begin{align}\begin{split}\label{e:A8.2}
\frac{C_1 t}{\Phi(\rho_n)} -  \frac{C_2 2^nr}{\rho_n}
&\le \Big(\frac{c_F\Phi(\Phi^{-1}(C_1t))}{\Phi(\rho)} - \frac{C_2  r}{2\rho}\Big) + \frac{C_2  r}{2\rho} - \frac{2^{n(1-\alpha)}C_2  r}{\rho} \\
&\le \Big(\frac{ c_F F(\Phi^{-1}(C_1t))}{ F(\rho)} - \frac{C_2 r}{2\rho}\Big) + \Big(\frac 12 -2^{n(1-\alpha)}\Big) \frac{C_2r}{\rho}  \\
&\le -c_F F_1( \frac{C_2 r}{2c_F} , F(\Phi^{-1}( C_1t))  )  + c_F + \Big(\frac 12 -2^{n(1-\alpha)}\Big) \frac{C_2r}{\rho}\\ 
&\le -c_3 F_1(r,F (\Phi^{-1}(t))) - \frac{ 1}{2} 2^{n(1-\alpha)} \frac{C_2 r}{2\Phi^{-1}(C_1t)} + c_F \\ 
&\le -c_3 F_1(r,F (\Phi^{-1}(t)))  - c_4 2^{n(1-\alpha)} \frac{r}{\Phi^{-1}(t)} + c_F,
\end{split}\end{align}
which is the counterpart of \eqref{e:A8.1}. Since we have Lemma \ref{l:4.241} and \eqref{e:A8.2}, following the proof of Theorem \ref{t:main10} we obtain \eqref{3tail}.
Thus, we can apply Lemma \ref{l:tail} with $f(r,t):=F_1 (r, F(\Phi^{-1}(t)))$. Note that condition (i) in Lemma \ref{l:tail}  follows from \eqref{e:sT2} and \eqref{e:ginv}. Thus,
$$ p(t,x,y) \le \frac{c_0 t}{V(x,d(x,y))\psi(d(x,y))} + \frac{c_0}{V(x,\Phi^{-1}(t))} \Big( 1+ \frac{d(x,y)}{\Phi^{-1}(t)}\Big)^{d_2} e^{- a_2 k F_1(d(x,y)/(16k), F(\Phi^{-1}(t)))}. $$
The remainder is the same as the proof of Theorem \ref{t:main10}. \qed

\noindent\textit{{Proof of Theorem \ref{t:main32}.}} 
 By \eqref{e:Phi2}, $(3) \Rightarrow (2)$ holds trivially. 

Let $\varphi$ be a non-decreasing function satisfying $U(\delta_2,c_0)$ for some $\delta_2>0$ and $c_0 \ge 1$. We first note that $\varphi^{-1}(t)\to0$ as $t\to0$ by $U(\delta_2,c_0, \varphi)$ and Lemma \ref{l:ginv}.  
Using \eqref{Phi0lower} and \eqref{wsF}, we have that for any $t,r>0$ with $\varphi^{-1}(t)\le r$,
\begin{align*}
F_1(r, F(\varphi^{-1}(t)))\ge \Big(\frac{F(r)}{F(\varphi^{-1}(t))}\Big)^{1/(\gamma_2-1)}-1\ge c_F^{-1/(\gamma_2-1)}\Big(\frac{r}{\varphi^{-1}(t)}\Big)^{\gamma_1/(\gamma_2-1)}-1.
\end{align*}
Thus, by Lemma \ref{l:ginv} and $\VD$, we have that for any $t,r>0$ with $\varphi^{-1}(t)\le r$,
\begin{align}\label{pf_t2.19}
&\frac{1}{tV(x,\varphi^{-1}(t))}\exp\big(-a_UF_1(r, F(\varphi^{-1}(t)))\big)\nn\\
&\le c_0 \frac{\varphi(r)}{\varphi(\varphi^{-1}(t))}\frac{V(x,r)}{V(x,\varphi^{-1}(t))}\frac{1}{\varphi(r)V(x,r)}\exp\big(-a_UF_1(r, F(\varphi^{-1}(t)))\big)\nn\\
&\le c_0^2C_{\mu}\frac{1}{\varphi(r)V(x,r)}\Big(\frac{r}{\varphi^{-1}(t)}\Big)^{\delta_2+d_2}\exp\Big(-a_Uc_F^{-1/(\gamma_2-1)}\Big(\frac{r}{\varphi^{-1}(t)}\Big)^{\gamma_1/(\gamma_2-1)}+a_U\Big). 
\end{align}
The last term in above inequalities tends to zero as $t\to0$.
Thus, under $\GHK(\varphi,\psi)$, we see that
\begin{align*}
&\limsup_{t \to 0} \frac{p(t,x,y)}{t} \\
&\le \limsup_{t \to 0}\frac{c}{tV(x, \varphi^{-1}(t))}\wedge\left(\frac{c}{V(x,d(x,y))\psi(d(x,y))}+\frac{c}{tV(x,\varphi^{-1}(t))}e^{-a_UF_1(d(x,y), F(\varphi^{-1}(t)))}\right)\\
&=\frac{c}{V(x,d(x,y))\psi(d(x,y))}.
\end{align*}
Moreover, under $\GHK(\varphi,\psi)$, it is easy to see that
$$\liminf_{t \to 0} \frac{p(t,x,y)}{t}\ge \frac{c^{-1}}{V(x,d(x,y))\psi(d(x,y))}.$$
Using the above observations and following the proof of \cite[Proposition 3.3]{HKE}, we obtain that $\GHK(\varphi,\psi)$ implies $\J_{\psi}$.  

Now assume $\J_{\psi}$. Then, the upper bound in $\GHK(\Phi,\psi)$ follows from Theorem \ref{t:main31}. Also, using Theorem \ref{t:UHK-Phi} and Lemma \ref{l:PI} we have $\E_{\Phi}$ and $\PI(\Phi)$. Since all conditions in Proposition \ref{l:LHK-psi} holds, we obtain the lower bound of $\GHK(\Phi,\psi)$. 

When $\Phi$ in \eqref{e:dPhi} satisfies $L(\alpha_1,c_L)$ with $\alpha_1>1$ and $(M,d)$ satisfies $\Ch(A)$ for some $A \ge 1$, by Proposition \ref{l:LHK-exp1} we conclude that $\GHK(\Phi,\psi)$ implies $\SHK(\Phi,\psi)$. Since $\SHK(\Phi,\psi)$ is clearly
stronger than $\GHK(\Phi,\psi)$, we finish the proof. 
\qed

\noindent\textit{Proof of Corollary \ref{c:main33}.}
Note that we already proved (i) in the previous proof.
Recall that by Lemma \ref{l:Phiwsc} and the observation below \eqref{e:Phi2}, we see that $L(\alpha_1, c_L, \Phi)$ and $U(\alpha_2, c_L, \Phi)$ hold.
Now we assume that $\GHK(\Phi,\psi)$ holds and  that $\Phi$ satisfies $U(\alpha_2, c_U)$ with $\alpha_2<d_1$ and
we show Green function estimates.  

Fix $x,y \in M_0$ and let $r=d(x,y)$.   Thus, by $\GHK(\Phi,\psi)$, Lemma \ref{lem:inverse} and $\VD(d_2)$, we get
\begin{align*}
G(x,y)&\geq \int_{\Phi(r/\eta)}^\infty p(t,x,y)dt \geq c_1\int_{\Phi(r/\eta)}^{2\Phi(r/\eta)} V(x,\Phi^{-1}(t))^{-1}dt\geq c_2\frac{\Phi(r)}{V(x,r)}.
\end{align*}
For the upper bound, we first note that by the change of variable and $\RVD(d_1)$, 
\begin{align*}
&\int_{\Phi(r)}^\infty\frac{dt}{ V(x,\Phi^{-1}(t))}  =\int_{r}^{\infty}\frac{d\Phi(s)}{V(x,s)} \le\frac{c_3}{V(x,r)}\int_{r}^{\infty}\left(\frac{r}{s}\right)^{d_1} d\Phi(s)=\frac{c_3r^{d_1}}{V(x,r)}\int_{r}^{\infty}s^{-d_1} d\Phi(s).
\end{align*}
Using  the integration by parts and $U(\alpha_2, c_U, \Phi)$,
\begin{align*}
\int_{r}^{\infty}s^{-d_1} d\Phi(s)&=\Big[s^{-d_1}\Phi(s)\Big]^{\infty}_{r}+d_1\int_{r}^{\infty}s^{-d_1-1}\Phi(s)ds\\
&\le \Big[s^{-d_1}\Phi(s)\Big]^{\infty}_{r}+c_4\Phi(r)\int_{r}^{\infty}s^{-d_1-1}\left(\frac sr\right)^{\alpha_2}ds.
\end{align*} 
Thus, by using the condition $d_1>\alpha_2$, we get that
$
\int_{\Phi(r)}^\infty V(x,\Phi^{-1}(t))^{-1} dt \le c_5 \Phi(r)V(x,r)^{-1}.
$
By this inequality and  Theorem \ref{t:UHK-Phi}, we conclude that
\begin{align*}
G(x,y)&=\int_0^{\Phi(r)}p(t,x,y)\,dt +\int_{\Phi(r)}^\infty p(t,x,y)\,dt\\
&\leq \frac{c_6}{V(x,r)\Phi(r)}
\int_0^{\Phi(r)} t\,dt+c_6\int_{\Phi(r)}^\infty V(x,\Phi^{-1}(t))^{-1}dt \le c_7\frac{\Phi(r)}{V(x,r)}.
\end{align*}
Since $V(x,r)\asymp V(y,r)$ by \eqref{VD2}, we also have $G(x,y)\asymp \Phi(r)V(y,r)^{-1}$.

\qed

\noindent\textit{Proof of Corollary \ref{c:main32}.} Using Corollary \ref{c:main33}(i), $X$ satisfies $\J_\psi$, $\PI(\Phi)$ and $\E_{\Phi}$. Thus, the conclusion follows from Theorem \ref{t:main21}. \qed

\section{Generalized Khintchine-type law of iterated logarithm at infinity}\label{s:LIL}

In this section, using  our main result, we establish a generalized version of 
 the law of iterated logarithm at infinity on metric measure space.
 Throughout this section, as in Corollary \ref{c:main32} we assume that $(M,d,\mu)$ satisfies $\RVD(d_1)$, $\VD(d_2)$ and $\Diff(F)$, where $F$ is a strictly increasing and satisfies \eqref{wsF}. 

Let $(\EE,\FF)$ be a regular Dirichlet form given by \eqref{e:DF}, which satisfies $\J_\psi$ with a non-decreasing function $\psi$ satisfying  $L(\beta_1,C_L)$ and $U(\beta_2,C_U)$. 
Recall that $X= \{ X_t , t \ge 0; \P^x , x \in M \}$ is the $\mu$-symmetric Hunt process associated with  $(\EE, \FF)$ and \eqref{e:intcon} holds.
Recall that $\Phi$ is the function defined in \eqref{e:dPhi}. 

We first establish the zero-one law for tail events.
\begin{lemma}\label{l:tailevent} 
	Let $U$ be a tail event {(with respect to the natural filtration of $X$).} Then, either $\P^x(U)=0$ for all $x \in M$ or else $\P^x(U)=1$ for all $x \in M$.
\end{lemma}
\pf By \cite[Lemma 2.7]{HKE}, we have the constant $c>0$ such that for any $x \in M$ and $r,t>0$,
$$ \P^x(\tau_{B(x,r)} \le t) \le \frac{ct}{\Phi(r)}. $$
Let us fix $t_0,\eps>0$ and $x_0 \in M$. Then, by the above inequality and {$L(\alpha_1, c_L,\Phi)$ from Lemma \ref{l:Phiwsc},} there is $c_1>0$ such that
\begin{equation}\label{e:tailevent}
\P^{x_0} \big( \sup_{s \le t_0} d(x_0,X_s) > c_1 \Phi^{-1}(t_0) \big) \le \P^{x_0} \big( \tau_{B(x_0,c_1\Phi^{-1}(t_0))} \le t_0) \le \eps.
\end{equation}
Using \eqref{e:tailevent} and Theorem \ref{t:PHR}, the remainder part of the  proof is the same as those of \cite[Theorem 2.10]{KKW} and \cite[Theorem 5.1]{BKKL}.  Thus, we skip it.\qed

\noindent From \eqref{e:dPhi} and \eqref{e:intcon} with $\VD(d_2)$, we easily see that the following three conditions are equivalent:
\begin{equation}\label{e:comp1}
\sup_{x \in M} \left( {\rm or} \, \inf_{x \in M} \right) \int_M  F(d(x,y)) J(x,dy) < \infty;
\end{equation}
\begin{equation}\label{e:comp2}
\exists\, c>0 \quad \mbox{such that} \quad c^{-1}F(r) \le \Phi(r) \le cF(r), \quad \mbox{for all} \quad r>1;
\end{equation}
\begin{equation}\label{e:comp3}
\int_1^\infty \frac{dF(s)}{\psi(s)} < \infty.
\end{equation}
Since $X$ satisfies $\GHK(\Phi,\psi)$ by Theorem \ref{t:main32}, the above conditions \eqref{e:comp1}-\eqref{e:comp3} are also equivalent to the following moment condition.
\begin{lemma}\label{l:moment}
The following is also equivalent to \eqref{e:comp1}-\eqref{e:comp3}:
\begin{equation*}
 \sup_{x \in M} \left( {\rm or} \,\, \inf_{x \in M} \right) \bE^x[ F(d(x,X_t)) ] < \infty, \quad \forall ( {\rm or} \, \exists) \,\, t>0.
 \end{equation*}
\end{lemma}
\pf (i) Fix $t>0$ and assume \eqref{e:comp3}. 
  Using $\GHK(\Phi,\psi)$, we have for all $x \in M$,
	\begin{align*}
	&\bE^x[F(d(x,X_t))] = \int_M p(t,x,y) F(d(x,y)) \mu(dy) \\
	&= \int_{d(x,y) \le F^{-1}(t)} p(t,x,y) F(d(x,y))\mu(dy) + \int_{d(x,y) > F^{-1}(t)} p(t,x,y) F(d(x,y)) \mu(dy) \\
	&\le \int_{d(x,y) \le F^{-1}(t)} \frac{c_2F(d(x,y))}{V(x,\Phi^{-1}(t))} \mu(dy) + \int_{d(x,y) > F^{-1}(t)} \frac{c_2 tF(d(x,y)) }{V(x,d(x,y))\psi(d(x,y))} \mu(dy) \\
	&\quad + \int_{d(x,y) > F^{-1}(t)} \frac{c_2}{V(x,\Phi^{-1}(t))} e^{-a_0 F_1(d(x,y),F(\Phi^{-1}(t)) )} \mu(dy) := c_2(I + II + III).
	\end{align*}
	Using $\VD(d_2)$ we have
	\begin{align*}
	I &= \int_{d(x,y) \le F^{-1}(t)} \frac{1}{V(x,\Phi^{-1}(t))} F(d(x,y))\mu(dy) \\ & \le \frac{V(x,F^{-1}(t))}{V(x,\Phi^{-1}(t))} F(F^{-1}(t)) \le C_\mu t \Big( \frac{F^{-1}(t)}{\Phi^{-1}(t)} \lor 1\Big)^{d_2} = c_3(t) < \infty.
	\end{align*}
For $II$, we first observe that by $L(\gamma_1,c_1,F)$, there exists $c_4>1$ such that 
$ F(c_4 r) \ge 2F(r)$ for any $r>0$. Using this, $\VD(d_2)$ and \eqref{e:comp3} we obtain
\begin{align*}
II &= \int_{d(x,y) > F^{-1}(t)} \frac{t}{V(x,d(x,y))\psi(d(x,y))} F(d(x,y)) \mu(dy) \\ &= t\sum_{i=0}^\infty \int_{c_4^i F^{-1}(t) < d(x,y) \le c_4^{i+1}F^{-1}(t)} \frac{1}{V(x,d(x,y))\psi(d(x,y))} F(d(x,y)) \mu(dy) \\ &\le
t\sum_{i=0}^\infty \frac{V(x,c_4^{i+1}F^{-1}(t))}{V(x,c_4^{i}F^{-1}(t))}  \frac{F(c_4^{i+1} F^{-1}(t))}{\psi(c_4^{i} F^{-1}(t))}  \le c_5 t \sum_{i=0}^\infty \frac{F(c_4^{i+1} F^{-1}(t))}{\psi(c_4^{i+1} F^{-1}(t))} \\
&\le 2c_5 t \sum_{i=0}^\infty \frac{F(c_4^{i+1} F^{-1}(t)) - F(c_4^i F^{-1}(t))}{\psi(c_4^{i+1} F^{-1}(t))} \le 2c_5 t \sum_{i=0}^\infty \int_{c_4^i F^{-1}(t)}^{c_4^{i+1}F^{-1}(t)} \frac{dF(s)}{\psi(s)} \\
&= 2c_5 t \int_{F^{-1}(t)}^\infty \frac{dF(s)}{\psi(s)} = c_6(t)<\infty.
\end{align*}
For $III$, using \eqref{Phi0lower} for the fourth line, and $\VD(d_2)$ for the fifth line we have
\begin{align*}
III &= \int_{d(x,y)>F^{-1}(t)} \frac{1}{V(x,\Phi^{-1}(t))} \exp \big(-a_0 F_1(d(x,y),F(\Phi^{-1}(t))\big) \mu(dy)	 \\
& = \sum_{i=0}^\infty \int_{2^i F^{-1}(t) < d(x,y) \le 2^{i+1}F^{-1}(t)} \frac{1}{V(x,\Phi^{-1}(t))} \exp \big(-a_0 F_1(d(x,y),F(\Phi^{-1}(t))\big) \mu(dy)	 \\
& \le \sum_{i=0}^\infty  \frac{V(x,2^{i+1}F^{-1}(t))}{V(x,\Phi^{-1}(t))} \exp \big(-a_0 F_1(2^i F^{-1}(t),F(\Phi^{-1}(t))\big) \\
&\le  \sum_{i=0}^\infty  \frac{V(x,2^{i+1}F^{-1}(t))}{V(x,\Phi^{-1}(t))} \exp {\Big(-a_0 \Big(\frac{F(2^i F^{-1}(t))}{F(\Phi^{-1}(t))}\Big)^{\frac{1}{\gamma_2-1}} +a_0 \Big)} \\
&\le (c_7(t) )^{d_2}   \sum_{i=0}^\infty  ( 2^{i+1})^{d_2} \exp\left(- (a(t) 2^i )^{\frac{\gamma_1}{\gamma_2-1}}\right)= c_8(t) < \infty. 
\end{align*}
Combining all the estimates above, we conclude that for any $t>0$, 
	$$ \sup_{x \in M} \bE^x[F(d(x,X_t))] \le c_2(I + II + III) \le c_9(t) < \infty. $$
	
\noindent (ii) Assume that there exist $x \in M$ and $t>0$ such that
$ \bE^x[F(d(x,X_t))] < \infty.$ 
Note that by $\RVD(d_1)$ and $L(\gamma_1,c_1,F)$, we have constants $\theta, c>1$ such that
\begin{equation}\label{e:RVD}
V(x,\theta r) \ge cV(x,r) \mbox{ and } F(\theta r) \ge cF(r), \quad x \in M, r>0.
\end{equation}
Then using $\GHK(\Phi,\psi)$ for the first inequality and \eqref{e:RVD} for the third one, we have
\begin{align*}
&\bE^x[F(d(x,X_t))] = \int_M p(t,x,y) F(d(x,y)) \mu(dy) \\ &\ge c_{10}t \int_{d(x,y)>\eta \Phi^{-1}(t)} \frac{F(d(x,y))}{V(x,d(x,y))\psi(d(x,y))}\mu(dy) 
\\ &= c_{10}t \sum_{i=0}^\infty \int_{\eta \theta^i \Phi^{-1}(t) < d(x,y) \le \eta \theta^{i+1}\Phi^{-1}(t)} \frac{F(d(x,y))}{V(x,d(x,y))\psi(d(x,y))}\mu(dy) \\ &
\ge c_{10}t \sum_{i=0}^\infty \big( V(x, \eta \theta^{i+1}\Phi^{-1}(t)) - V(x, \eta \theta^i \Phi^{-1}(t)) \big) \frac{F(\eta \theta^{i}\Phi^{-1}(t))}{V(x,\eta \theta^{i+1}\Phi^{-1}(t)) \psi(\eta \theta^{i+1}\Phi^{-1}(t))} \\
&\ge  c_{11} t \sum_{i=0}^\infty \frac{F(\eta \theta^{i+1}\Phi^{-1}(t)) - F(\eta \theta^i \Phi^{-1}(t))}{\psi(\eta \theta^i \Phi^{-1}(t))} \ge c_{11} t \sum_{i=0}^\infty \int_{\eta \theta^i \Phi^{-1}(t)}^{\eta \theta^{i+1}\Phi^{-1}(t)} \frac{dF(s)}{\psi(s)} = c_{11}t \int_{\eta \Phi^{-1}(t)}^\infty \frac{dF(s)}{\psi(s)}.
\end{align*}
Here we have used that $\psi,F$ are non-decreasing for the last two lines. This implies \eqref{e:comp3}. 

Combining the above results in (i) and (ii), we conclude the lemma. \qed

Let us define an increasing function $h(t)$ on $[16,\infty)$ by
$$h(t) := (\log \log t) F^{-1} \Big( \frac{t}{\log \log t} \Big).$$ 

\begin{lemma}\label{l:h}
	For any $c_1>0$, $c_2 \in (0,1]$ and $t \in [16,\infty)$,
	\begin{equation}\label{e:h1}
	F_1((c_1+1)h(t),t) \ge c_1 \log \log t
	\end{equation}
	and
{	\begin{equation}\label{e:h2}
	F_1(c_2 h(t),t) \le c_F^{1/(\gamma_1 - 1)} c_2\log \log t.
	\end{equation}}
\end{lemma}
\pf  By the definition of $F_1$, letting $s = h(t)(\log \log t)^{-1}$ we have that for $t \ge 16$,
\begin{align*}
F_1 \big( (c_1+1)h(t), t \big)&= \sup_{s>0} \left( \frac{ (c_1 +1)h(t)   }{s} - \frac{t}{F(s)} \right) \ge \frac{(c_1 +1)h(t)   }{h(t)(\log \log t)^{-1}} - \frac{t}{F(h(t)(\log \log t)^{-1})} \\ &= c_1 \log \log t.
\end{align*}

For \eqref{e:h2}, we fix $t>0$  and let $s_0 :=   {c_F}^{-1/(\gamma_1-1)} h(t) (\log \log t)^{-1} \le h(t)(\log \log t)^{-1}$. If $s \le s_0$,  using $L(\gamma_1,c_F^{-1},F)$ we have
$$ \frac{s}{F(s)} \ge  c_F^{-1} \left(\frac{h(t)(\log \log t)^{-1}}{s}\right)^{\gamma_1 -1} \frac{h(t)(\log \log t)^{-1}}{F(h(t)(\log \log t)^{-1})} \ge c_F^{-1} ( {c_F}) \frac{h(t)}{t} = \frac{h(t)}{t} \ge \frac{c_2 h(t)}{t}. $$
Thus, we obtain $\frac{c_2h(t)}{s} - \frac{t}{F(s)} \le 0$ for $s \le s_0$. Since $F_1(r,t)>0$ for all $r,t>0$, we have
\begin{align*} F_1(c_2h(t),t) &= \sup_{s>0} \left( \frac{c_2 h(t)}{s} - \frac{t}{F(s)} \right) = \sup_{s \ge s_0} \left( \frac{c_2 h(t)}{s} - \frac{t}{F(s)} \right) \le \frac{c_2 h(t)}{s_0} = c_F^{1/(\gamma_1 - 1)} c_2 \log \log t.\end{align*}
\qed

Note that  if $(M,d,\mu)=(\R^n, |\cdot|, dm)$, we have $F(r)=r^2$ and so $h(t) = (t \log \log t)^{1/2}$. 
Thus, the next theorem is the counterpart of \cite[Theorem 5.2]{BKKL}.
\begin{thm}\label{t:lil}
	(i) Assume that \eqref{e:comp1} holds. Then there exists a constant $c \in (0,\infty)$ such that for all $x \in M$,
	\begin{equation}\label{e:lil}
	\limsup_{t \to \infty} \frac{d(x,X_t) }{h(t)} = c  \quad \P^x\text{-a.e}.
	\end{equation}
	(ii) Suppose that \eqref{e:comp1} does not hold. Then for all $x \in M$, \eqref{e:lil} holds with $c = \infty$.
\end{thm}
\pf Fix $x \in M$. We first observe that by \eqref{wsF}, there exist constants $a>16$ and $c_1(a)>1$ such that for any $t \ge 16$,
\begin{equation}\label{e:ha}
(2c_F)^{1/\gamma_1}h(t) \le h(at) \le c_1 h(t). 
\end{equation}
In particular, combining \eqref{e:ha} and $L(\gamma_1,c_F^{-1},F)$ we have
\begin{equation}\label{e:F1}
2F(h(t)) \le F(h(at)). 
\end{equation}
Also, using $L(\gamma_1,c_F^{-1},F)$, we have for $t \ge 16$,
$ \frac{F(h(t))}{t/\log \log t}= \frac{F(h(t))}{F(h(t)/ \log \log t)} \ge c_F^{-1} (\log \log t)^{\gamma_1} $.
Thus, 
\begin{equation}\label{e:h}
c_F^{-1} t (\log \log t)^{\gamma_1 -1} \le F(h(t)), \quad t \ge 16.
\end{equation}
Using \eqref{e:F1}, \eqref{e:ha} and $U(\beta_2,C_U,\psi)$ we obtain that for $n \ge 1$,
\begin{align*}
\int_{h(a^n)}^{h(a^{n+1})} \frac{dF(s)}{\psi(s)} &\ge \big( F(h(a^{n+1})) - F(h(a^n)) \big) \frac{1}{\psi(h(a^{n+1}))} \ge  c_2 \frac{F(h(a^{n+1}))}{\psi(h(a^n))} \\
&\ge c_3 a^{n+1} \frac{(\log\log a^{n+1})^{\gamma_1-1}}{\psi(h(a^n))} \ge c_3 \int_{a^n}^{a^{n+1}} \frac{(\log \log t)^{\gamma_1 -1}}{\psi(h(t))} dt.
\end{align*}
In particular, this and $a>16$ imply that
\begin{equation}\label{e:lil1}
\int_{h(a)}^\infty \frac{dF(s)}{\psi(s)}ds  \ge  c_3 \int_{a}^\infty \frac{1}{\psi(h(t))} dt.
\end{equation}
\noindent (i) Let $k_0 \in \N$ be a natural number satisfying $2^{k_0} \ge a$. By \eqref{e:comp3} and \eqref{e:lil1},
\begin{equation}\label{e:lil3}
\sum_{k=k_0}^\infty \frac{2^k}{\psi(h(2^k))} \le c_4 \sum_{k=k_0}^\infty \int_{2^k}^{2^{k+1}} \frac{dt}{\psi(h(t))} \le c_4 \int_{a}^\infty \frac{dt}{\psi(h(t))} < \infty.
\end{equation}
By \eqref{e:comp2}, we have {$c_8^{-1} t \ge F(\Phi^{-1}(u))$} for any $u \ge 16$ and $t \le u \le 4t$. Thus, we have
\begin{equation}\label{e:lil3.2}
F_1(d(x,y),F(\Phi^{-1}(u))) \ge  F_1 (d(x,y),c_8^{-1} t) = c_8^{-1} F_1(c_8 d(x,y),t). 
\end{equation}
Using $\GHK(\Phi,\psi)$, $\VD(d_2)$ and \eqref{e:lil3.2}  we have 
\begin{align}
&\P^x(d(x,X_u)>Ch(t)) = \int_{\{y: d(x,y)>Ch(t)\}} p(u,x,y)\mu(dy) \nn \\
& \le c_5  t \int_{\{d(x,y)>Ch(t)\}} \frac{\mu(dy)}{V(x,d(x,y))\psi(d(x,y))} + \frac{c_5}{V(x,F^{-1}(t))} \int_{\{d(x,y)>Ch(t)\}} e^{-c_7 F_1(c_8 d(x,y), t)}  \mu(dy) \nn\\
& := c_5(I + II). \label{e:lil3.1}
\end{align}

Let us choose $C={c_{8}^{-1}} (1+ {4}{c_7^{-1}})$ for \eqref{e:lil3.1}. By \cite[Lemma 2.1]{HKE}, we have $I \le c_{11} \frac{t}{\psi(h(t))}$. For $II$, using $\VD(d_2)$ and \eqref{Phi0lower} we have
\begin{align*}
II &= \frac{1}{V(x,F^{-1}(t))} \int_{\{d(x,y)>Ch(t)\}}  \exp(-c_7 F_1(c_{8}d(x,y), t))\mu(dy) \\
&= \frac{1}{V(x,F^{-1}(t))} \sum_{i=0}^\infty \int_{\{C 2^{i}h(t) < d(x,y) \le C2^{i+1} h(t)\}}\exp(-c_7 F_1(c_8 d(x,y), t))\mu(dy) \\
&\le \sum_{i=0}^\infty \frac{V(x,C2^{i+1} h(t))}{V(x,F^{-1}(t))} \exp \Big( - c_7   F_1\big( (1+{4}{c_7^{-1}}) 2^{i} h(t),t \big) \Big) \\
&\le e^{c_7/2}c_{9}\exp \Big( - \frac{c_7}{2}   F_1\big( (1+{4}{c_7^{-1}})  h(t),t \big) \Big) \sum_{i=0}^\infty \Big( \frac{2^i h(t)}{F^{-1}(t)} \Big)^{d_2} \exp \Big( - c_{10}   \Big(\frac{F(2^i h(t))}{t}\Big)^{\frac{1}{\gamma_2-1}} \Big) \\
&\le e^{c_7/2}c_{9}\exp \Big( - \frac{c_7}{2}   F_1\big( (1+{4}{c_7^{-1}})  h(t),t \big) \Big) \sum_{i=0}^\infty \Big( \frac{2^i h(t)}{F^{-1}(t)} \Big)^{d_2} \exp \Big( - c_{11}   \Big(\frac{2^i h(t)}{F^{-1}(t)}\Big)^{\frac{\gamma_1}{\gamma_2-1}} \Big) \\
&\le e^{c_7/2}c_{9}\exp \Big( - \frac{c_7}{2}   F_1\big( (1+{4}{c_7^{-1}})  h(t),t \big) \Big) \sup_{s \ge 1} \sum_{i=0}^\infty (2^i s)^{d_2} \exp\big(-c_{11} (2^i s)^{\frac{\gamma_1}{\gamma_2 -1}}\big) \\
&\le c_{12}\exp \Big( - \frac{c_7}{2}   F_1\big( (1+{4}{c_7^{-1}})  h(t),t \big) \Big).
\end{align*}
Note that by \eqref{e:h}, we have $h(t) \ge cF^{-1}(t)$. Using \eqref{e:h1}, we obtain
$$ II \le  c_{13}\exp \Big( - \frac{c_7}{2}   F_1\big( (1+{4}{c_7^{-1}}) h(t),t \big) \Big)  \le c_{13} \exp \big( -2 \log \log t \big) \le c_{13} (\log t)^{-2}. $$ 
Thus, for $C= {c_{8}^{-1}} (1+ {4}{c_7^{-1}})$ and any $t \ge 16$ and $t \le u \le 4t$, we have
$$ \P^x(d(x,X_u) > Ch(t) ) \le c_{14} \big( \frac{t}{\psi(h(t))} + (\log t)^{-2} \big). $$
Using this and the strong Markov property, for $t_k =2^k$ with $k \ge k_0+1$ we get
\begin{align*}& \P^x (  d(x,X_s) > 2Ch(s) \mbox{ for some } s \in [t_{k-1},t_k]  ) \le \P^x (\tau_{B(x,Ch(t_{k-1}))} \le t_k ) \\
&\le 2 \sup_{s \le t_k, z \in M} \P^z \big( d(z, X_{t_{k+1} -s}) > Ch(t_{k-1}) \big) \le c_{15} \Big(\frac{1}{k^2} + \frac{2^k}{\psi(h(2^k))} \Big).
\end{align*}
Thus, by \eqref{e:lil3} and the Borel-Cantelli lemma, the above inequality implies that
$$\P^x (  d(x,X_t) \le 2Ch(t) \mbox{ for all sufficiently large } t) =1.$$
Thus, $\limsup_{t \to \infty} \frac{d(x,X_t)}{h(t)} \le 2C$.

On the other hand, by \eqref{e:comp2} and $L(\gamma_1,c_F^{-1},F)$, we have $L^1(\gamma_1,c_L,\Phi)$ with some $c_L>0$. Also, by \eqref{e:Phi2} we have $U(\gamma_2,c_F,\Phi)$. Since $\gamma_1>1$, using \eqref{proplowinf} we have for any $x,y \in M$ and $t \ge T$,
\begin{align}\label{e:1}
 p(t,x,y) &\ge  c_{16}V(x,\Phi^{-1}(t))^{-1} \exp\big( -a_L \wt \Phi_1(d(x,y),t)),
 \end{align}
 where $\wt \Phi(r) = r^{\gamma_2}\Phi(1) \1_{\{r < 1\}} + \Phi(r) \1_{\{r \ge 1\}}$ and $\wt \Phi_1(r,t) = \sT(\wt \Phi)(r,t)$ are the functions defined in \eqref{d:wePhi} and \eqref{e:wtPhi1}. Note that by $U(\gamma_2,c_F,F)$ we have $\wt \Phi(r) = r^{\gamma_2} \Phi(1) \le c_F \frac{\Phi(1) F(r)}{F(1)}$ for $r <1$. Using this and \eqref{e:comp2}, we obtain that $\wt \Phi(r) \le cF(r)$ for all $r>0$. Thus, by the definitions of $\wt \Phi_1$ and $F_1$ we obtain
 \begin{equation}\label{e:phif}
 \wt \Phi_1(r,t) \le F_1\big(r,\frac{t}{c}\big), \quad r,t>0.
 \end{equation}
 Combining \eqref{e:1} and \eqref{e:phif}, we have that for all $c_0 \in (0,1), t \ge 16$ and $t \le u \le 4t$,
\begin{align}
\P^x(d(x,X_u)>c_0 h(t)) &= \int_{\{ d(x,y)>c_0 h(t)\}} p(u,x,y)\mu(dy) \nn \\
&\ge  \frac{c_{16}}{V(x,\wt\Phi^{-1}(u))} \int_{\{ d(x,y)>c_0 h(t)\}} e^{-a_L \wt\Phi_1( d(x,y), u ) }   \mu(dy) \nn \\
& \ge \frac{c_{16}}{V(x,F^{-1}(t))} \int_{\{ d(x,y)>c_0 h(t)\}} e^{-a_L F_1(d(x,y), \frac{u}{c})}  \mu(dy).\nn 
\end{align}
Note that by $\RVD(d_1)$, we have a constant $c_{17}>0$ such that
$$ V(x,c_{17}r) \ge 2V(x,r), \quad \mbox{for all} \quad x \in M,\,\, r>0.$$
Thus, using this and \eqref{e:h} we have that for $u\ge t$,
\begin{align*}
& \frac{1}{V(x,F^{-1}(t))} \int_{\{d(x,y)>c_0 h(t)\}} e^{-a_L F_1(d(x,y), \frac{u}{c})}  \mu(dy)\\
& \ge \frac{1}{V(x,F^{-1}(t))} \int_{\{c_0 h(t) < d(x,y) \le c_0c_{17} h(t)\}} e^{-a_L F_1(d(x,y), \frac{u}{c})} \mu(dy) \\
& \ge \frac{V(x,c_0 h(t))}{V(x,F^{-1}(t))} \exp \big( -a_L F_1(c_0 c_{17} h(t), {t}{c^{-1}}) \big) \\
& \ge c_0^{d_2} C_\mu^{-1} \exp \big( -c_{18} F_1(c_0 c_{19} h(t), t) \big).
\end{align*}
Since the constants $c_{16},c_{18}, c_{19}$ are independent of $c_0$, provided $c_0>0$ small and using \eqref{e:h2}, we have
\begin{align}\begin{split}\label{e:lil41}
 &\P^x(d(x,X_u)>c_0h(t)) \ge c_{16}\exp \big( -c_{18} F_1(c_0 c_{19} h(t), t) \big) \\ &\ge c_0^{d_2}c_{16}C_\mu^{-1}\exp\big( - c_{20} c_0 \log \log t\big) \ge c_0^{d_2}c_{16}C_\mu^{-1}(\log t)^{-1/2}. \end{split}\end{align}
 Thus, 
  by the strong Markov property and \eqref{e:lil41}, we have
 $$ \sum_{k=1}^\infty\P^x (d(X_{t_k}, X_{t_{k+1}})  \ge c_0 h(t_{k}) | \FF_{t_k}  )  \ge \sum_{k=4}^\infty c_0^{d_2} c_{16} C_{\mu}^{-1} (\log t_k)^{-1/2} = \infty. $$ 
 Thus, by the second Borel-Cantelli lemma, 
 $\P^x (\limsup \{ d(X_{t_k},X_{t_{k+1}})   \ge c_0h(t_{k})  \}) =1. $ Whence, for infinitely many $k \ge 1$, $d(x,X_{t_{k+1}}) \ge \frac{c_0h(t_{k})}{2}$ or $d(x,X_{t_k}) \ge \frac{c_0h(t_{k})}{2}.$ Therefore, for all $x \in M$,
 $$ \limsup_{ t \to \infty} \frac{d(x,X_t)}{h(t)} \ge  \limsup_{k \to \infty } \frac{d(x,X_{t_k})}{h(t_k)} \ge c_{21}, \quad \P^x\text{-a.e.}, $$
 where $c_{21}>0$ is the constant satisfying $c_{21}h(2t) \le \frac{c_0}{2} h(t)$ for any $t \ge 16$. Since
 $$ \P^x\Big( c_{21} \le \limsup_{t \to \infty} \frac{d(x,X_t)}{h(t)}  \le 2C\Big) = 1, $$
 by Lemma \ref{l:tailevent} there exists a constant $c>0$ satisfying \eqref{e:lil}.

\noindent (ii) Recall that by Lemma \ref{l:Phiwsc}, $L(\alpha_1, c_L, \Phi)$ holds for some $\alpha_1, c_L>0$. Let $t_k =2^k$.  Note that using the assumption $\int_0^\infty \frac{dF(s)}{\psi(s)} = \infty$ to \eqref{e:dPhi} we obtain {$\lim_{t \to \infty} \frac{F(t)}{\Phi(t)} = \infty$}, which implies
\begin{equation}\label{e:lil4.1}
\lim_{s \to \infty} \frac{\Phi^{-1}(s)}{F^{-1}(s)} = \infty.
\end{equation}
Let $\eta>0$ be the constants in \eqref{e:GLHK1}. Let $C_0 \in (0,1)$ be a constant which will be determined later. Define $N :=N(k):= \lceil C_0  \log \log t_k \rceil + 1$. Then, by \eqref{e:lil4.1} we have $\displaystyle  \lim_{k \to \infty} \frac{\Phi^{-1}(t_k/N)}{F^{-1}(t_k/N)} = \infty.$
Thus, there exists $k_1 \in \N$ such that for any $k \ge k_1$, we have $N(k) \ge 3$ and $ \frac{\Phi^{-1}(t_k/N)}{F^{-1}(t_k/N)} \ge \frac{6\theta Ac_F}{\eta C_0} (2C_0c_F)^{1/\gamma_2}$, where $\theta$ is the constant in \eqref{e:RVD}.

Using this and $L(1/\gamma_2,c_F^{-1/\gamma_2},F^{-1})$, we have for $k \ge k_1$,
\begin{align*}
\frac{\eta}{3} \Phi^{-1}(t_k/N) &= \frac{\eta}{3} \frac{\Phi^{-1} \big( \frac{t_k}{N} \big)}{F^{-1}\big( \frac{t_k}{N} \big)} F^{-1}\Big( \frac{t_k}{N} \Big) \ge \frac{2\theta Ac_F}{C_0} (2C_0c_F)^{1/\gamma_2} F^{-1}\Big(\frac{t_k}{N}\Big) \\ 
& \ge \frac{2\theta Ac_F}{C_0} (2C_0c_F)^{1/\gamma_2} F^{-1}\Big(\frac{t_k}{2C_0\log \log t_k}\Big) \\
&\ge \frac{2\theta Ac_F}{C_0} (2C_0c_F)^{1/\gamma_2} (2C_0c_F)^{-1/\gamma_2} F^{-1}\Big(\frac{t_k}{\log \log t_k}\Big) = \frac{2\theta Ac_F h(t_k)}{C_0\log \log t_k} \ge \frac{2\theta Ac_F h(t_k)}{N}.
\end{align*}
Note that by $\Ch(A)$, we have a sequence $\{z_l\}_{l=0}^N$ of points in $M$ such that $z_0=x$, $z_N=y$ and $d(z_{l-1},z_l) \le A \frac{d(x,y)}{N}$ for any $l\in\{1,\dots, N\}$.
Thus, following the chain argument in \eqref{3} equipped with \eqref{e:GLHK1} and $\RVD(d_1)$, we have for $k \ge k_1$ and $2c_F h(t_{k}) \le d(x,y) \le 2\theta c_F h(t_{k})$,
\begin{align*}
&p(t_k,x,y) \geq \int_{B(z_{N-1},\frac{\eta}{3}\Phi^{-1}(t_k/N))} \cdots\int_{B(z_{1},\frac{\eta}{3}\Phi^{-1}(t_k/N))}p(\tfrac{t_k}{N},x,\xi_1)\cdots p(\tfrac{t_k}{N},\xi_{N-1},y)\mu(d\xi_1)\cdots \mu(d\xi_{N-1})\nn \\
&\geq  \prod_{l=0}^{N-1} C_3 V(z_l,\Phi^{-1}(t_k/N))^{-1} \prod_{l=1}^{N-1} c_\mu (\frac{\eta}{3})^{d_1} V(z_l,\Phi^{-1}(t_k/N))  \\
&\geq c_1\left(\frac{C_3 \eta^{d_1}}{3^{d_1}}\right)^N V(x,\Phi^{-1}(t_k))^{-1}
{\geq c_1V(x,\Phi^{-1}(t_k))^{-1}\exp\left({-N\log \Big(\frac{3^{d_1}}{C_3 \eta^{d_1}}\Big) }\right).}
\end{align*}
Since $C_3$ and $\eta$ are the constants in \eqref{e:GLHK1} which are independent of $N$, letting {$C_0 = \frac{1}{4} (\log (\frac{3^{d_1}}{C_3 \eta^{d_1}}))^{-1}$} we obtain that for any $k \ge k_1$ and $x,y \in M$ with $2c_F h(t_k) \le d(x,y) \le 2\theta c_F h(t_k)$,
\begin{align}\begin{split}\label{e:lil5}
 p(t_k,x,y) &\ge \frac{c_1}{V(x,\Phi^{-1}(t_k))}\exp\left({-N\log \Big(\frac{3^{d_1}}{C_3 \eta^{d_1}}\Big) }\right) \\ &\ge  \frac{c_1}{V(x,\Phi^{-1}(t_k))} \exp(-\frac{1}{2} \log \log t_k) = \frac{c_1}{V(x,\Phi^{-1}(t_k))} k^{-1/2}. 
\end{split}\end{align}
Now we claim that for every $C>1$,
\begin{equation}\label{e:lil4}
\sum_{k=1}^\infty \inf_{x \in M} \int_{ d(x,y) \ge Ch(t_{k+1}) } p(t_k,x,y) \mu(dy) = \infty,
\end{equation}
which implies the theorem. Indeed, the strong Markov property and \eqref{e:lil4} imply that for all $C>0$,
$$ \sum_{k=1}^\infty\P^x (d(X_{t_k}, X_{t_{k+1}})  \ge Ch(t_{k+1}) | \FF_{t_k}  ) \ge \sum_{k=1}^\infty\inf_{ x \in M} \int_{ d(x,y) \ge Ch(t_{k+1})} p(t_k,x,y)\mu(dy)= \infty. $$ 
Thus, by the second Borel-Cantelli lemma, 
$\P^x (\limsup \{ d(X_{t_k},X_{t_{k+1}})   \ge Ch(t_{k+1})  \}) =1. $ Whence, for infinitely many $k \ge 1$, $d(x,X_{t_{k+1}}) \ge \frac{Ch(t_{k+1})}{2}$ or $d(x,X_{t_k}) \ge \frac{Ch(t_{k+1})}{2} \ge \frac{Ch(t_k)}{2}.$ Therefore, for all $x \in M$,
$$ \limsup_{ t \to \infty} \frac{d(x,X_t)}{h(t)} \ge  \limsup_{k \to \infty } \frac{d(x,X_{t_k})}{h(t_k)} \ge \frac{C}{2}, \quad \P^x\text{-a.e.} $$
Since the above holds for every $C >1$, the theorem follows.

We now prove the claim \eqref{e:lil4} by considering two cases separately. Let $\eta>0$ be the constant in Proposition \ref{l:LHK-psi}.
Using $\GHK(\Phi,\psi)$ and $\RVD(d_1)$, there is $\lambda \in (0,1)$ such that
\begin{equation}\label{e:lambda}
 \sup_{ t \ge 1} \sup_{x \in M} \int_{d(x,y) < \lambda \Phi^{-1}(t)} p(t,x,y) \mu(dy) \le c_2 \sup_{ t \ge 1} \sup_{x \in M} \frac{V(x, \lambda\Phi^{-1}(t))}{V(x,\Phi^{-1}(t))} < \frac{1}{2}. \end{equation}

\noindent \textit{Case 1} : If there exist infinitely many $k \ge 1$ such that $Ch(t_{k+1}) \le \lambda \Phi^{-1}(t_k)$, then, by \eqref{e:lambda}, for infinitely many $k \ge 1$,
$$ \inf_{x \in M} \int_{d(x,y) \ge \lambda \Phi^{-1}(t_k)} p(t_k,x,y)\mu(dy) = 1 - \sup_{x \in M} \int_{d(x,y) < \lambda \Phi^{-1}(t_k)} p(t_k,x,y) \mu(dy)> 1/2.$$
Thus we get \eqref{e:lil4}.

\noindent \textit{Case 2} : Assume that there is $k_2 \ge 3$ satisfying that $C h(t_{k+1}) \ge \lambda \Phi^{-1}(t_k)$ for all $k \ge k_2$. 
Then, using \eqref{e:lil5}, \eqref{e:RVD} and  $\VD(d_2)$ we have that for every $k \ge k_1 \lor k_2$, 
\begin{align*}
& \inf_{y \in M} \int_{\{z : Ch(t_{k+1}) \le d(y,z) < \theta Ch(t_{k+1})\}} p(t_k,y,z)\mu(dz) \\ &\ge c_3 \frac{(c-1)V(x,Ch(t_{k+1}))}{V(x,\Phi^{-1}(t_k)) k^{1/2}}\ge c_3 \frac{V(x,\eta \Phi^{-1}(t_k))}{V(x,\Phi^{-1}(t_k))) k^{1/2}} \ge  \frac{c_4}{k^{1/2}}. 
\end{align*}
This proves \eqref{e:lil4}. \qed

\section{Examples}\label{s:examples}

In this section, we give some examples which are covered by our results. In this section, we will consider  a metric measure space   $(M, d, \mu)$ satisfying $\Ch(A)$, $\RVD(d_1)$, $\VD(d_2)$ and $\Diff(F)$.

Typical examples of metric measure spaces satisfying the above conditions are unbounded Sierpinski gasket and unbounded generalized Sierpinski carpet in $\R^d$ with $d\ge2$.  In the following two examples, we observe that unbounded Sierpinski gasket and unbounded generalized Sierpinski carpet satisfy $\Ch(A)$, $\RVD(d_1)$, $\VD(d_2)$ and $\Diff(F)$. 

For any $x\in \R^d$, $\lambda>0$ and $U\subset \R^d$, we use notation $x+U:=\{x+y:y\in U\}$ and $\lambda U:=\{\lambda y:y\in U\}$.
\begin{example}\label{ex:Sg} 
   {\rm
{ First, let us define unbounded Sierpinski gasket in $\R^2$. Let $b_0=(0,0)$, $b_1=(1,0)$ and $b_2=(\frac12, \frac{\sqrt{3}}{2})$ be the points in $\R^2$.   Define inductively, $\frF_0:=\{b_0, b_1, b_2\}$  and 
$$\frF_{n+1}:=\frF_n\cup (2^nb_1+\frF_n)\cup (2^nb_2+\frF_n),\quad\text{for}\;\; n\in\bN\cup\{0\}.$$
Let
$\scG'_0=\cup^{\infty}_{n=0}\frF_n$ and  $\scG_0=\scG'_0\cup \scG''_0,$
where $\scG''_0$ is the reflection of $\scG'_0$ in the $y$-axis. Let $\scG_n=2^{-n}\scG_0$ for  $n\in\bZ$, and define
$\scG_{\infty}=\bigcup^{\infty}_{n=0}\scG_n$. 
$\scG:=\overline{\scG_{\infty}}$ is called the unbounded Sierpinski gasket in $\R^2$ (c.f. \cite{BP88}). It is easy to see that $\scG$ is closed and connected.
}

Let us check that unbounded Sierpinski gasket in $\R^2$ satisfies $\Ch(A)$, $\RVD(d_1)$, $\VD(d_2)$ and $\Diff(F)$ with suitable metric $d$ and measure $\mu$. Let $(M_{SG}, d_{SG}, \mu_{SG})$ be the unbounded Sierpinski gasket in $\R^2$. Here $d_{SG}(x,y)$ denotes the length of the shortest path in $M_{SG}$ from $x$ to $y$, and $\mu_{SG}$ is
a multiple of the $d_f$-dimensional Hausdorff measure on  $M_{SG}$ with $d_f=\log 3/\log2$ (see \cite[Lemma 1.1]{BP88}). By \cite[(1.13)]{BP88}, $d_{SG}(x,y)$ is comparable to $|x-y|$ which is the Euclidean distance, which implies $\Ch(A)$. 
Indeed, for any $x,y\in M_{SG}$, let $\frL=\frL_{x,y}\subset M_{SG}\subset \R^2$ be the shortest path in $M_{SG}$. Then, $|\frL|=d_{SG}(x,y)\le A|x-y|$. For any $N\in\bN$ and $i=0,1,\dots, N$, let $z_i$'s be the points in $\frL$ satisfying that $z_0=x$, $z_N=y$ and $|\frL_{z_{i-1}, z_i}|=|\frL|/N$ for all $i=1, 2, \dots, N$, where $\frL_{z_{i-1}, z_i}$ is the path from $z_{i-1}$ to $z_i$ along $\frL$. Then, $d_{SG}(z_{i-1}, z_i)\le |\frL_{z_{i-1}, z_i}|=|\frL|/N\le \frac{A}{N}|x-y|\le \frac{A}{N}d_{SG}(x,y)$ for all $i=1, 2, \dots, N$. 
Also, by \cite[Theorem 1.5]{BP88}, $\Diff(F)$ holds for $F(r)=r^{d_w}$ with $d_w=\log5/\log2>2$. 
{ Moreover, by \cite[Theorem 1.5]{BP88} and \cite[Theorem 3.2]{GHL03}, we see that $\mu_{SG}(B(x,r))\asymp r^{d_f}$ for any $x\in M_{SG}$ and $r>0$. 
Thus, $(M_{SG}, d_{SG}, \mu_{SG})$ satisfies $\Ch(A)$, $\VD(d_f)$, $\RVD(d_f)$ and $\Diff(F)$.}

 The above result also holds for unbounded Sierpinski gaskets constructed in $\R^n$($n\ge3$) with different $d_f(n)>0$ and $d_w(n)>1$ (see \cite[Section 10]{BP88}). 
}
\end{example}

 Now, we consider unbounded generalized Sierpinski carpet and pre-Sierpinski carpet in $\R^d$($d\ge2$).
\begin{example}\label{ex:gpSg} 
{\rm
 Let $\scF_0=[0,1]^d$ and $l\in\bN$, $l\ge3$ be fixed. For $n\in\bZ$, let $\sS_n$ be the collection of closed cubes of side length $l^{-n}$ with vertices in $l^{-n}\bZ^d$. For $S\in \sS_n$, let $\Psi_S$ be the orientation preserving affine map which maps $\scF_0$ onto $S$.
We define a decreasing sequence $(\scF_n)_{n\ge0}$ of closed subsets of $\scF_0$ as follows. 
For $U\subset \R^d$, define
$\sS_n(U)=\{S:S\subset U, S\in \sS_n\}.$
Let $2\le m<l^d$ be an integer and let $\scF_1$ be the union of $m$ distinct elements of $\sS_1(\scF_0)$. We assume the following conditions on $\scF_1$:
\begin{itemize}
\item[(H1)](Symmetry) $\scF_1$ is preserved by all the isometries of the unit cube $\scF_0$.
\item[(H2)](Connectedness) int$(\scF_1)$ is connected, and contains a path connecting the hyperplanes $\{x_1=0\}$ and $\{x_1=1\}$.
\item[(H3)](Non-diagonality) Let $B$ be a cube in $\scF_0$ which is the union of $2^d$ distinct elements of $\sS_1$.(So $B$ has side length $2l^{-1}$.) If int$(\scF_1\cap B)$ is non-empty, then it is connected.
\item[(H4)](Borders included) $\scF_1$ contains the line segment $\{x:0\le x_1\le1, x_2=\dots=x_d=0\}$.
\end{itemize}
We may think of $\scF_1$ as being derived from $\scF_0$ by removing the interiors of $l^d-m$ squares in $\sS_1(\scF_0)$. Given $\scF_1$ satisfying the above conditions, we define $\scF_{n+1}$ inductively
$$\scF_{n+1}=\bigcup_{S\in \sS_n(\scF_n)}\Psi_S(\scF_1)=\bigcup_{S\in \sS_1(\scF_1)}\Psi_S(\scF_n),\quad n\in\bN.$$
For $n\in\bN\cup\{0\}$, let 
$$\wt\scF_n=\bigcup^{\infty}_{k=0}l^k \scF_{n+k}, \qquad \wt\scF=\bigcap^{\infty}_{n=0}\wt\scF_n.$$
$\wt\scF_0$ and  $\wt\scF$ are generalized pre-Sierpinski carpet and unbounded generalized Sierpinski carpet in $\R^d$, respectively (c.f. \cite{BB}). If $l=3$ and $\scF_1=\sS_1(\scF_0)\setminus D$, where $D=\{x\in [0,1]^d: x_i\in(\frac{1}{3}, \frac{2}{3}), i=1,2,\dots, d\}$ (so, $m=3^d-1$), $\wt\scF_0$ is called pre-Sierpinski carpet (c.f. \cite{AB}).

Let $(M_{SC}, d_{SC}, \mu_{SC})$ be the unbounded generalized Sierpinski carpet  in $\R^n$ ($n\ge2$). Here, $d_{SC}(x,y)$ is the length of the shortest path in $M_{SC}$ from $x$ to $y$ and $\mu_{SC}$ is a multiple of the $d_f$-dimensional Hausdorff measure on $M_{SC}$ with $d_f=\log m/\log l$. Then, by \cite[Remark 2.2]{BB},  $d_{SC}(x,y)\asymp |x-y|$ for all $x,y\in M_{SC}$. Thus, by the same argument as in the Sierpinski gasket case,  $(M_{SC}, d_{SC})$ satisfies $\Ch(A)$. Also, by \cite[Theorem 1.3]{BB},  $\Diff(F)$ holds for $F(r)=r^{d_w}$ with $d_w\ge2$.
{ Moreover, by \cite[Theorem 1.3]{BB} and \cite[Theorem 3.2]{GHL03}, we see that  $\mu_{SC}(B(x,r))\asymp r^{d_f}$ for any $x\in M_{SC}$ and $r>0$. Thus,  $(M_{SC}, d_{SC}, \mu_{SC})$ satisfies $\Ch(A)$, $\VD(d_f)$, $\RVD(d_f)$ and $\Diff(F)$.} 

 Let $(\tilde M, \tilde d, \tilde \mu)$ be a pre-Sierpinski carpet in $\R^n$($n\ge3$). Here, $\tilde d(x,y)$ denotes the shortest path distance in $\tilde M$ and $\tilde \mu$ is Lebesgue measure restricted to $\tilde M$. By \cite[(7.2)]{AB} and the same argument as in the Sierpinski gasket case, we see that $\Ch(A)$ holds. The conditions $\RVD(d_1)$, $\VD(d_2)$ and $\Diff(F)$ follow from \cite[(7.3) and Theorem 7.1(1)]{AB}. In this case, $F(r)=r^2$ if $0<r\le1$ and $F(r)=r^{d_w}$ with $d_w>2$ if $r>1$. }
\end{example}

Throughout the remainder of this section,   we will assume $(M, d, \mu)$ is a metric measure space satisfying $\Ch(A)$, $\RVD(d_1)$, $\VD(d_2)$, and $\Diff(F)$, where the function $F:(0,\infty)\to(0,\infty)$ is strictly increasing function satisfying $L(\gamma_1,c_F^{-1})$ and $U(\gamma_2,c_F)$ with some constants $1<\gamma_1 \le \gamma_2$. 
We also assume that $\psi$ is a non-decreasing function satisfying \eqref{e:intcon},  $L(\beta_1,C_L)$ and $U(\beta_2,C_U)$.

 Let $X$ be  the symmetric pure-jump Hunt process  on $(M,d,\mu)$, which is
 associated with the regular Dirichlet form $(\sE,\sF)$ in \eqref{e:DF} 
 satisfying $\J_{\psi}$ and 
 $p(t,x,y)$ be the transition density of $X$.
In this section, we will use the notation $f(\cdot)\simeq g(\cdot)$ at $\infty$ (resp. $0$) if
 $\frac{f(t)}{g(t)}\to 1$ as $t\to \infty$ (resp. $t\to 0$).

 \begin{definition}\label{D:E1} {\rm \begin{itemize}
\item[(i)] 
We denote $\sR_0^\infty$ (resp. $\sR_0^0$) by the class of slowly varying functions at $\infty$ (resp. $0$). 
\item[(ii)] 
For $\ell\in \sR_0^\infty$, we denote $\Pi_\ell^\infty$ (resp. $\Pi_\ell^0$) by the class of { positive} measurable function $f$  on $[c,\infty)$ (resp. $(0,c)$) 
such that for all $\lambda>0$,
$f(\lambda \cdot)-f(\cdot)
\simeq (\log \lambda)  \ell(\cdot)$ at $\infty$ (resp. $0$). 
  $\Pi_\ell^\infty$ (resp. $\Pi_\ell^0$) is called de Haan class at $\infty$ (resp. $0$) determined by $\ell$. 
  \item[(iii)]
   For $\ell \in \sR_0^\infty$ (resp. $\sR_0^0$), we say $\ell_\#$ is de Bruijn conjugate of $\ell$ if both $\ell(t)\ell_\#(t\ell(t)) \simeq 1$ and $ \ell_\#(t)\ell(t\ell_\#(t))\simeq 1$ at $\infty$ (resp. $0$).
\end{itemize}}
\end{definition}

Note that $f\in \sR_0^\infty$ if  $f\in \Pi^\infty_{\ell}$ (see \cite[Theorem 3.7.4]{BGT}). {  By \cite[Lemma 1.3.2]{BGT}, we see that for $\ell\in \sR_0^\infty$, there exists $C>0$ such that $\ell\in L^1_{loc}[C,\infty)$. 
For $0<c<C$, by letting $\ell(s)=1$ for $0<c\le s\le C$, we may assume that $\ell\in\sR_0^\infty$ belongs to $L^1_{loc}[c,\infty)$. If $\ell\in\sR_0^\infty\cap L^1_{loc}[c,\infty)$ satisfies $\int^\infty_c \ell(s)s^{-1}ds=\infty$, then by \cite[Proposition 1.5.9a and Theorem 3.7.3]{BGT}, we see that for any $f\in \Pi^\infty_{\ell}$,
\begin{align}\label{dhequiv}
f(s)^{-1}\int^s_c \ell(u)u^{-1}du\to1\quad\text{as}\;s\to\infty. 
\end{align}
Note that for $c'\ge0$, $s\mapsto c'+\int^s_c \ell(u)u^{-1}du$ belongs to $\Pi^\infty_{\ell}$ (see \cite[Theorem 3.6.6]{BGT}).
}

In the following corollaries, their proofs and examples, $a_i=a_{i,L}$ or $a_i=a_{i,U} $ depending on whether we consider lower or upper bound.

\begin{corollary}\label{c:example}
Suppose $F$ is differentiable function satisfying  $F(s)\asymp s F'(s)$. Let  $T\in (0,\infty)$. 

\noindent
{ (i) Suppose there exists $\ell\in \sR_0^0$ such that  $\psi(s) \asymp \frac{F(s)}{\ell(s)}$ for $s<T$. Let $f(s):=\int_0^{s\wedge T} {\ell(u)}{u^{-1}}du\in \Pi^0_{\ell}$.  Then, for $t<T$,}
\begin{align*}
&p(t,x,y) \asymp \\
&\;\;\frac{1}{V(x, (F/f)^{-1}(t))}\wedge \left(\frac{ t}{V(x,d(x,y))\psi(d(x,y))} + \frac{1}{V(x, (F/f)^{-1}(t))} \exp{\left(\frac{-a_1\, d(x,y)}{(F'/f)^{-1}(t/d(x,y))}\right)}\hspace{-1mm}\right). 
\end{align*}

\noindent
{ (ii)   Suppose  there exists $\ell\in \sR_0^\infty\cap L^1_{loc}[T,\infty)$ such that $\int_{T}^\infty {\ell(u)}{u^{-1}}du=\infty $ and $\psi(s) \asymp \frac{F(s)}{\ell(s)}$ for $s>T$. Let $f(s):=\int^{T\wedge s}_{0}\frac{F'(u)}{\psi(u)}du+\int^{s}_{T\wedge s}{\ell(u)}{u^{-1}}du \in \Pi^\infty_{\ell}$. Then, for $t>T$,}
\begin{align*}
&p(t,x,y) \asymp \\
&\;\;\frac{1}{V(x, (F/f)^{-1}(t))}\wedge \left(\frac{ t\ell(d(x,y))}{V(x,d(x,y))F(d(x,y))} + \frac{1}{V(x, (F/f)^{-1}(t))} \exp{\left(\frac{-a_2\, d(x,y)}{(F'/f)^{-1}(t/d(x,y))}\right)}\hspace{-1mm}\right).
\end{align*}
\end{corollary}

\pf Let $r=d(x,y)$. Note that  $\Phi(u)= {F(u)}/{\int_0^{u} \frac{F'(s)}{\psi(s)} ds}$. By \eqref{e:Phi2} and Lemma \ref{l:Phiwsc}, we see that $\Phi$ satisfies $U(\gamma_2, c_F)$ and $L(\alpha_1, c_L)$ for some $\alpha_1, c_L>0$. Without loss of generality, we assume that $T=1$.

\noindent(i)  Since $\Phi(s)\asymp F(s)/f(s)$ for $s<1$ and $f\in\sR^{0}_{0}$, we observe that $\Phi$ satisfies $L_1(\delta, \wt C_L)$ for some $\delta>1$. Thus, we can apply  Corollary \ref{c:main32}(i)  and  obtain that for $t<1$,
\begin{align*}
p(t,x,y) \asymp \frac{1}{V(x, (F/f)^{-1}(t))}\wedge \left(\frac{ t}{V(x, r)\psi(r)} + \frac{1}{V(x, (F/f)^{-1}(t))} \exp{\big(-a_1 \Phi_1(r,t)\big)}\hspace{-1mm}\right). 
\end{align*}
By Lemma \ref{lem:inverse} and $L_1(\delta, \wt C_L, \Phi)$, there exists $c_1>0$ such that $\Phi^{-1}(t)>c_1 t^{1/\delta}$ holds for all $t<1$.  Thus, for $r>2c_F^2\Phi^{-1}(t)$ and $t<1$, we have $r>2c_F^2c_1 t^{1/\delta}>2c_F^2c_1t$ which shows $t/r<\frac12 c_1^{-1}c_F^{-2}$. Note that by \eqref{e:scK1}, $\scK(s)\asymp \Phi(s)/s\asymp (F'/f)(s)$ for $s<\frac12 c_1^{-1}c_F^{-2}$.
From this and Lemma \ref{l:scK}, we see that $\Phi_1(r,t)\asymp r/\scK^{-1}(t/r)\asymp r/(F'/f)^{-1}(t/r)$ for  $r>2c_F^2\Phi^{-1}(t)$ and $t<1$. This,  combined with Lemma \ref{l:sT}(ii),
completes the proof of (1).

\vspace{3mm}

\noindent(ii) 
 Note that by the definition of $f$ and the assumption  $\psi(s) \asymp \frac{F(s)}{\ell(s)}$ for $s>1$, we have 
\begin{align}\label{eex1}
f(s)\asymp \int^{s}_{0}\frac{F'(u)}{\psi(u)}du\quad\text{for}\;\;s>0.
\end{align}
Thus, $\Phi(s)\asymp F(s)/f(s)$ for $s>0$. Using this and $f\in\sR^{\infty}_{0}$, we observe that $\Phi$ satisfies $L^1(\delta, \wt C_L)$ for some $\delta>1$. Let $\wt\Phi$ be the function defined in \eqref{d:wePhi}. Then, $U(\gamma_2, c_F, \wt\Phi)$, $L(\delta, \wt C_L, \wt\Phi)$ hold and $\wt\Phi(s)=\Phi(s)$  for $s>1$. By applying Corollary \ref{c:main32}(ii)
, we obtain that for $t\ge1$,
\begin{align*}
p(t,x,y) \asymp \frac{1}{V(x, (F/f)^{-1}(t))}\wedge \left(\frac{ t\ell(r)}{V(x, r)F(r)} + \frac{1}{V(x, (F/f)^{-1}(t))} \exp{\big(-a_2 \wt\Phi_1(r,t)\big)}\hspace{-1mm}\right). 
\end{align*}
Choose small $\theta_0>0$ such that $\frac{1}{\delta}+\theta_0(\frac{1}{\delta}-\frac{1}{\beta_2})=:\eps_0<1$ and let $c_2=c_2(\theta_0)>0$ be a constant satisfying $s^{-d_2-\beta_2-\beta_2/\theta_0}\ge c_2\exp(-a_2 s)$ for all $s>0$. For $r >  2c_F^2\frac{\wt\Phi^{-1}(t)^{1+\theta_0}}{\psi^{-1}(t)^{\theta_0}}$, there exists $\theta>\theta_0$ such that $r=2c_F^2\frac{\wt\Phi^{-1}(t)^{1+\theta}}{\psi^{-1}(t)^{\theta}}$ since $\wt\Phi^{-1}(t)>\psi^{-1}(t)$.  Then, by using $\VD(d_2)$, $U(\beta_2, C_U, \psi)$ and the same argument in the proof of Lemma \ref{l:exp1},
\begin{align*}
\frac{t}{V(x,r)\psi(r)}
&\ge \frac{c_3}{V(x, \wt\Phi^{-1}(t))}\left(\frac{\wt\Phi^{-1}(t)}{\psi^{-1}(t)}\right)^{\theta(-d_2-\beta_2-\beta_2/\theta)}
\ge \frac{c_3}{V(x, \wt\Phi^{-1}(t))}\left(\frac{\wt\Phi^{-1}(t)}{\psi^{-1}(t)}\right)^{\theta(-d_2-\beta_2-\beta_2/\theta_0)}\\
&\ge \frac{c_2c_3}{V(x, \wt\Phi^{-1}(t))} \exp{\left(-a_2 \left(\frac{\wt\Phi^{-1}(t)}{\psi^{-1}(t)}\right)^{\theta} \right)}
= \frac{c_2c_3}{V(x, \wt\Phi^{-1}(t))} \exp{\left(-a_2 \big(r/\wt\Phi^{-1}(t)\big)\right)}\\
&\ge\frac{c_4}{V(x, \wt\Phi^{-1}(t))} \exp{\left(-a_2\wt\Phi_1(r,t)\right)},
\end{align*}
where the last inequality follows from the definition of $\wt\Phi_1(r,t)$. Thus, by this and Lemma \ref{l:sT}(ii), it suffices to estimate $\wt\Phi_1(r,t)$ for $2c_F^2\wt\Phi^{-1}(t) <r \le  2c_F^2\frac{\wt\Phi^{-1}(t)^{1+{\theta_0}}}{\psi^{-1}(t)^{\theta_0}}$ and $t\ge1$. By Lemma \ref{lem:inverse},  $L(\delta, \wt C_L, \wt\Phi)$ and $U(\beta_2, C_U, \psi)$, we see that there exists $c_5>0$ such that $r \le  c_5 t^{\eps_0}$ holds for  $r \le  2c_F^2\frac{\wt\Phi^{-1}(t)^{1+{\theta_0}}}{\psi^{-1}(t)^{\theta_0}}$ and $t\ge1$.
Thus, we have $t/r\ge c_5^{-1}t^{1-\eps_0}\ge c_5^{-1}$.
Define $ \scK_{\infty}(r):= \sup_{0<s \le r} \frac{\wt\Phi(s)}{s}.$
Then, by \eqref{e:scK1} and \eqref{eex1},  $\scK_{\infty}(s)\asymp \wt\Phi(s)/s\asymp (F'/f)(s)$ for $s>c_5^{-1}$. Therefore, by Lemma \ref{l:scK}, we have $\wt\Phi_1(r,t)\asymp r/\scK_{\infty}^{-1}(t/r)\asymp r/(F'/f)^{-1}(t/r)$ for  $2c_F^2\wt\Phi^{-1}(t) <r \le  2c_F^2\frac{\wt\Phi^{-1}(t)^{1+{\theta_0}}}{\psi^{-1}(t)^{\theta_0}}$ and $t\ge1$. This completes the proof.
\qed

\begin{remark}\label{r:ex}
{\rm 

\noindent(1) Let $H$ be a strictly increasing function with weak scaling property. Suppose that $f\in \sR^{0}_{0}$ and there  exists $\wt f\in\sR^{0}_{0}$  such that
$H(r)/f(r)\asymp H(r\wt f(r))=:h(r)$. Note that $h$ also satisfies weak scaling property at $0$. By the weak scaling property of $H$, we see that for $r<1$,
\begin{align*}
h\Big(H^{-1}(r){\wt f}_{\#}(H^{-1}(r))\Big)=H\Big(H^{-1}(r){\wt f}_{\#}(H^{-1}(r))\wt f(H^{-1}(r){\wt f}_{\#}(H^{-1}(r))\Big)\asymp r.
\end{align*}
Thus, by using the weak scaling property of $h$, 
$$(H/f)^{-1}(r)\asymp H^{-1}(r){\wt f}_{\#}(H^{-1}(r))\quad\text{for}\;\;r<1.$$
 
Similarly, suppose that $g\in \sR^{\infty}_{0}$ and there exists $\wt g\in\sR^{\infty}_{0}$ such that $H(r)/g(r)\asymp H(r\wt g(r))$. Then,
$$(H/g)^{-1}(r)\asymp H^{-1}(r){\wt g}_{\#}(H^{-1}(r))\quad\text{for}\;\;r\ge1.$$
(c.f. \cite[Proposition 1.5.15]{BGT}.)

\vspace{3mm}

\noindent(2) Let $\ell\in \sR_0^\infty$(resp. $ \sR_0^0$). Suppose that $\ell$ satisfies 
\begin{align}\label{ellcheck}
\lim_{s\to \infty (\text{resp.}s\to 0)}\frac{\ell(s\ell^\eta(s))}{\ell(s)} = 1  \text{ for some } \eta \in \R\setminus\{0\}. 
\end{align}
Then by \cite[Corollary 2.3.4]{BGT}, $(\ell^\eta)_\# \simeq 1/\ell^\eta$ at $\infty$ (resp. $0$).

}
\end{remark}

Combining Corollary \ref{c:example} and Remark \ref{r:ex}, we have the following:
\begin{corollary}\label{c:example2}
Suppose $F$ is differentiable function satisfying $F(s)\asymp sF'(s)$. Let  $T\in (0,\infty)$ and $\gamma_3, \gamma_4>1$.

\noindent
{ (i) Suppose there exist $\ell_0, \ell_1\in \sR_0^0$ such that $F(s)\asymp s^{\gamma_3}\ell_1(s)$ and $\psi(s)\asymp {F(s)}/{\ell_0(s)}$ for $s<T$. Suppose further that  $f_1(s):=\ell_1(s)^{-1}\int_0^{s\wedge T} {\ell_0(u)}{u^{-1}}du\in \sR^0_{0}$ satisfies \eqref{ellcheck} for $\eta\in\{1/\gamma_3,1/(\gamma_3-1)\}$. Then,  for $t<T$,}
\begin{align*}
&p(t,x,y) \asymp \; \frac{1}{V(x, t^{1/\gamma_3}f_1^{1/\gamma_3}(t^{1/\gamma_3}))}\\
	&\wedge \left(\frac{ t}{V(x,d(x,y))\psi(d(x,y))} + \frac{1}{V(x, t^{1/{\gamma_3}}f_1^{1/\gamma_3}(t^{1/\gamma_3}))} \exp{\left(-a_1\left(\frac{d(x,y)^{\gamma_3}}{tf_1((\frac{t}{d(x,y)})^{\frac{1}{\gamma_3-1}})}\right)^{\frac{1}{\gamma_3-1}}\right)}\right). 
\end{align*}
Furthermore, if  $f_1$ is monotone and satisfies $f_1(s^{\gamma_3})\asymp f_1(s)$ for $s<1$, then for  $t <T$,
\begin{align}\begin{split}\label{e:smallt}
&p(t,x,y)\asymp \;\frac{1}{V(x, t^{1/\gamma_3}f_1^{1/\gamma_3}(t))}\\
	&\wedge \left(\frac{ t}{V(x,d(x,y))\psi(d(x,y))} + \frac{1}{V(x, t^{1/\gamma_3}f_1^{1/\gamma_3}(t))} \exp{\left(-a_2\left(\frac{d(x,y)^{\gamma_3}}{tf_1(t)}\right)^{\frac{1}{\gamma_3-1}}\right)}\right). 
\end{split}\end{align}

\noindent
{ (ii) Suppose there exist $\ell_2\in \sR_0^\infty$ and $\ell_{3}\in \sR_0^\infty\cap L^1_{loc}[T,\infty)$  such that  $\int_{T}^{\infty} {\ell_{3}(u)}{u^{-1}}du=\infty$, $\ell_2(s)\asymp s^{-\gamma_4}F(s) $ and ${\ell_3(s)}\asymp {F(s)}/\psi(s)$ for $s>T$. Suppose further that $f_2(s) :=\ell_2(s)^{-1}\int_0^{T\wedge s}\frac{F'(u)}{\psi(u)}du+\ell_2(s)^{-1}\int_{T\wedge s}^{s}{\ell_3(u)}{u^{-1}}du\in\sR^{\infty}_{0}$  satisfies \eqref{ellcheck} for $\eta\in\{1/\gamma_4,1/(\gamma_4-1)\}$. Then,  for $t>T$,}
\begin{align*}
&p(t,x,y) \asymp \;\frac{1}{V(x, t^{1/\gamma_4}f_2^{1/\gamma_4}(t^{1/\gamma_4}))}\\
&\wedge\left(\frac{ t\ell(d(x,y))}{V(x,d(x,y))d(x,y)^{\gamma_4}} + \frac{1}{V(x, t^{1/\gamma_4}f_2^{1/\gamma_4}(t^{1/\gamma_4}))} \exp{\left(-a_3\left(\frac{d(x,y)^{\gamma_4}}{tf_2((\frac{t}{d(x,y)})^{\frac{1}{\gamma_4-1}})}\right)^{\frac{1}{\gamma_4-1}}\right)}\right). 
\end{align*}
\vspace{-0.2cm}
Furthermore, if  $f_2$ is monotone and satisfies $f_2(s^{\gamma_4})\asymp f_2(s)$ for $s>1$, then for  $t >T$,
\begin{align}\begin{split}\label{e:larget}
&p(t,x,y) \asymp \;\frac{1}{V(x, t^{1/\gamma_4}f_2^{1/\gamma_4}(t))}\\
&\wedge \left(\frac{ t\ell(d(x,y))}{V(x,d(x,y))d(x,y)^{\gamma_4}} + \frac{1}{V(x, t^{1/\gamma_4}f_2^{1/\gamma_4}(t))} \exp{\left(-a_4\left(\frac{d(x,y)^{\gamma_4}}{tf_2(t)}\right)^{\frac{1}{\gamma_4-1}}\right)}\right). 
\end{split}\end{align}
\end{corollary}

\pf Let $r=d(x,y)$, $\Phi(u):= {F(u)}/{\int_0^{u} \frac{F'(s)}{\psi(s)} ds}$. Without loss of generality, we assume that $T=1$.

\noindent(i) We will apply Corollary \ref{c:example}(i) to prove the claim. Since $\Phi(s)\asymp s^{\gamma_3}/f_1(s)$ for $s<1$ and $f_1\in\sR^{0}_{0}$, we observe that $\Phi$ satisfies $L_1(\delta, \wt C_L)$ for some $1<\delta<\gamma_3+1$ and $U_1(\gamma_3+1, c_U)$. By Lemma \ref{lem:inverse} and $L_1(\delta, \wt C_L, \Phi)$, there exists $c_1>0$ such that $\Phi^{-1}(t)>c_1 t^{1/\delta}$ holds for all $t<1$. Thus, for $r>\Phi^{-1}(t)$ and $t<1$, we have $r>c_1 t^{1/\delta}>c_1 t$ which shows $t/r\le c_1^{-1}$.

Define $f(s):=\int_0^{s\wedge1} {\ell_0(u)}{u^{-1}}du$. By Remark \ref{r:ex} and the condition that $f_1$ satisfies \eqref{ellcheck} for $\eta\in\{ 1/\gamma_3,1/(\gamma_3-1)\}$, we see that for $s<1$, 
$$(F/f)^{-1}(s)\asymp(s^{\gamma_3}/f_1)^{-1}(s)\asymp s^{1/\gamma_3}(1/f_1^{1/{\gamma_3}})_{\#}(s^{1/{\gamma_3}})\asymp  s^{1/\gamma_3}f_1^{1/\gamma_3}(s^{1/\gamma_3})$$ 
and 
$$
(F'/f)^{-1}(s)\asymp(s^{\gamma_3-1}/f_1)^{-1}(s)\asymp s^{1/{(\gamma_3-1)}}(1/f_1^{1/{(\gamma_3-1)}})_{\#}(s^{1/{(\gamma_3-1)}})\asymp  s^{1/(\gamma_3-1)} f_1^{1/(\gamma_3-1)}(s^{1/(\gamma_3-1)}).$$
Using this, volume doubling property and the fact that $t/r\le c_1^{-1}$ for $r>\Phi^{-1}(t)$ and $t<1$, we can apply Corollary \ref{c:example}(i) to obtain the first claim. 

Now we prove \eqref{e:smallt}. By using  $f_1(s^{\gamma_3})\asymp f_1(s)$ for $s<1$ and volume doubling property, we have ${V(x, t^{1/\gamma_3}f_1^{1/\gamma_3}(t^{1/\gamma_3}))}\asymp {V(x, t^{1/\gamma_3}f_1^{1/\gamma_3}(t))}$ for $t<1$. Thus, it is enough to show that
$f_1((t/r)^{1/(\gamma_3-1)})$ in the exponential term is comparable to  $f_1(t)$. To show this, we choose  small $\theta_1>0$ such that $\frac{1}{\gamma_3+1}+\theta_1(\frac{1}{\gamma_3+1}-\frac{1}{\beta_1})=:\eps_1\in(0,\delta^{-1})$. 
Let $c_2=c_2(\theta_1)>0$ be a constant satisfying $s^{-d_2-\beta_2-\beta_2/\theta_1}\ge c_2\exp(-a_1 s)$ for all $s>0$. For $r >  \frac{\Phi^{-1}(t)^{1+\theta_1}}{\psi^{-1}(t)^{\theta_1}}$, there exists $\theta>\theta_1$ such that $r=\frac{\Phi^{-1}(t)^{1+\theta}}{\psi^{-1}(t)^{\theta}}$ since $\Phi^{-1}(t)>\psi^{-1}(t)$.  Then, by using $\VD(d_2)$, $U(\beta_2, C_U, \psi)$ and the same argument in the proof of Lemma \ref{l:exp1},
\begin{align*}
\frac{t}{V(x,r)\psi(r)}
&\ge \frac{c_3}{V(x, \Phi^{-1}(t))}\left(\frac{\Phi^{-1}(t)}{\psi^{-1}(t)}\right)^{\theta(-d_2-\beta_2-\beta_2/\theta)}
\ge \frac{c_3}{V(x, \Phi^{-1}(t))}\left(\frac{\Phi^{-1}(t)}{\psi^{-1}(t)}\right)^{\theta(-d_2-\beta_2-\beta_2/\theta_1)}\\
&\ge \frac{c_2c_3}{V(x, \Phi^{-1}(t))} \exp{\left(-a_1 \left(\frac{\Phi^{-1}(t)}{\psi^{-1}(t)}\right)^{\theta} \right)}
= \frac{c_2c_3}{V(x, \Phi^{-1}(t))} \exp{\left(-a_1 \big(r/\Phi^{-1}(t)\big)\right)}\\
&\ge\frac{c_4}{V(x, \Phi^{-1}(t))} \exp{\left(-a_1\Phi_1(r,t)\right)},
\end{align*}
where the last inequality follows from the definition of $\Phi_1(r,t)$. Thus, it is enough to show $f_1((\tfrac{t}{r})^{1/(\gamma_3-1)}) \asymp f_1(t)$ for $\Phi^{-1}(t) < r \le  \frac{\Phi^{-1}(t)^{1+{\theta_1}}}{\psi^{-1}(t)^{\theta_1}}$.  By  Lemma \ref{lem:inverse}, $L_1(\delta, \wt C_L, \Phi)$, $U_1(\gamma_3+1, c_U, \Phi)$ and $L(\beta_1,C_L, \psi)$, we have $\frac{\Phi^{-1}(t)^{1+{\theta_1}}}{\psi^{-1}(t)^{\theta_1}} \le c_5 t^{\frac{1}{\gamma_3+1}+\theta_1(\frac{1}{\gamma_3+1}-\frac{1}{\beta_1})}=c_5 t^{\varepsilon_1}$ and $\Phi^{-1}(t)\ge c_1 t^{1/\delta}$ for $t < 1$. Thus, $c_1t^{1/\delta} < r \le  c_5t^{\eps_1}$ for $t<1$ and $\Phi^{-1}(t) < r \le  \frac{\Phi^{-1}(t)^{1+{\theta_1}}}{\psi^{-1}(t)^{\theta_1}}$. Note that $0<\eps_1<\delta^{-1}<1$. Using that $f_1$ is monotone and $c_1t^{1/\delta} < r \le  c_5t^{\eps_1}$,  we see that 
$$f_1(c_5't^{(1-\eps_1)/(\gamma_3-1)})\le f_1((t/r)^{1/(\gamma_3-1)})\le f_1(c_1't^{(1-1/\delta)/(\gamma_3-1)})$$ 
or 
$$f_1(c_5't^{(1-\eps_1)/(\gamma_3-1)})\ge f_1((t/r)^{1/(\gamma_3-1)})\ge f_1(c_1't^{(1-1/\delta)/(\gamma_3-1)}).$$  
By Lemma \ref{comp_mon} and the fact that $f_1$ is slowly varying at $0$, we have $f_1(c_6 s^{c_7})\asymp f_1(s)$ for $s<1$. Thus, the above inequalities give that $f_1((t/r)^{1/(\gamma_3-1)})\asymp f_1(t)$ for $t<1$.

\vspace{3mm}

\noindent(ii) The proof is similar to that of (i).  We will apply Corollary \ref{c:example}(ii) to prove the claim. 

Define $f(s):=\int_0^{1\wedge s} \frac{F'(u)}{\psi(u)}du+\int_{1\wedge s}^s {\ell_3(u)}{u^{-1}}du\in \Pi^\infty_{\ell_3}$. Then, by \eqref{eex1}, $\Phi(s)\asymp {F(s)}/{f(s)}\asymp s^{\gamma_4}/f_2(s)$ for $s>1$. { Since $f(s)^{-1} \int_{1}^s {\ell_3(u)}{u^{-1}}du\to1$ as $s\to\infty$, we see that $f_2=f/\ell_2\in \sR^{\infty}_{0}$.} Using these observation, we have that $\Phi$ satisfies $L^1(\delta, \wt C_L)$ for some $1<\delta<\gamma_4+1$ and $U^1(\gamma_4+1, c_U)$. Choose small $\theta_2>0$ such that $\frac{1}{\delta}+\theta_2(\frac{1}{\delta}-\frac{1}{\beta_2})=:\eps_2\in((1+\gamma_4)^{-1},1)$, where $\beta_2>0$ is the global upper scaling index of $\psi$. Then, by Lemma \ref{lem:inverse},  $L^1(\delta, \wt C_L, \Phi)$, $U^1(\gamma_4+1, c_U, \Phi)$ and $U(\beta_2, C_U, \psi)$, we see that there exist $c_1, c_2>0$ such that  for $2c_F^2\Phi^{-1}(t) <r \le  2c_F^2\frac{\Phi^{-1}(t)^{1+{\theta_2}}}{\psi^{-1}(t)^{\theta_2}}$ and $t\ge1$,
\begin{align}\label{range}
c_1 t^{1/(\gamma_4+1)}<r \le  c_2 t^{\eps_2}.
\end{align}
As we have seen in the proof of Corollary \ref{c:example}(ii), for the estimates of $(F'/f)^{-1}(t/r)$, it suffices to consider the case of $2c_F^2\Phi^{-1}(t)<r\le 2c_F^2\frac{\Phi^{-1}(t)^{1+{\theta_2}}}{\psi^{-1}(t)^{\theta_2}}$ and $t\ge1$ since $\wt\Phi(s)$ defined in \eqref{d:wePhi} is equal to $\Phi(s)$ for $s>1$. In this case, we have $t/r\ge c_2^{-1}t^{1-\eps_2}\ge c_2^{-1}$ by \eqref{range}.

By Remark \ref{r:ex} and the condition that $f_2$ satisfies \eqref{ellcheck} for $\eta\in\{ 1/\gamma_4,1/(\gamma_4-1)\}$, we have for $s\ge 1$, 
$$(F/f)^{-1}(s)\asymp(s^{\gamma_4}/f_2)^{-1}(s)\asymp s^{1/\gamma_4}(1/f_2^{1/{\gamma_4}})_{\#}(s^{1/{\gamma_4}})\asymp s^{1/{\gamma_4}}f_2^{1/{\gamma_4}}(s^{1/{\gamma_4}})$$ 
and 
$$(F'/f)^{-1}(s)\asymp(s^{\gamma_4-1}/f_2)^{-1}(s)\asymp s^{1/{(\gamma_4-1)}}(1/f_2^{1/{(\gamma_4-1)}})_{\#}(s^{1/{(\gamma_4-1)}})\asymp s^{1/({\gamma_4}-1)} f_2^{1/({\gamma_4}-1)}(s^{1/({\gamma_4}-1)}).$$
Using this, volume doubling property and the fact that $t/r\ge c_2^{-1}$ for  $2c_F^2\Phi^{-1}(t)<r\le 2c_F^2\frac{\Phi^{-1}(t)^{1+{\theta_2}}}{\psi^{-1}(t)^{\theta_2}}$ and $t\ge1$, we can apply Corollary \ref{c:example}(ii) to obtain the first claim.

For \eqref{e:larget}, it is enough to show that $f_2((t/r)^{1/(\gamma_4-1)})\asymp f_2(t)$ for $2c_F^2\Phi^{-1}(t)<r\le 2c_F^2\frac{\Phi^{-1}(t)^{1+{\theta_2}}}{\psi^{-1}(t)^{\theta_2}}$ and $t\ge1$. Using that $f_2$ is monotone and \eqref{range},  we see that 
$$f_2(c_2't^{(1-\eps_2)/(\gamma_4-1)})\le f_2((t/r)^{1/(\gamma_4-1)})\le f_2(c_1't^{(1-(\gamma_4+1)^{-1})/(\gamma_4-1)})$$ 
or 
$$f_2(c_2't^{(1-\eps_2)/(\gamma_4-1)})\ge f_2((t/r)^{1/(\gamma_4-1)})\ge f_2(c_1't^{(1-(\gamma_4+1)^{-1})/(\gamma_4-1)}).$$  
By Lemma \ref{comp_mon} and the fact that $f_2$ is slowly varying at $\infty$, we have $f_2(c_3 s^{c_4})\asymp f_2(s)$ for $s>1$. Thus, the above inequalities give that $f_2((t/r)^{1/(\gamma_4-1)})\asymp f_2(t)$ for $t\ge1$.
\qed

\begin{example}\label{ex:log1} 
{\rm Suppose that $F$ is differentiable function satisfying $F(s)\asymp sF'(s)$ and $F(s)\1_{\{s<1\}}\asymp s^{\gamma}(\log\tfrac1s)^{\kappa}\1_{\{s<1\}}$ for $\gamma>1$ and $\kappa\in\R$. 
Suppose further that  $\psi:(0,\infty)\to(0,\infty)$ is a
non-decreasing function which satisfies $L(\beta_1,C_L)$ and $U(\beta_2,C_U)$.
Define
$f_{\alpha,\beta}(s):=(\log \frac1s)^{1-\alpha}(\log\log \frac1s)^{-\beta}$
and $D:=\{(a,b)\in\R^2: a>1, b\in\bR\}\cup \{(1,b)\in\R^2: b>1\}$. Then, we observe that for $(\alpha, \beta)\in D$,
$\ell_{\alpha,\beta}(s):=sf_{\alpha,\beta}'(s)\asymp 
(\log \frac1s)^{-\alpha}(\log\log \frac1s)^{-\beta}$.
In particular,  $\ell_{\alpha,\beta}\in\sR^{0}_{0}$ and $f_{\alpha,\beta}(s)=\int^{s}_{0}\ell_{\alpha,\beta}(u)u^{-1}du$. Assume that  for $(\alpha, \beta)\in D$
\begin{align*}
\psi(\lambda) \asymp 
F(\lambda)\Big(\log \frac{1}{\lambda}\Big)^{-\alpha}\Big(\log\log \frac{1}{\lambda}\Big)^{-\beta}, \quad 0<\lambda<2^{-4}.
 \end{align*}
Then, $f_{\alpha-\kappa,\beta}$ satisfies \eqref{ellcheck} for all $\eta\in\bR\setminus\{0\}$. Moreover, there exist $T=T(\alpha-\kappa, \beta)\le 2^{-4}$ such that for $s\leq T$, $f_{\alpha-\kappa,\beta}$ is monotone and satisfies $(f_{\alpha-\kappa,\beta})(s^{\gamma})\asymp (f_{\alpha-\kappa,\beta})(s)$. Thus, by the above observation and \eqref{e:smallt}, we have the following  heat kernel estimates for $t<T$:
\begin{align*}
&p(t,x,y)\asymp \;\frac{1}{V(x, t^{1/\gamma}(f_{\alpha-\kappa,\beta})(t)^{1/\gamma})}\\
&\wedge
\left(\frac{t}{V(x,d(x,y))\psi(d(x,y))} + \frac{1}{V(x, t^{1/\gamma}(f_{\alpha-\kappa,\beta})(t)^{1/\gamma})}\exp\left({-a_1\left(\tfrac{d(x,y)^\gamma}{(f_{\alpha-\kappa,\beta})(t)}\right)^{1/(\gamma-1)}}\right)\right).
\end{align*}

}
\end{example}
 
\begin{example}\label{ex:log2} 
{\rm Suppose that $F$ is differentiable function satisfying $F(s)\asymp sF'(s)$ and $F(s)\1_{\{s>2\}}\asymp s^{\gamma'}(\log s)^{\kappa}\1_{\{s>2\}}$ for $\gamma'>1$ and $\kappa\in\R$. Suppose further that  $\psi:(0,\infty)\to(0,\infty)$ is a
non-decreasing function which satisfies \eqref{e:intcon}, $L(\beta_1,C_L)$, $U(\beta_2,C_U)$ and
$\psi(r){\1_{\{r>16\}}}\asymp F(r)(\log r)^\beta{\1_{\{r>16\}}}$ for $\beta\in \R$. 
Let $\ell(s)=(\log s)^{-\beta}$. Then for $\beta\leq 1$, 
$\int^{\infty}_{16}\frac{\ell(s)}{s}ds=\infty$. For  $s>16$,   let 
\begin{align*}
f(s)=
\begin{cases}
\frac{1}{1-\beta}(\log s)^{1-\beta} \quad & \text{if} \quad \beta<1,\\
\log\log s \quad & \text{if} \quad \beta=1.
\end{cases}
\end{align*}
Then, $f\in \Pi_{\ell}^\infty$ and $f(s)/(\log s)^{\kappa}$ satisfies \eqref{ellcheck} for all $\eta \in \R\setminus\{0\}$. Moreover, there exists $T=T(\beta, \kappa)\ge 16$ such that for $s\ge T$, $f(s)/(\log s)^{\kappa}$ is monotone and $f(s)/(\log s)^{\kappa}\asymp f(s^{\gamma'})/(\log s^{\gamma'})^{\kappa}$. 
Thus, by  \eqref{e:larget}, we have the following heat kernel estimates  for $t\ge T$:

\vspace{3mm}

\noindent(i) If $\beta<1$:
\begin{align*}
p(t,x,y)&\asymp
\;\frac{1}{V(x, t^{1/{\gamma'}}(\log t)^{{(1-\beta-\kappa)}/{\gamma'}})}\wedge \Biggl(\frac{t}{V(x,d(x,y))d(x,y)^{\gamma'}(\log(1+d(x,y)))^{\beta+\kappa}}\\
&\qquad\qquad+ \frac{1}{V(x, t^{1/{\gamma'}}(\log t)^{{(1-\beta-\kappa)}/{\gamma'}})} \exp\left(-a_2\left(\tfrac{d(x,y)^{\gamma'}}{t(\log t)^{1-\beta-\kappa}}\right)^{\frac{1}{\gamma'-1}} \right)\Biggl),
\end{align*}

\vspace{2mm}

\noindent(ii) If $\beta=1$:}
\begin{align*}
p(t,x,y) &\asymp \frac{1}{V(x, t^{1/{\gamma'}}(\log t)^{-\kappa/{\gamma'}}(\log \log t)^{1/{\gamma'}})} \wedge \Biggl(\frac{t}{V(x,d(x,y))d(x,y)^{\gamma'}(\log(1+d(x,y)))^{1+\kappa}}\\
&\qquad\qquad+ \frac{1}{V(x, t^{1/{\gamma'}}(\log t)^{-\kappa/{\gamma'}}(\log \log t)^{1/{\gamma'}})} \exp\left(-a_3\left(\tfrac{d(x,y)^{\gamma'}}{t(\log\log t)(\log t)^{-\kappa}}\right)^{\frac{1}{\gamma'-1}} \right)\Biggl).
\end{align*}
\end{example}

\begin{example}
{\rm Recall that $\gamma_1, \gamma_2>1$ are the constants in \eqref{wsF}. Suppose $F$ is differentiable function such that there exists $c>0$ satisfying  $\gamma_1 F(s)\le sF'(s)\le c F(s)$ for all $s>0$.
Let $T>0$ and  $\psi(r)=r^{\alpha}\1_{\{r\le1\}}+r^{\beta}\1_{\{r>1\}}$, where $\alpha<\gamma_1\le \gamma_2< \beta$. Then, by Corollary \ref{t:main32}, we see that for $t\le T$,
\begin{align}\label{ex21}
p(t,x,y) \asymp \frac{1}{V(x, t^{1/\alpha})}
\wedge \frac{ t}{V(x,d(x,y))\psi(d(x,y))} . 
\end{align}
Indeed, for $d(x,y)<1$, \eqref{ex21} follows from Theorem \ref{t:UHK-Phi}. If 
$d(x,y)\ge1$, then $\frac{ t}{V(x,d(x,y))\psi(d(x,y))}$ dominates the upper bound of off-diagonal term in \eqref{GHKjum}.

On the other hand, by the condition $\gamma_2<\beta$, we have $\int^{\infty}_{0}\frac{dF(s)}{\psi(s)}\le c+c\int^{\infty}_{1}\frac{s^{\gamma_2}}{s^{1+\beta}}ds<\infty$. Thus, for $r>1$, $\Phi(r)$ defined in \eqref{e:dPhi} is comparable to $F(r)$ and $\Phi(r)/r\asymp F(r)/r\asymp F'(r)$. Now, by the same argument as in the proof of Corollary \ref{c:example}(ii), we see that for $t>T$,
\begin{align*}
p(t,x,y) \asymp \frac{1}{V(x, F^{-1}(t))}
\wedge \left(\frac{ t}{V(x,d(x,y))d(x,y)^{\beta}} + \frac{1}{V(x, F^{-1}(t))}
\exp{\left(-a_5\frac{d(x,y)}{F'(t/d(x,y))}\right)}
 \right). 
\end{align*}
}\end{example}

\appendix

\section{Appendix}\label{s:App}

We first observe simple consequences of weak scaling conditions. 

\begin{remark}\label{mwsc}
	\rm Suppose $g:(0,\infty)\to(0,\infty)$ is non-decreasing.
If $g$ satisfies $L_a(\beta,  c)$, then $g$ satisfies $L_b(\beta,  c(ab^{-1})^{\beta})$ for any $b > a$. Indeed, for $r\leq a \leq R\leq b$, 
$$
g(R)\geq g(a)\geq c\left(\frac{a}{r}\right)^{\beta}g(r)\geq c\left(\frac{a}{b}\right)^{\beta}\Big(\frac{R}{r}\Big)^{\beta}g(r)
$$
and for $a \leq r\leq R\leq b$,
$$
g(R)\geq g(r)\geq c\left(\frac{a}{b}\right)^{\beta}\Big(\frac{R}{r}\Big)^{\beta}g(r).
$$
Similarly, if $g$ satisfies $L^a(\beta, c)$, then $g$ satisfies $L^b(\beta,  c(a^{-1}b)^{\beta})$ for $b<a$.
\end{remark}

\begin{lemma}\label{lem:inverse}
Let $g:(0,\infty) \rightarrow (0,\infty)$ be a 
non-decreasing function with $g(\infty)=\infty.$
\begin{enumerate}
\item[(1)] If $g$ satisfies $L_a(\beta, c)$ $($resp. $U_a(\beta, C))$,  then $g^{-1}$ satisfies $U_{g(a)}(1/\beta,c^{-1/\beta})$\\  $($resp. $L_{g(a)}(1/\beta,C^{-1/\beta}))$. 
\item[(2)] If $g$ satisfies $L^a(\beta, c)$  (resp. $U^a(\beta, C)$),  then $g^{-1}$ satisfies $U^{g(a)}(1/\beta,c^{-1/\beta})$\\ $($resp. $L^{g(a)}(1/\beta,C^{-1/\beta}))$.
\end{enumerate}
\end{lemma}

We now give the proof of Remark \ref{r:hk}.

\vspace{3mm}
\noindent \textit{Proof of Remark \ref{r:hk}.} Assume that the condition $\HK(\Phi,C\Phi)$ holds. If $d(x,y) \le \eta \Phi^{-1}(t)$, by $\VD(d_2)$ and $L(\alpha_1,c_L,\Phi)$ we have 
\begin{align}\begin{split}\label{e:rhk1-1}
\frac{1}{V(x,\Phi^{-1}(t))} &\le \frac{c_1}{V(x,\eta \Phi^{-1}(t))} \le \frac{c_1}{V(x,d(x,y))} \\ &\le \frac{c_1t}{V(x,d(x,y)) \Phi(\eta^{-1} d(x,y))} \le \frac{c_2 t}{V(x,d(x,y)) \Phi( d(x,y))} . 
\end{split}\end{align}
Also, if $d(x,y) \ge \eta \Phi^{-1}(t)$, using $U(\alpha_2,c_U,\Phi)$, we obtain that
\begin{align*}
\frac{1}{V(x,\Phi^{-1}(t))} &\ge \frac{c_3}{V(x,\eta \Phi^{-1}(t))} \ge \frac{c_3}{V(x,d(x,y))} \\ &\ge \frac{c_3t}{V(x,d(x,y)) \Phi(\eta^{-1} d(x,y))} \ge  \frac{c_4 t}{V(x,d(x,y)) \Phi( d(x,y))} .  
\end{align*}
By above two inequalities, we have that lower bounds in \eqref{HKjum} and \eqref{e:rhk1} are equivalent. 

For the upper bound, it suffices to verify the existence of constant $c>0$ satisfying
$$ c^{-1} \frac{t}{V(x,d(x,y))\Phi(d(x,y))} \le  \sG(a_0,t,x,d(x,y)) 
\le c\frac{t}{V(x,d(x,y))\Phi(d(x,y))} $$
for all $t>0$ and $x,y \in M$ with $d(x,y) \ge 2c_U\Phi^{-1}(t)$. Indeed, when $d(x,y) \le 2c_U\Phi^{-1}(t)$, following the calculations in \eqref{e:rhk1-1} we have 
$$ \frac{c_5}{V(x,\Phi^{-1}(t))} \le \frac{t}{2V(x,d(x,y))\Phi(d(x,y))} \le \sG(a_0,t,x,d(x,y)). $$
Note that the second inequality immediately follows from the definition of $\sG$. Now we observe that from
$ \Phi_1(r,t) = \sup_{s>0} \big[ \frac{r}{s} - \frac{t}{\Phi(s)} \big] $ and $r \ge 2c_U\Phi^{-1}(t)$,
$$ \Phi_1(r,t) \ge \frac{r}{\Phi^{-1}(t)} - \frac{t}{\Phi(\Phi^{-1}(t))} \ge \frac{r}{2\Phi^{-1}(t)}, $$
where we have used \eqref{e:ginv} for the second inequality. Thus,
\begin{align*}
&\frac{1}{V(x,\Phi^{-1}(t))} \exp(-a_0\Phi_1(d(x,y), t)) \le  \frac{c_6}{V(x,\Phi^{-1}(t))} \exp\Big(-\frac{a_0}{2} \frac{d(x,y)}{\Phi^{-1}(t)}\Big) \\
&\le \frac{c_7}{V(x,\Phi^{-1}(t))} \Big( \frac{\Phi^{-1}(t)}{d(x,y)} \Big)^{d_2+\beta_2} 
 \le \frac{c_8}{V(x,\Phi^{-1}(t))} \frac{V(x,\Phi^{-1}(t))}{V(x,d(x,y))} \frac{\Phi(\Phi^{-1}(t))}{\Phi(d(x,y))} = \frac{c_8 t}{V(x,d(x,y))\Phi(d(x,y))}, 
\end{align*}
where we have used $\VD(d_2)$ and $U(\beta_2,C_U,\Phi)$ in the last inequality. This finishes the proof. \qed

Suppose that the function $\Phi$ satisfies $L_a(\delta, \wt C_L)$ for some $a \in (0,\infty]$, $\wt C_L>0$ and $\delta>1$. Now, define
$$ \scK(r):= \sup_{0<s \le r} \frac{\Phi(s)}{s}. $$
Then, for any $R_0>0$, letting $c_1=\wt C_L^{-1} (aR_0^{-1}\wedge 1)^{-\delta}\ge1$ we have for any $r \in (0, R_0]$,
\begin{equation}\label{e:scK1}
\frac{\Phi(r)}{r} \le \scK(r) \le  c_1 \frac{\Phi(r)}{r},
\end{equation}
 (c.f. \cite[Lemma 2.5]{BKKL}). Note that if $a =\infty$, \eqref{e:scK1} holds for every $r>0$.  Let $\scK^{-1}$ be the generalized inverse of non-decreasing function $\scK$.
 
  The following lemma yields that \cite[Theorem 1.4]{BKKL} is the special case of Corollary \ref{c:main32}.

\begin{lemma}\label{l:scK} 
Suppose $\Phi$ is non-decreasing function satisfying $L(\alpha_1, c_L)$, $U(\alpha_2, c_U)$ and $L_{a}(\delta, \wt C_L)$ for $\delta>1$. Let $T \in (0,\infty)$. Then, there exists a constant $c>1$ such that for any $t \in (0, T]$ and $r \ge 2c_U^2\Phi^{-1}(t)$, 
	\begin{equation}\label{e:scK}
c^{-1}\Phi_1(r,t) \le 	\frac{r}{\scK^{-1}(t/r)} \le c \Phi_1(r,t).
	\end{equation}
Moreover, if $L(\delta, \wt C_L)$ holds, then \eqref{e:scK} holds for any $t \in (0,\infty)$ and $r \ge  2c_U^2 \Phi^{-1}(t)$.
\end{lemma} 
\pf   
 Without loss of generality we may and do assume $a = \Phi^{-1}(T)$.     Note that $\alpha_2>1$. Let $R_0:=\Phi^{-1}(T)$ and $c_1= \wt C_L^{-1}$ so that $L_{R_0}(\delta, c_1^{-1}, \Phi)$ and \eqref{e:scK1} hold. Denote $\eps:= \frac{1}{\alpha_2 - 1}$. 
Since $r \ge 2c_U^2 \Phi^{-1}(t)$, we have
$$c_U^{2\eps} \frac{\Phi^{-1}(t)^{1+\eps}}{r^{\eps}}\le \Phi^{-1}(t) \le R_0.$$
It follows from \eqref{e:scK1}, Lemma \ref{l:ginv} and $U(\alpha_2,c_U,\Phi)$ that
\begin{align*}
\scK\left(c_U^{2\eps} \frac{\Phi^{-1}(t)^{1+\eps}}{r^{\eps}}\right)&\ge c_U^{-2\eps} \frac{r^{\eps}}{\Phi^{-1}(t)^{1+\eps}}\Phi\left(\Phi^{-1}(t)c_U^{2\eps}\frac{\Phi^{-1}(t)^{\eps}}{r^{\eps}}\right) \ge c_U^{-1-2\eps}\frac{r^{\eps}t}{\Phi^{-1}(t)^{1+\eps}}\frac{\Phi(\Phi^{-1}(t)c_U^{2\eps}\frac{\Phi^{-1}(t)^{\eps}}{r^{\eps}})}{\Phi(\Phi^{-1}(t))}\\
&\ge c_U^{-2-2\eps}\frac{t}{r}\frac{r^{1+\eps}}{\Phi^{-1}(t)^{1+\eps}}\left(c_U^{2\eps}\frac{\Phi^{-1}(t)^{\eps}}{r^{\eps}}\right)^{\alpha_2}
= \frac{t}{r}.
\end{align*}
Thus,
\begin{align*}
\rho:=\scK^{-1}\Big(\frac{t}{r}\Big)\le c_U^{2\eps}\frac{\Phi^{-1}(t)^{1+\eps}}{r^{\eps}}\le 2^{-\eps}\Phi^{-1}(t) \le R_0.
\end{align*}
By \eqref{e:scK1},  $\scK$ satisfies $U_{R_0}(\alpha_2-1, c_1c_U)$ and $L_{R_0}(\delta-1, c_1^{-1}\wt C_L)$. Thus, using Lemma \ref{l:ginv} we have
\begin{align}\label{rhocomp}
(c_1c_U)^{-1}\frac{t}{r}\le \scK(\rho)\le c_1c_U \frac{t}{r}.
\end{align} 
Using \eqref{rhocomp} and \eqref{e:scK1}, we have
\begin{align*}
(c_1c_U)^{-1}\frac{t}{r}\le \scK(\rho)\le c_1\frac{\Phi(\rho)}{\rho}.
\end{align*}
 Then, letting $c_2= c_1^2 c_U$, the above inequality and \eqref{e:sT1} imply that there exists $c_3>0$ such that
\begin{align*}
c_3 \Phi_1(r,t) \ge \Phi_1(2c_2 r,t) \ge \frac{2c_2r}{\rho} - \frac{t}{\Phi(\rho)} \ge \frac{r}{\rho} \Big(2c_2 - \frac{t}{r} \frac{\rho}{\Phi(\rho)} \Big) \ge  \frac{c_2 r}{\rho} = \frac{c_2 r}{\scK^{-1}(t/r)}.
\end{align*}
 This proves the second inequality in \eqref{e:scK}. For the first one, we take a $s>0$ such that 
\begin{align}\label{sckeq00}
0 \le \frac{r}{s} - \frac{t}{\Phi(s)} \le \Phi_1(r,t) \le 2\Big( \frac{r}{s} - \frac{t}{\Phi(s)} \Big).
\end{align}
Since $\Phi_1(r,t) \ge 0$, we have $\Phi(s)/s \ge t/r$. Using this, \eqref{e:scK1} and \eqref{rhocomp} we have
\begin{align}\label{sckeq0}  \frac{\Phi(\rho)}{\rho} \le \scK(\rho)\le c_1c_U \frac{t}{r}\le c_1 c_U \frac{\Phi(s)}{s}.\end{align} 
Thus, if $s<\rho \le R_0$, using $L_{R_0}(\delta, c_1^{-1}, \Phi)$ and \eqref{sckeq0}
$$c_1^{-1}\Big(\frac{\rho}{s}\Big)^{\delta-1}\le \frac{\Phi(\rho)}{\rho}\Big/\frac{\Phi(s)}{s}\le c_1 c_U.$$
Thus, we conclude that there is $c_4>0$ such that $s>c_4\rho$. Using this and \eqref{sckeq00}, we have 
\begin{align*} 
\Phi_1(r,t) \le 2 \frac{r}{s} \le 2 c_4^{-1}\frac{r}{\rho} = 2 c_4^{-1} \frac{r}{\scK^{-1}(t/r)}.
\end{align*}
When $L(\delta, \wt C_L,\Phi)$ holds, we may take $R_0=\infty$ and $c_1=\wt C_L^{-1}$. Then, the proof is the same with the above since \eqref{e:scK1} holds for all $r>0$ and \eqref{e:sT1} holds for all $t>0$ and $r \ge 2c_U^2\Phi^{-1}(t)$. This completes the proof. \qed

\begin{lemma}\label{comp_mon}
Let $h:(0,\infty)\to (0,\infty)$ be a monotone function. Suppose that there exist $k>1$ and $c_1>1$ such that $c_1^{-1}h(r)\le h(r^k)\le c_1h(r)$ for all  $r<1$ ({resp.} $r>1$). Then, for any $m>0$, there exists $c_2=c_2(k, c_1, m)>1$ such that  
$$c_2^{-1}h(r)\le h(r^m)\le c_2h(r)\quad\text{for all}\;\;r<1\,\, (\text{resp.}\; r>1).$$  
\end{lemma}

\pf Since proofs for other cases are all similar, we only prove the case that $h$ is non-decreasing and the comparability condition holds for $r<1$. 
When $1\le m\le k$, using  the comparability condition, clearly $h(r^m)\le h(r)\le c_1h(r^k)\le c_1h(r^m)$.

If $m<1$, we let $n\in\bN$ be a constant satisfying $mk^n\ge1$. Then, we get $h(r)\le h(r^m)\le c_1h(r^{mk})\le c_1^2h(r^{mk^2})\le\dots\le c_1^{n}h(r^{mk^n})\le c_1^{n}h(r).$

For $m>k$, let $l\in\bN$ be a constant satisfying $k^l\ge m$. Then, $h(r^m)\le h(r)\le c_1h(r^{k})\le c_1^2h(r^{k^2})\le\dots\le c_1^{l}h(r^{k^l})\le c_1^{l}h(r^m).$ 
\qed

\noindent
{\bf Acknowledgements:}

\smallskip 
After finishing the main part of this paper, the result was announced
by the fourth  named author during 7th Bielefeld-SNU Joint Workshop in Mathematics at Bielefeld University on February 27, 2019.
 We thank  Center for Interdisciplinary Research and the Department of Mathematics of Bielefeld University for the hospitality.

\end{doublespace}
\begin{singlespace}
	
\end{singlespace}

\end{document}